\documentclass[10pt]{amsart}
\usepackage{amsmath}
\usepackage{amsxtra}
\usepackage{amscd}
\usepackage{amsthm}
\usepackage{amsfonts}
\usepackage{amssymb}
\theoremstyle{definition}
\theoremstyle{remark}

\newcommand{\thmref}[1]{Theorem~\ref{#1}}
\newcommand{\secref}[1]{\S \ref{#1}}
\newcommand{\lemref}[1]{Lemma~\ref{#1}}
\newcommand{\defref}[1]{Definition~\ref{#1}}
\newcommand{\propref}[1]{Proposition~\ref{#1}}
\newcommand{\corref}[1]{Corollary~\ref{#1}}

\newcommand{\remref}[1]{Remark~\ref{#1}}
\newcommand{\exref}[1]{Example~\ref{#1}}

\newcommand{\nc}{\newcommand}
\nc{\renc}{\renewcommand}
\nc{\ssec}{\subsection}
\nc{\sssec}{\subsubsection}
\nc{\on}{\operatorname}

\nc\ol{\overline}
\nc\wt{\widetilde}
\nc\wh{\widehat}
\nc\tboxtimes{\wt{\boxtimes}}

\renc{\d}{{\delta}}
\nc{\Aa}{{\mathbb{A}}}
 \nc{\Gg}{{\mathbb{G}}}  
\nc{\Hh}{{\mathbb{H}}}
 \nc{\Nn}{{\mathbb{N}}}
\nc{\Pp}{{\mathbb{P}}}
\nc{\Rr}{{\mathbb{R}}}
\nc{\BV}{{\mathbb{V}}}
\nc{\BW}{{\mathbb{W}}}
\nc{\Zz}{{\mathbb{Z}}}
\nc{\Qq}{{\mathbb{Q}}}
\nc{\Ss}{{\mathbb{S}}}
\nc{\Cc}{{\mathbb{C}}}
\nc{\Ff}{{\mathbb{F}}}

\nc{\CA}{{\mathcal{A}}}
\nc{\CB}{{\mathcal{B}}}

\nc{\CE}{{\mathcal{E}}}
\nc{\CF}{{\mathcal{F}}}
\nc{\CG}{{\mathcal{G}}}
\nc{\CL}{{\mathcal{L}}}
\nc{\CC}{{\mathcal{C}}}
\nc{\CM}{{\mathcal{M}}}
\def\Mm{\CM}
\nc{\CN}{{\mathcal{N}}}
\nc{\Oo}{{\mathcal{O}}}
\nc{\CP}{{\mathcal{P}}}
\nc{\CQ}{{\mathcal{Q}}}
\nc{\CR}{{\mathcal{R}}}
\nc{\CS}{{\mathcal{S}}}
\nc{\CT}{{\mathcal{T}}}
\nc{\CU}{{\mathcal{U}}}
\nc{\CV}{{\mathcal{V}}}
\nc{\CK}{{\mathcal{K}}}
\nc{\CW}{{\mathcal{W}}}
\nc{\CZ}{{\mathcal{Z}}}

\nc{\cM}{{\check{\mathcal M}}{}}
\nc{\csM}{{\check{\mathcal A}}{}}
\nc{\oM}{{\overset{\circ}{\mathcal M}}{}}
\nc{\obM}{{\overset{\circ}{\mathbf M}}{}}
\nc{\oCA}{{\overset{\circ}{\mathcal A}}{}}
\nc{\obA}{{\overset{\circ}{\mathbf A}}{}}
\nc{\ooM}{{\overset{\circ}{M}}{}}
\nc{\osM}{{\overset{\circ}{\mathsf M}}{}}
\nc{\vM}{{\overset{\bullet}{\mathcal M}}{}}
\nc{\nM}{{\underset{\bullet}{\mathcal M}}{}}
\nc{\oD}{{\overset{\circ}{\mathcal D}}{}}
\nc{\obD}{{\overset{\circ}{\mathbf D}}{}}
\nc{\oA}{{\overset{\circ}{\mathbb A}}{}}
\nc{\op}{{\overset{\bullet}{\mathbf p}}{}}
\nc{\cp}{{\overset{\circ}{\mathbf p}}{}}
\nc{\oU}{{\overset{\bullet}{\mathcal U}}{}}
\nc{\oZ}{{\overset{\circ}{\mathcal Z}}{}}
\nc{\ofZ}{{\overset{\circ}{\mathfrak Z}}{}}
\nc{\oF}{{\overset{\circ}{\fF}}}

\nc{\fa}{{\mathfrak{a}}}
\nc{\fb}{{\mathfrak{b}}}
\nc{\fg}{{\mathfrak{g}}}
\nc{\fgl}{{\mathfrak{gl}}}
\nc{\fh}{{\mathfrak{h}}}
\nc{\fj}{{\mathfrak{j}}}
\nc{\fm}{{\mathfrak{m}}}
\nc{\fn}{{\mathfrak{n}}}
\nc{\fu}{{\mathfrak{u}}}
\nc{\fp}{{\mathfrak{p}}}
\nc{\fr}{{\mathfrak{r}}}
\nc{\fs}{{\mathfrak{s}}}
\nc{\fsl}{{\mathfrak{sl}}}
\nc{\hsl}{{\widehat{\mathfrak{sl}}}}
\nc{\hgl}{{\widehat{\mathfrak{gl}}}}
\nc{\hg}{{\widehat{\mathfrak{g}}}}
\nc{\chg}{{\widehat{\mathfrak{g}}}{}^\vee}
\nc{\hn}{{\widehat{\mathfrak{n}}}}
\nc{\chn}{{\widehat{\mathfrak{n}}}{}^\vee}

\nc{\fA}{{\mathfrak{A}}}
\nc{\fB}{{\mathfrak{B}}}
\nc{\fD}{{\mathfrak{D}}}
\nc{\fE}{{\mathfrak{E}}}
\nc{\fF}{{\mathfrak{F}}}
\nc{\fG}{{\mathfrak{G}}}
\nc{\fK}{{\mathfrak{K}}}
\nc{\fL}{{\mathfrak{L}}}
\nc{\fM}{{\mathfrak{M}}}
\nc{\fN}{{\mathfrak{N}}}
\nc{\fP}{{\mathfrak{P}}}
\nc{\fU}{{\mathfrak{U}}}
\nc{\fV}{{\mathfrak{V}}}
\nc{\fZ}{{\mathfrak{Z}}}

\nc{\bb}{{\mathbf{b}}}
\nc{\bc}{{\mathbf{c}}}
\nc{\bd}{{\mathbf{d}}}
\nc{\be}{{\mathbf{e}}}
\nc{\bj}{{\mathbf{j}}}
\nc{\bn}{{\mathbf{n}}}
\nc{\bp}{{\mathbf{p}}}
\nc{\bq}{{\mathbf{q}}}
\nc{\bF}{{\mathbf{F}}}
\nc{\bu}{{\mathbf{u}}}
\nc{\bv}{{\mathbf{v}}}
\nc{\bx}{{\mathbf{x}}}
\nc{\bs}{{\mathbf{s}}}
\nc{\by}{{\mathbf{y}}}
\nc{\bw}{{\mathbf{w}}}
\nc{\bA}{{\mathbf{A}}}
\nc{\bK}{{\mathbf{K}}}
\nc{\bI}{{\mathbf{I}}}
\nc{\bB}{{\mathbf{B}}}

\def\mbG{{\mathbf \Gamma}}
\nc{\bG}{{\mathbf{G}}}
\nc{\bC}{{\mathbf{C}}}
\nc{\bD}{{\mathbf{D}}}
\nc{\bP}{{\mathbf{P}}}
\nc{\bH}{{\mathbf{H}}}
\nc{\bM}{{\mathbf{M}}}
\nc{\bN}{{\mathbf{N}}}
\nc{\bV}{{\mathbf{V}}}
\nc{\bU}{{\mathbf{U}}}
\nc{\bL}{{\mathbf{L}}}
\nc{\bT}{{\mathbf{T}}}
\nc{\bW}{{\mathbf{W}}}
\nc{\bX}{{\mathbf{X}}}
\nc{\bY}{{\mathbf{Y}}}
\nc{\bZ}{{\mathbf{Z}}}
\nc{\bS}{{\mathbf{S}}}

\nc{\sA}{{\mathsf{A}}}
\nc{\sB}{{\mathsf{B}}}
\nc{\sC}{{\mathsf{C}}}
\nc{\sD}{{\mathsf{D}}}
\nc{\sF}{{\mathsf{F}}}
\nc{\sG}{{\mathsf{G}}}
\nc{\sK}{{\mathsf{K}}}
\nc{\sM}{{\mathsf{M}}}
\nc{\sO}{{\mathsf{O}}}
\nc{\sQ}{{\mathsf{Q}}}
\nc{\sP}{{\mathsf{P}}}
\nc{\sZ}{{\mathsf{Z}}}
\nc{\sfp}{{\mathsf{p}}}
\nc{\sr}{{\mathsf{r}}}
\nc{\sg}{{\mathsf{g}}}
\nc{\sff}{{\mathsf{f}}}
\nc{\sfb}{{\mathsf{b}}}
\nc{\sfc}{{\mathsf{c}}}
\nc{\sd}{{\ltimes}}

\nc{\tN}{{\widetilde{{N}}}}
\nc{\tA}{{\widetilde{{A}}}}
\nc{\tB}{{\widetilde{{B}}}}
\nc{\tg}{{\widetilde{\mathfrak{g}}}}
\nc{\tG}{{\widetilde{G}}}
\nc{\tH}{{\widetilde{H}}} 
\nc{\tO}{{\widetilde{\mathsf{O}}}{}} 
\nc{\tU}{\widetilde{U}}
\nc{\TZ}{{\tilde{Z}}}
\nc{\tx}{{\tilde{x}}}
\nc{\tq}{{\tilde{q}}}

\nc{\tfP}{{\widetilde{\mathfrak{P}}}{}}
\nc{\tz}{{\tilde{\zeta}}}
\nc{\tmu}{{\tilde{\mu}}}

 \def\e{\epsilon}

\def\e{\epsilon}

\def\bdr{\partial}
  \nc{\Ob}{{\mathop{\operatorname{\rm Ob}}}}
  \nc{\Sym}{{\mathop{\operatorname{\rm Sym}}}}
   \nc{\Aut}{{\mathop{\operatorname{\rm Aut}}}}
 \nc{\Spec}{{\mathop{\operatorname{\rm Spec}}}}
  \nc{\spec}{{\mathop{\operatorname{\rm Spec}}}}
\nc{\Ker}{{\mathop{\operatorname{\rm Ker}}}}
 \nc{\dom}{{\mathop{\operatorname{\rm dom}}}}
\nc{\End}{{\mathop{\operatorname{\rm End}}}}
 \nc{\Hom}{\on{\Hom}}
 \nc{\GL}{{\mathop{\operatorname{\rm GL}}}}
 \nc{\Id}{{\mathop{\operatorname{\rm Id}}}}
 \nc{\rk}{{\mathop{\operatorname{\rm rk}}}} 
 \nc{\length}{{\mathop{\operatorname{\rm length}}}}
\nc{\supp}{{\mathop{\operatorname{\rm supp}}}}
\nc{\val}{{\rm val}}
\nc{\res}{{\mathop{\operatorname{\rm res}}}}

\def\meet{\cap}
\def\union{\cup}

\def\si{\sigma}
\def\g{\gamma}
\def\G{\Gamma}

\def\<{\begin}
 \def\>{\end}

\def\m{\smallsetminus}

\nc{\seq}[1]{\stackrel{#1}{\sim}}

\def\inv{^{-1}}
\def\claim#1{{\noindent \bf Claim #1.\ }}

\def\beq#1{{\begin{equation} \label{#1}}  }
\def\normal{\trianglelefteq}
\def\Uu{\mathbb U}

\def\prf{\begin{proof}}
\def\pv{\end{proof} }
 \def\eprf{\end{proof} }

\def\liminv{\underset{\longleftarrow}{lim}\,}

 \def\lbl#1{\underline{#1}  \label{#1}  }
\def\lbl#1{  \label{#1}  }
\def\a{\alpha}

 \renc{\b}{{\beta}}

\begin{document}

  \def\tE{\widetilde{E}}  \def\ta{\tilde{a}}

\def\e{\epsilon}
 \def\vol{\rm vol}
 \def\ml{{\mathfrak l}}
 \def\ben{\begin{enumerate}} \def\enn{\end{enumerate}}  \def\een{\enn}
 \def\inft{\bigwedge}
 \def\mE{{\mathcal E}}
 
  \def\Ups{\Upsilon}
 \def\cU{\bar{U}}
 \def\mQ{{\mathcal Q}}
  \title{Stable group theory and approximate subgroups}
\author{Ehud Hrushovski  }
\address{Hebrew University at Jerusalem}
\email{ehud@math.huji.ac.il}
\thanks{ISF grant 1048/07}

 \<{abstract}  We note a parallel between some ideas of stable model theory and certain topics in  finite combinatorics  related to the sum-product phenomenon.    For a simple linear group $G$,  we show that a finite subset $X$ with $|X X \inv X |/ |X|$ bounded is close to a finite subgroup, or else to a subset of a proper algebraic subgroup of $G$.   We also find a connection with Lie groups, and use it to obtain some
 consequences suggestive of topological nilpotence.  Model-theoretically we prove the independence theorem and 
 the stabilizer theorem in a general first-order setting.   \>{abstract}
 
 \maketitle   

\<{section}{Introduction}

   Stable group theory, as developed in the 1970's and 80's, was an effective bridge between   definable sets  and objects of more geometric categories.   One of the reasons was a body of results showing that groups can be recognized from their traces in softer categories.  The first and simplest example is Zilber's stabilizer. Working with an integer-valued
 dimension theory on the definable subsets of a group $G$, Zilber considered
 the dimension-theoretic stabilizer of a definable set $X$:  this is the group $S$ of elements $g \in G$ with $gX \triangle X$ of smaller dimension than $X$.  Let $XX$ be the product set  $XX=\{xy: x, y \in X \}$.
 If $X$ differs little from $XX$ in the sense that $\dim(XX \triangle X) < \dim(X)$, Zilber showed that $X$ differs little from a coset of $S$.   
    
 In the 90's, Zilber's theory was generalized to the `simple theories" of \cite{shelah-simple},
 again initially in a definable finite dimensional context (\cite{CH}, \cite{PAC}).   
   Here the definable sets $X_t$ in a definable family   $(X_t: t  \in T)$ are viewed as   ``differing little from each other" if  simply the  pairwise intersections $X_t \meet X_{t'}$ have the same dimension as each $X_t$.   Nevertheless it is shown that when  the family of translates $(Xa: a \in X)$ satisfies this condition, there is a group $H$ of the same dimension as $X$ and with a large intersection with some translate of $X$; this 
 group was still, somewhat inappropriately, called the stabilizer, and we will keep this terminology.   
 
In the present paper we prove the stabilizer theorem in a general first-order setting.  A definition is
given of being a``near-subgroup" (\defref{near}), generalizing the stable and simple cases.   We 
then prove the existence of a nearby group (\thmref{stabilizer}.)    
In outline, the  proof remains the  same as in \cite{PAC}; the definability condition on the dimension was removed in \cite{kim-pillay}.   The key  is a general amalgamation statement for definable ternary relations, dubbed the ``Independence Theorem" (see \cite{CH}, p. 9 and p. 185).  Roughly speaking, 
 in maximal dimension, consistent relations among each pair of types determine consistent relations on a triple; see \thmref{ind}.  
 
The stabilizer obtained in \thmref{stabilizer} is  not a definable group but an  $\inft$-definable one; it is defined by a countable set of formulas
in a saturated model, or  alternatively as a group object in the category of   projective limits of definable sets.
In the finite dimensional setting of \cite{PAC} the construction of the stabilizer was complemented by a proof
that $\inft$-definable groups are limits of definable groups.    This last step is not true at the level of generality
considered here:  the group of infinitesimals of a Lie groups  provide  counterexamples.  We show however that all counterexamples are closely associated with Lie groups:  see \thmref{ap1}.  
The proof uses
the Gleason-Yamabe-Montgomery-Zippin structure theory for locally compact groups.

A very interesting dictionary between this part of model theory, and certain parts of finite combinatorics, 
can be obtained by making the model-theoretic
 ``dimension n" correspond  to the combinatorial    ``cardinality of order $c^n$" (cf. \cite{CH}, 8.4).  
    Near-subgroups in the above sense
 then correspond to  asymptotic families of finite subsets $X$ of a group (or a family of groups), with
 $(X \union X \inv)^3/|X|$ bounded.  Equivalently (see \cite{tao-nc}, Lemma 3.4, and \corref{3-to-n} below)   $|X^k|/|X|$ is bounded for any given $k$.    Subsets of groups with weak closure conditions  were   considered in combinatorics at least since
 \cite{erdos-sz}.  An excellent survey centering on rings can be found in the first pages of  \cite{tao-rings};
see also \cite{tao-nc} for more general non-commutative groups.      The parallels
to the model-theoretic development are striking.   
We turn now to a description of some consequences of the stabilizer theorem in this combinatorial setting.

For the sake of the introduction we consider finite subsets of $G$ (more general
situations will be allowed later.)  We recall Terence Tao's notion of an {\em approximate subgroup}.  A finite subset
$X \subseteq G$ is said to be a $k$-approximate group if $1 \in X, X=X \inv$,
and $XX$ is contained in $k$ right cosets of $X$. Say $X,Y$ are {\em commensurable}
if each is contained in finitely many right cosets of the other, with the number bounded in terms of $k$.    It is felt that approximate subgroups should be commensurable to actual subgroups, except in situations involving Abelian groups in some way.  See
 \cite{tao-blog} for a compelling exposition of the issue.

 Gromov's theorem \cite{gromov} on finitely generated
groups of polynomial growth fits into this framework, taking $X$ to be a ball of size $2^n$ in the Cayley graph,
for large $n$; then $X$ is  a $2^d$-approximate subgroup, where $d$ is the growth exponent.  Gromov shows
that the group is nilpotent, up to finite index.

\thmref{ap1} says nothing about a fixed finite approximate subgroup, but it does have asymptotic consequences
to the family of all $k$-approximate subgroups for fixed $k$.    In particular, we obtain:

 \<{thm}\lbl{gc0} Let $f: \Nn^2 \to \Nn$ be any function, and fix $k \in \Nn$.   Then there exist
$e^*,c^*,N \in \Nn$    such that the following holds.  

 Let $G$ be any  group, $X$ a finite subset, and assume  $|X X \inv X| \leq k |X|$.
 
 Then there are $e \leq e^*, c \leq c^*$, and  subsets 
$X_N \subseteq X_{N-1} \subseteq \cdots \subseteq X_1    \subseteq  X \inv X X \inv X$ such that
  $X,X_1$  are $e$-commensurable, 
  and for $1 \leq m,n < N$ we have:
\<{enumerate} 
\item  $X_{n}=X_{n} \inv$
\item   $X_{n+1}X_{n+1} \subseteq X_{n}$
\item $X_{n}$ is contained in the union of  $c$ translates of $X_{n+1}$.
\item  $[X_{n},X_{m}] \subseteq X_{k}$ whenever $k \leq N$ and $k < n+m$.
\item  $N> f(e,c)$.
\>{enumerate}
\>{thm}  

 Roughly speaking,  this is deduced as a special case of the following principle:  
if a sentence of a certain logic holds of all compact neighborhoods of the identity in all finite-dimensional Lie groups, then it holds of all approximate subgroups.  We have not explicitly determined the relevant logic; 
  \propref{sop1} hints that, given further work
on the first order theory of Lie groups with distinguished closed subsets, much stronger transfer principles may be possible than what we have used.  

  The first three clauses of \thmref{gc0} suggest a part of  a non-commutative Bourgain system as defined 
in \cite{gs}, and conjectured by   Ben Green in \cite{tao-blog} to exist for approximate subgroups.   Green's
conjecture was in part intended to show that ``one can do a kind of approximate representation
theory", which can be viewed as a description of \thmref{ap1} and the deduction between the two.

The fourth clause suggests  a kind of topological nilpotence.    Note that (4) implies that $[X_1,X_1] \subseteq X_1$.   For a set of generators of a finite simple group, this in itself seems to be a curious property.

The use of the structure theory of locally compact groups here follows Gromov  \cite{gromov}.    But the bridge to locally compact groups is a different one:  Gromov's is metric, while ours is measure-theoretic.

It is natural to consider a somewhat more general framework.   Call a pair $(X,G,\cdot,\inv,1)$ a {\em  Freiman approximate   group} if 
$X$ is a finite subset of $G$, $\cdot: X^{(2)} \to G$ and $\inv: X \to X$ are functions, such that 
for any  $(x_1,\cdots,x_{12}) \in X^{12}$, the iterated products $((x_1 \cdot x_2) \cdot (x_3 \cdot \ldots ))$ are 
defined and independent of the placing of the parentheses; $x x\inv = x \inv x = 1 \in X$, and $1 \cdot x = x \cdot 1 =x$; and $|X X \inv X | /|X| \leq k |X|$.  Then \thmref{gc0} is also valid for Freiman approximate groups.    In particular
$X$ has a large subset $X_2$ closed under $[,]$ if not under $\cdot$, and in fact with   $[X_2,X_2]^2 \subseteq X_2$.  This again suggests that approximateness can only really enter via an Abelian part of a structure.  This "local"
version uses local versions of the theory of  locally compact groups  due to Goldbring \cite{goldbring}. 

  The finiteness assumption on $X$ in the above results is really only used via the counting measure ``at the top dimension", so they remain valid  in a measure-theoretic setting, see  \thmref{gc}.   
  
The remaining corollaries of \thmref{stabilizer}  attempt to make a stronger use of finiteness.  They
are proved directly, without the Lie theory, and   go in a somewhat complementary direction.  
The first assumes that the group generated by an approximate subgroup $X$ is perfect in a certain strong
statistical sense.  The 
conclusion is that $X$ is close to an actual subgroup.  We write 
$a^X = \{x \inv a x:  x \in X \}$.

\<{cor} \lbl{anti-nil} For any $k,l,m \in \Nn$,  
  for some    $p<1$, $K \in \Nn$,  we have the following statement.

Let $G$ be a group, $X_0$ a finite subset, $X=X_0 \inv X_0$.  Assume  $|X_0 X| \leq k |X_0|$.
Also assume that 
 with probability $\geq p$, 
an $l$-tuple $(a_1,\ldots,a_l) \in X^l$ satisfies: $|a_1^X \cdots a_l^X | \geq |X| /m$.  

Then there exists a subgroup $S$ of $G$, $S \subseteq X^2$, such that $X$ is contained in $\leq K$ cosets of $S$.  \>{cor}

We could use $(a^X \union (a \inv)^X)^{(l)}$ (or $a^{X_0}$) in place of $a^X$ above.    See \thmref{anti-nil2} for a weaker   alternative version of the hypotheses.    
$p$  can be taken to be a recursive functions of $k,l,m$,  but I have made no attempt to estimate it.  As Ward Henson pointed out, the proof does give an explicit estimate for $K$.   
 The proof also shows that $X$ normalizes $S$.   Laci 
  Pyber remarked that with this strengthening (but not without it), the conclusion implies small tripling for $X$.
  
Here and later on, when confusion can arise between iterated set product and Cartesian power, 
we use $Y^l$ to denote the former, and $Y^{(l)}$ for the latter.  

The assumption of \corref{anti-nil} may be strong in a general group theoretic setting, but
it does hold for sufficiently dense subgroups of simple linear groups.     The proof uses an idea originating in the Larsen-Pink classification of large finite simple linear groups, \cite{larsen-pink},  somewhat generalized and formulated as a  dimension-comparison lemma in \cite{hw}.  We obtain:

 \<{thm} \lbl{linear}  Let $G$ be a semisimple algebraic group, defined over $\Zz$, and $k$ an integer.    Then for some integer $k'$, the following statement holds.
 Let $K$ be a field, $X$ a finite subset of  $G(K)$ with $|XX  \inv X | \leq k |X |$.  Then  there exists a 
   subgroup  $H $ of $G(K )$  such that $|X  /H | \leq k'$, and either $H $ is (the set of $K $-points of) a connected proper algebraic subgroup of $G$ of   degree $\leq k'$, or 
 $H \subseteq (X \inv X)^2$.  
  \>{thm}  

Here "degree $\leq k'$" means   that if we view $G$ as a subset of
 the $n \times n$ matrices $M_n$, then $H_i$ is the intersection of $G$ with a subvariety 
 of $M_n$ cut out by polynomials of   degree $\leq k'$.   Thus if the group generated by $X$ is sufficiently   Zariski dense, 
 $X$ will not be contained in such an algebraic subgroup, so that $X \inv X$ must be commensurable to a subgroup.  A special case:

\<{cor} \lbl{1.4}   For any $n \in \Nn$, for sufficiently large $n' \in \Nn$, the following holds.  Let $X$ be a finite subset of $GL_n(K)$, $K$ a field, with $|X| \geq n'$.    Assume  $|X  X  \inv X | \leq k |X |$, and that $X $ generates an almost simple group $S $.   Then $(X  X  \inv)^2=S $.  \>{cor}
 
Here $S $ is not {\em assumed} to be finite.   ``Almost simple" means:  perfect, and simple modulo a center  of bounded size.  The proof also shows that $X  X  \inv$ contains $99\%$ of the elements of $S$;   and that $X  X  \inv X  = S $; see proof and remarks following \propref{1.3p}.

For $S=SL_2(F_q)$ and $SL_3(F_p)$, \thmref{linear} follows from results 
Helfgott \cite{helfgott}, \cite{helfgott3} and Dinai;  for $G=SL_2(\Cc)$ and $G=SL_3(\Zz)$, \thmref{linear} follows from 
     \cite{chang} and \cite{eg}.  These authors all make a  much
weaker  assumptions on a subset $X$ of a group, namely $|X X \inv X | \leq |X|^{1+\e}$ for a small $\e$.  The combinatorial regime they work in is also meaningful model-theoretically  (cf. \exref{dim}), but we do not study it at present.

 Stable group theory includes a family of related results; for instance, the group law may be given by a multi-valued or partial function.  The partial case has antecedents in algebraic geometry, in Weil's group chunk theorem.  
 A  version of the partial case, including the Freiman approximate groups mentioned above, will be briefly noted
 in the paper.      
  It is   likely that the multi-valued case too admits finite combinatorial  translations along similar lines.    
 
In \secref{1}  we introduce the model-theoretic setting, and prove the independence theorem  and the stabilizer theorem in a rather general context.  In the presence of a $\si$-additive measure the stabilizer sounds 
   close to Tao's noncommutative Balog-Szemeredi-Gowers theorem (\cite{tao-nc}), while the independence
 theorem is, in the finite setting, extremely close  to the Komlos-Simonovitz corollary \cite{komlos-s}
 to Szemeredi's lemma (as I realized recently while listening to a talk by M. Malliaris.)  It is thus quite possible that combinatorialists can find other proofs of the results
 of  \secref{1}  and skip to the next section.   I find the   independent, convergent development of the two fields rather 
fascinating.  

   All  the  results we need from stability will be explicitly defined and proved.    \thmref{gc0} (and the more detailed
   \corref{gc}) are proved in \secref{lie}.  The   methods here are very close to \cite{hpp}; however we do not assume NIP.  This is in line with a sequence of realizations in recent years that
 tools discovered first in the stable setting are in fact often valid, when appropriately formulated, for first
 order theories in general.       \thmref{linear} is proved in \S 5.   
 
 \S 6 contains a proof that the topology on the associated Lie group is generated by the image of a definable family of definable sets.  
 
 In \S 7, we use the techniques of this paper along with Gromov's proof of the polynomial
growth theorem, to show (for any $k$)  that if a finitely generated group is {\em not} nilpotent-by-finite, it has a finite set of generators contained in no $k$-approximate subgroup.

This paper was prepared for a talk at the Cherlin Bayram{\i}  in Istanbul in June 2009.  I am grateful to Dugald Macpherson for a valuable conversation on this subject in Leeds.   Thanks also to Zo\'e Chatzidakis, 
 Lou van den Dries, Ward Henson,  Itay Kaplan, Krzysztof Krupinski, Elon Lindenstrauss, Dugald, Anand Pillay, Fran\c coise Point, Laci Pyber, Tom Scanlon, Pierre Simon, Terry Tao,   Alex Usvyatsov, and   two anonymous referees, for many useful comments.

\ssec{Basic model theory:  around compactness}  \lbl{basic}

We recall the basic setup of model theory, directed to a large extent at an efficient use of the compactness theorem.
We  refer to the reader to a book such as 
\cite{CK}, \cite{marker},\cite{poizat} or the lecture notes in \cite{pillay} for a fuller treatment.   We assume knowledge of  the the definition of a  first-order formula,
and of the compactness theorem, asserting that a finitely satisfiable set of formulas is satisfiable in some structure.   
 
 Let $L$ be a fixed language, $T$ a theory, $M$ a model. 
 We will occasionally use notation
 as if the language is countable (e.g. indices named $n$), but  this will not be
  really assumed unless explicitly indicated.  At all events for much of this paper, a language with a symbol
  for multiplication and an additional unary predicate will be all we need.

  $A$ will refer to a subset of  $M$.   We will assume $L,A$ are countable (this is quite inessential, and will be used only to avoid the need for
 cardinal parameters in discussing saturation below.)       We   expand $L$ to a language $L(A)$ with an additional constant symbol for each element of $A$.   The $L$-structure   $M$  is tautologically expanded to an $L(A)$-structure, and the result
 is still denoted $M$, by abuse of notation.   $T(A)$ is the $L(A)$- theory of $M$.
 $L_x(A)$ denotes the Boolean algebra of formulas of $L(A)$ with free variables $x$,
  up to $T(A)$-equivalence.  $S_x(A)=Hom(L_x(A),2)$ is the Stone space, or the 
  space of {\em types}.     A subset of $L_x(A)$ is {\em finitely satisfiable} if each finite subset has a common solution in $M$.   
  A {\em type} in a variable $x$, over $A$, is a maximal finitely satisfiable  subset of   $L_x(A)$.
  For   an element or tuple $a$ over a subset $A$ of a model $M$, $tp(a/A) = \{\phi(x) \in L(A): M \models \phi(a) \}$;
  if $tp(a/A)=p$ we say that $a$ realizes $p$.     An $A$-{\em definable set} is the solution set of some
 $\phi \in L(A)$.       
 It is an easy corollary of the compactness theorem that every  theory $T$ has models $\Uu$ with the following  properties holding for every small substructure $B$ of $\Uu$.    Here let us say $B$ is {\em small} if $2^{|B|} \leq |\Uu|$.  

 \<{enumerate}   
\item{Saturation:}   Every type over $B$  is realized in $\Uu$.
 \item{Homogeneity:} For $c,d$ tuples from $M$,   $tp(c/B)=tp(d/B)$ iff there exists $\si \in Aut(M/B)$  with $\si(c)=d$. 
\>{enumerate}

(In fact (1) implies (2) if the generalized continuum hypothesis holds; moreover in this case $\Uu$ is determined
up to isomorphism by $T$ and by $|\Uu|$, provided $T$ is complete.)

Given a complete  theory $T$, we fix a model $\Uu$ of $T$
with the above properties and with $|\Uu|>>\aleph_0$ (if it is not finite),  and interpret definable sets as subsets of $\Uu^n$.   
We will occasionally consider elementary submodels $M$ of $\Uu$; these will just be referred to as {\em models}.   
 We write $A \leq M$ to mean that $A$ is a substructure of $M$.      As stated above,  substructures $A$ will be taken to be countable;
"small" would be enough.  

A partial type over $A$
is any collection of formulas over $A$, in some free variable $x$, and closed under implication
in the $L_A$-theory of $M$.    

 The solution sets $D$ of partial types $r$ (over various countable sets $A$) are called  $\inft$-definable (read: $\infty$-definable) sets; so an $\inft$-definable set over $A$ is any intersection of $A$-definable sets.  The correspondence
 $r \mapsto D$ is bijective, because of the saturation property (1) above.   Complements of $\inft$-definable sets are called $\bigvee$-definable.   An equivalence relation is called $\inft$-definable if it has an $\inft$-definable graph.   It follows from saturation that 
 an $\inft$-definable set is either finite or
has size $|\Uu|$; a $\bigvee$-definable set is either countable or has size $\Uu$; 
  an $\inft$-definable equivalence relation  has either $\leq 2^{\aleph_0}$ classes or 
$|\Uu|$-classes.  Since $|\Uu|$ is taken to be large, this gap lends sense to the notion
of {\em bounded size} for   sets and quotients at these various levels of definability.

Another consequence of countable saturation is that projections commute with countable decreasing intersections:
$$(\exists x) \bigwedge_{i=1}^\infty \phi_i(x,y) \iff \bigwedge_{i=1}^{\infty} (\exists x) \phi_i(x,y)$$
provided that $\phi_{i+1} $ implies $\phi_i$ for each $i$.   The condition on the left beginning with $(\exists x)$
seems to be stronger, but compactness assures that the weaker condition on the right suffices for the existence of $x$ in some model, and countable saturation implies that such an $x$ exists in the given model.  In particular, the projection of an $\inft$-definable set is $\inft$-definable.  We will use this routinely in the sequel.  
Specifically,
if $Q$ is an $\inft$-definable subset of a definable group (see below), then the product set $QQ=\{x: (\exists y,z \in Q)(x=yz) \}$  is also $\inft$-definable.

By a {\em definable group} we mean a definable set $G$ and a definable subset $\cdot$ of $G^3$,
such that $(G(\Uu),\cdot(\Uu))$ is a group.     An $\inft$-definable subgroup is an $\inft$-definable set
which is a subgroup.  It need not be an intersection of definable subgroups.   We insert here a lemma
that may clarify these concepts.

A subset of a set $X$ is {\em relatively definable} if it has the form $X \meet Z$ for some definable $Z$.  

 \<{lem} \lbl{finiteindex} Let $G$ be a definable group.  Let $X$ be an $\inft$-definable subset of $G$, $Y$ a $\bigvee$- definable subset of $G$, and assume $X$ and $X \meet Y$ are subgroups of $G$, and $X \meet Y$
 has bounded index in $X$.  Then $X \meet Y$ is relatively definable in $X$, and has   finite index in $X$.  \>{lem}
 
 \prf   By compactness, $[X: X \meet Y] < \infty$:  otherwise one can find an infinite sequence
 $(a_i)$ of elements of $X$ such that $a_i a_j \inv \notin Y$ for $i \neq j$; but since these are $\inft$-definable conditions,  arbitrarily long sequences with the same property exist.  So $X \meet Y$ has finitely many distinct cosets
 $C_1,\ldots,C_n$ in $X$.   Note that $X \m C_i$ is $\inft$-definable.  Hence
 $C_j = \meet_{i \neq j} (X \m C_j)$ is   $\inft$-definable for each $j$.  Since $X_i$ and $C \m X_i$ are $\inft$-definable,
 they are relatively definable in $X$.  \eprf

A $\Uu$-definable set is $A$-definable iff it is $Aut(\Uu/A)$-invariant. 
  The same
is true for $\inft$-definable sets  
and for $\bigvee$-definable sets.  

Types over $\Uu$ are also called global types.

  A sequence $(a_{i}: i  \in \Nn )$ of elements of $\Uu$ is called {\em $A$- indiscernible} if  any order-preserving map $f: u \to u'$ between two finite subsets of $\Nn$ extends to an automorphism of $\Uu$ fixing $A$.  The same applies to sequences of $n$-tuples.   Using Ramsey's theorem and compactness,
  one shows that if $(b_i: i \in \Nn)$ is any sequence, there exists an indiscernible
  sequence $(a_i: i \in \Nn)$ such that for any formula $\phi(x,y)$, if $\phi(b_i,b_j)$
 holds for all $i<j$, then $\phi(a_i,a_j)$ holds for all $i<j$.  
  A theorem of Morley's \cite{morley}
 asserts the same thing with the formulas $\phi(x,y)$ replaced by {\em types}, provided 
 $\Nn$ is replaced with a sufficiently large cardinal.  For certain points (outside the main line),
 we will use Morley's theorem as follows.  Let $q$ be a global type, and construct a sequence
 $a_i$ inductively, letting $A_i = \{a_j: j < i \}$, and choosing $a_i$ such that $a_i \models q|A_i$.
 By Morley's theorem, there exists an indiscernible sequence $(b_0,b_1,\ldots)$ such that
 for any $n$, $b_n \models q_n | \{b_0,\ldots,b_{n-1}\}$ for some $Aut(\Uu)$-conjugate $q_n$ of $q$.  
 
We will say in this situation that $(b_0,b_1,\ldots)$ are $q$-indiscernibles.  
The main case is that $q$ is an invariant type,  
and then Morley's theorem is not needed, for the original $(a_j)$ are automatically indiscernible; see
\cite{pillay-book}.    A global type finitely satisfiable in $M$ is always $M$-invariant.  In particular,
given any type over $M$, this yields an $M$-indiscernible sequence $(a_i)$ such that $tp(a_i/M \union A_i)$ does not
fork over $M$.  (cf. \secref{stability-sec} for the definition.)  We remark that 
Morley's theorem uses more infinite cardinals than the rest of the paper (namely, not only infinite sets but arbitrary countable iterations of the power set operation.)

 In all notations, if $A$ is absent we take $A=\emptyset$.    Generally a statement made for $T_A$
  over $\emptyset$ is equivalent to the same statement for $T$ over $A$, so no generality is lost.
 
We will occasionally refer to ultraproducts of a family $M_i$ of $L$-structures.  They are specific way of
constructing models $M$  of the set of all sentences holding in all but finitely many    $M_i$, and they have the saturation property (1).  No other properties of ultraproducts will be needed.

\section{Independence theorem}  \lbl{1}

\ssec{Stability}  \lbl{stability-sec}

The material in this subsection is a  presentation of \cite{kim-pillay}, Lemma 3.3, here   \lemref{stat}; compare
also  \cite{lazy} \S3, and the stability section in \cite{BU}.
 
 Let $T$ be a first-order theory, $\Uu$ a universal domain.   One of the main lessons of stability is the usefulness of
 $A$-invariant types, meaning $Aut(\Uu/A)$-invariant types.  We note that
 if a global type $p$ is finitely satisfiable in some $A \leq M$, then $p$ is 
 $A$-invariant:  if $a,a'$ are $Aut(\Uu/A)$-conjugate, then $\phi(x,a) \& \neg \phi(x,a')$
 cannot be satisfied in $A$.  
 
 We say $A$ is  an {\em elementary submodel} of $\Uu$ (written:  $A \prec \Uu$) if any nonempty $A$-definable set
 has points in $A$.  If $A \prec \Uu$, then  any type over $A$ extends to a global type, finitely satisfiable in $A$.  (\cite{pillay}).

Consider two partial  types $r(x,y), r'(x,y)$ over $A$.  Say $r,r'$ are  {\em stably separated} if there is no   sequence $((a_i,b_i): i \in \Nn)$
 such that $r(a_i,b_j)$ holds for $i<j$, and $r'(a_i,b_j)$ holds for $i>j$.    Note that if arbitrarily long
such sequences exist then by compactness an infinite one exists, and in fact one can take the $(a_i,b_i)$ to form an 
 $A$- indiscernible sequence.   Moreover $r,r'$ are stably separated iff they contain formulas $\phi,\phi'$ that are stably separated.  By reversing the ordering one sees that stable
 separation is  a symmetric property. 
 
We say $r'$ is {\em equationally separated} from $r$ if  there is no   sequence $((a_i,b_i): i \in \Nn)$ such that 
$r(a_i,b_i)$ holds for all $i$, and 
$r'(a_i,b_j)$ holds for $i<j$.  This is an asymmetric  condition, that implies stable separation:  if stable separation fails, so that $r(a_i,b_j)$ holds for $i<j$, and $r'(a_i,b_j)$ holds for $i>j$ in  some sequence $(a_i,b_i)$,
 the shifted subsequence $(a_{2i}, b_{2i-1})$ shows that equational separation fails too.


If $r,r'$ are stably separated then they are mutually inconsistent, since if $r(a,b)$ and $r'(a,b)$ we can
let $a_i=a,b_j=b$.  In stable theories, the converse holds.

Note that the set of stably separated pairs is open in the space $S_2^2$ of pairs of 2-types. Any extension
of a stably separated pair to a larger base set remains stably separated. 

A partial type  $r'(x,b)$ is said to divide over $A$ if there exists an indiscernible sequence $b_0,b_1,\ldots$ over $A$ such that   $\union_i r'(x,b_i)$ is inconsistent, and $tp(b/A)=tp(b_i/A)$.   Equivalently, for some $k$,  $\{r'(x,b_i):
i \in w \}$ is inconsistent for any $k$-element subset $w$ of $\Nn$.   By compactness,  $r'(x,b)$ divides over $A$
iff some formula $R(x,b) \in r'(x,b)$ divides over $A$.   The ideal generated by all formulas that divide over $A$ is called the {\em forking ideal};
thus $\phi(x,c)$ forks over $A$ if it implies a disjunction of formulas that divide over $A$.

If   $q=q(y)$ is a global type, we say that $r'(x,y)$ $q$-divides over $A$ if for some $n$, if $b_i \models q | A(b_0,\ldots,b_{i-1})$ for $i \leq n$,
then $\union_{i  \leq n}  r'(x,b_i)$   inconsistent.   This is equivalent to dividing, with the additional requirement that the indiscernible sequence be $q$-indiscernible.


 \<{lem}  
 \lbl{ws2} Let $r,r'$ be stably separated formulas over $A$.  
 Let $q(y)$ be an $A$-invariant global   type.  Assume $r'(a,y) \in q$, $p=tp(a/A)$.  
 Then $p(x) \union r(x,y)$   $q$-divides over $A$.
\>{lem}
 
\prf      Suppose it does not.  

 Define  $a_1,\ldots,  c_1,\ldots$ inductively:   
 given $a_1\ldots,a_{n-1},c_1,\ldots,c_{n-1}$, 
 choose $c_n$ such that 
$c_n \models q| \{a_1,\ldots,a_{n-1},c_1,\ldots,c_{n-1}\}$, and $a_n \models p$ chosen with $r(x,c_i)$ for $i < n$.   The latter choice is possible 
since  $p(x) \union r(x,y)$ does not $q$-divide over $A$.

Then $r'(a_i,c_j)$ holds if $i<j$, but $r(a_i,c_j)$ holds when $i>j$.  This    contradicts the stable separation of $r,r'$.  \eprf

  We say that an $A$- invariant relation $R$ is a {\em stable relation} over $A$ if whenever $(a,b)\in R$ and $(a',b') \notin R$, 
 $tp((a,b)/A)$ and $tp((a',b')/A)$ are stably separated.   If $R$ is stable, so is the complement of $R$; but we are interested mostly in $\inft$-definable $R$.  
 
We will also encounter the   condition of {\em equationality}.  $R$ is {\em equational} if whenever $(a,b)\in R$ and $(a',b') \notin R$, 
 $tp((a',b')/A)$ is equationally separated from  $tp((a,b)/A)$.   As we have seen that equational separation implies stable separation,
 equationality implies stability.
 
 We will say:  ``$R(a,b)$ holds" for ``$(a,b) \in R$".  
   When $q$ is a global type, write ``$R(a,y) \in q(y) $" to mean:  $R(a,b)$  holds when $b \models q|A(a)$.
 
\<{lem} \lbl{stat}  Let $p(x)$ be a type over $A$,
and $q(y)$ be a global, $A$-invariant type.    Let $R$ be a stable relation over $A$. 

(1)  Assume $R(a,b)$ holds with $a \models p, b \models q|A(a)$.   Then $R(a',b)$
holds whenever $a' \models p$ and $tp(a'/Ab)$ does not divide over $A$.

(2)  Assume $tp(a/A)=tp(a'/A)$, $b \models q$,  and neither  $tp(a/Ab)$ nor $tp(a'/Ab)$ 
divides over $A$.    Then $R(a,b)$ implies $R(a',b)$.     

(3)  Assume $p$ too extends to a global, $A$-invariant type.  Let $E =\{(a,b): a \models p, b \models q|A \}$.
  Then the eight conditions:

{ $R(a,b)$ holds for  {some}/{all} pairs $(a,b) \in E$ such that $tp(a / A(b))$ /$tp(b/A(a))$  does not fork / divide over $A$ }

are all equivalent.
\>{lem}  

\prf    (1)  Suppose $R(a',b)$ fails to hold.  So $tp(a',b)$ and $tp(a,b)$ are stably separated,
say by formulas $r',r$.   By \lemref{ws2}, since $r$ holds for $b \models q|A(a)$,
$r'(x,b) \union p(x)$ divides, so $tp(a'/Ab)$ divides over $A$, a contradiction.  

(2)   Let $R'$ be the complement of $R$; it is also a stable relation.  
 Let $c \models q|A(a)$.  
If $R(a,c)$ holds then by (1) we have $R(a',b)$ and $R(a,b)$.  If   $R(a,c)$ holds
then similarly $R'(a',b)$ and $R'(a,b)$.  In any case we have $R(a,b) \iff  R(a',b)$, so the stated implication holds.

(3)  
Let $E'$ be the set of pairs $(a,b) \in E$ such that $tp(a/A(b))$ does not divide over $A$, and $E''$
 the set of pairs $(a,b) \in E$ such that $tp(b/A(a))$ does not divide over $A$.  
The equivalence between the four conditions for $tp(a/A(b))$   follows from (2):
if $R(a,b)$ holds for some pair such that $tp(a/A(b))$ does not fork, then in particular it holds for a pair in $E'$ (the same pair); 
 by (2), it holds for al such pairs; hence certainly for all pairs for which $tp(a/A(b))$ does not fork over $A$.

Thus a single truth value   for $R$ is associated with pairs $(a,b) \in E'$.    Similarly, as the conditions are symmetric,
 a single truth value    for $R$ is associated with pairs $(a,b) \in E''$.   
  It remains
to show that these truth values are equal.  Replacing $R$ be its complement if necessary, we may assume $R(a,b)$ holds in the situation of (1),
where $b \models q | A(a)$.  In particular $tp(b/A(a))$ does not fork over $A$; so $R(a',b')$ holds for all $(a',b') \in E''$.
But (1) asserts that $R(a',b)$ holds for all $(a',b) \in E'$.  Hence $R$ holds for all pairs in $E' \union E''$.     \eprf

 Non-dividing in \lemref{stat} (3) can be replaced by any stronger condition;  non-forking was mentioned above; we will later use smaller ideals.

\<{remark}     Let $p,q,R$ be as in \lemref{stat}, with $R(x,y)$   equational.   Let $Q =\{b: b \models q | A \}$, $P= \{a: a \models p|A \}$.
If $R(a,b)$ holds with $a \models p, b \models q|A(a)$, then $P \times Q \subset R$.  \>{remark}

 \<{lem} \lbl{bounded-nf}     Let $S=S^{nf}_z$ be the set of global types that do not fork over $\emptyset$.  Define an equivalence
 relation $E=E_{st}$ on $S$:
 $p E_{st} p'$ iff for any stable invariant relation $R$, and any $b$, we have $R(b,z) \in p \iff R(b,z) \in p'$.  Then
 $|S/E| \leq 2^{|T|}$.   \>{lem}
 
 \prf  Let $M$ be a model.  It suffices to show that if $p|M = p'|M$ then $pE_{st} p' $.   Let $R(x,z)$ be a stable relation.
 Let $q=tp(b/M)$, and let $q^*$ be any $M$-invariant global type extending $q$.  Let $c \models p|M$.
    By \lemref{stat}, since $p,p'$ do  not fork over $M$, $R(b,z) \in p$ iff $R(x,c) \in q^*$ iff $R(b,z) \in p'$.   
 \eprf
     
\ssec{Making measures definable}  
 \lbl{measurequantifiers}
 
A Keisler measure $\mu_x$  is a finitely additive real-valued probability measure on the formulas (or definable sets) $\phi(x)$ over the universal domain $\Uu$.   
 See   \cite{hpp}.  
 
 We say $\mu$ is {\em $A$-invariant} if for any formula $\phi(x,y)$, for some function $g: S_y(A) \to \Rr$,
 we have  
 $\mu ( \phi(x,b)) = g(tp(b/A))$ for all $b$.  If in addition $g$ is continuous, we say that $\mu$ is an $A$-definable
 measure.     
 
 Let $M_i$ be a family of finite $L$-structures.     We wish to expand
 $L$ to a richer language $L[\mu]$, such that each $L[\mu]$ structure admits
 a canonical definable measure $\mu$.  For each formula $\phi(x,y)$ and $\a \in \Qq$
 we introduce   a formula $\theta(y) =(Q_\a x)\phi(x,y)$ whose intended intepretation is:   $\theta(b)$ holds iff $\mu_x \phi(x,b) \leq \a$.  If we wish $\mu$ to measure
 new formulas as well as $L$-formulas, this can be iterated.
 
 We can expand each $M_i$ canonically to $L[\mu]$,  interpreting the formulas
 $(Q_\a x)\phi(x,y)$ recursively  using the counting measure. 
  
    Let $N$ be any
 model of the set of sentences true in all $M_i$ (such as ultraproduct of the $M_i$ with respect to some ultrafilter.)
 Define $\mu \phi(x,b)   = \inf \{ \a \in \Qq: (Q_\a x)\phi(x,b) \}$.   Then $\mu$ is a Keisler measure.   The 
 formulas $(Q_\a x) \phi$ may not have their intended interpretation with respect to $\mu$ exactly,
 but very nearly so:  $(Q_\a x) \phi(x,b)$ implies $\mu_x \phi(x,b) \leq \a$, and is implied by
 $\mu_x \phi(x,b) < \a$.  Thus $\mu$ is a definable measure on $N$.   
 
 We will actually only use the corollary that the 0-ideal of $\mu$ is an invariant ideal, see below.

\ssec{Ideals}   Let $X$ be a definable set, over $A$.

$L_X(\Uu)$ denotes the Boolean algebra of $\Uu$-definable subsets of $X$. 
An ideal $I$ of this Boolean algebra is {\em $A$-invariant} if it is $Aut(\Uu/A)$-invariant;
equivalently  $I$ is 
a collection of formulas of the form $\{\phi(x,a): tp(a/A) \in E_\phi \}$, where
for each $\phi(x,y)$, $E_\phi$ is a subset of $S_y(A)$, and $\phi(x,a)$ implies $x \in X$.  
To emphasize the variable, we use the notation $I_x$.

 We say $I$ is {\em $\inft$-definable}   if for any $\theta(x,y)$, the set 
 $\{b: \theta(x,b) \in I \}$ is $\inft$-definable.   Similarly for $\bigvee$-definable.

  We say a partial type $Q$ over $A$  is $I$-wide if it implies no formula in $I$.
  
In case $X$ is $\bigvee$-definable, i.e. a countable union of  $A$-definable sets $X=\bigvee_i X_i$, we let $L_x(\Uu) = \union_i L_{X_i}(\Uu)$.  
An {\em ideal} of $L_x(\Uu)$ is a subset $I$ such that $I \meet L_{X_i}(\Uu)$ is an ideal for each $i$; it is called
   $A$-invariant, $\inft$-definable or $\bigvee$-definable  if  $I \meet L_{X_i}(\Uu)$ has the corresponding property, for each $i$.  
 
 By analogy with measures, we will sometimes denote ideals in a variable $x$ by 
 $\mu$, and write $\mu(\phi)=0$ for $\phi \in \mu$, and $\mu(\phi)>0$ for $\phi \notin \mu$.

The following definition is the defining property of $S1$-rank, \cite{PAC}, relativized
to an arbitrary ideal (so within a definable set of finite S1-rank, the definable sets of
smaller S1-rank form an S1-ideal.)  The terms {\em invariant, formula, indiscernible} are understood over
some fixed base set $A$.

\<{defn} \lbl{S1} An invariant ideal $I=I_x$ on $X$ is {\em S1} if for any formula $D(x,y)$ and indiscernible
$(a_i: i \in \Nn)$ with $D(x,a_i) \in L_X(\Uu)$, if  
   $D(x,a_{i}) \meet D(x,a_{j}) \in I$ for $i \neq j$,   then some $D(x,a_{i})  \in I$. 
   \>{defn}

The forking ideal is 
contained in any S1-ideal:  
 
 \<{lem}\lbl{fork} Let $I$ be an invariant S1 ideal over $A$.  
   If  $\phi(x,b)$ forks over $A$ then $\phi(x,b) \in I$.  \>{lem}
 
 \prf  It suffices to show that if $\phi(x,b)$ divides over $A$, then  $\phi(x,b) \in I$.
 Let $(b_i)$ be an $A$-indiscernible sequence, with $\{\phi(x,b_i)\}$ inconsistent;
 so for some $k$, $\phi(x,b_1) \wedge \ldots \wedge \phi(x,b_k) = \emptyset$.  If
 $\phi(x,b_1) \in I$ we are done.   Otherwise let $m$ be maximal such that $\phi(x,b_1) \wedge \ldots \wedge \phi(x,b_m) \notin I$.  Let $c_i=(b_1,\ldots,b_{m-1},b_{m+i})$, and let
 $\psi(x,c_i) = \phi(x,b_1) \wedge \ldots \wedge \phi(x,b_{m-1}) \wedge \phi(x,b_{m+i})$.  Then the intersection of any two $\psi(x,c_i)$ is  in $I$, but no $\psi(x,c_i)$ is in $I$.  This contradicts \defref{S1}.
  \eprf
  
The forking ideal over $A$ is also invariant under all $A$-definable bijections; in particular  for subsets of  a group $G$
under left and right translations by elements of $G(A)$, i.e. by elements of $G$ definable over $A$.  This will not be of real use to us however as we will be interested in translation invariance, right and left, by elements not necessarily
defined over $A$.
   
A  fundamental observation from \cite{CH}, \cite{PAC}, and \cite{kim-pillay}:

\<{lem}\lbl{stable} 
 Let $I_z$ be an invariant S1-ideal. 
Let $P=P(x,z),Q=Q(y,z)$ be formulas.   Define:
$$R(a,b) \iff (P(a,z) \wedge Q(b,z)) \in I_z$$
     Then    $R$ is a stable  invariant relation.  \>{lem}

\prf  We show indeed that $R$ is   equational: if $R(a_i,b_j)$ holds for $i<j$, where $(a_i,b_i)_i$ is indiscernible,
then $R(a_i,b_i)$ holds too.  

  Otherwise, let $C_i = \{z:   P(a_i,z) \wedge Q(b_i,z) \}$.  Then $ C_i \notin I_z$ but $\mu_z(C_i \meet C_j) = 0$.  This contradicts the S1 property of \defref{S1}.

  \eprf

  \<{example}  \lbl{measures-epsilon}  \rm Let $\mu(z)$ be a Keisler measure on $\Uu$-definable
subsets of a set $Z$, with $\mu(Z)=1$.   Let $e \in \Nn, \e=1/e > 0$.   Let $\phi(x,z), \phi'(y,z)$ be formulas, and
write $D(a,b)= \{z \in Z: \phi(a,z) \meet \phi'(b,z) \}$.  Let $r(x,y), r'(x,y)$ be 
formulas such that
if $r(a,b)$ then  $\mu(D(a,b)) \geq   \e$, while $r'(a,b)$ implies 
$\mu(D(a,b))  < \e^2/2$.  
Then $r,r'$ are stably separated;  indeed $r$ is equationally separated from $r'$.   \rm  For suppose $r'(a_i,b_j)$ holds for $i=1,\ldots, 2e$.  Let $D_i = D(a_i,b_i)$.  Then $\mu(D_i) \geq \e$,
but $\mu (\union_{1 \leq i < j \leq 2e} D_i \meet D_j) <   (2e(2e-1)/2 ) (\e^2 /2) < 1$.  So $\mu(\union_i D_i) >
2e \e -1 =1$, a contradiction.  
 \>{example}

 \<{example} \lbl{measures} \rm Let $\mu$ be an $Aut(\Uu/A)$-invariant, real-valued, finitely additive  measure on $\Uu$-definable sets.    Then $I=\{\phi(x,b): \mu(\phi(x,b)) =0 \}$ is an  $Aut(\Uu/A)$-invariant S1-ideal.  It is $\inft$-definable
 if $\mu$ is definable. \>{example}

 \<{example} \lbl{dim} \rm  Let $X$ have nonstandard finite size $\a$, and let $I$ be the ideal
 of all definable sets  with nonstandard size $\b$, where $\log(\b) \leq  (1-\e) \log(\a) $
 for some standard $\e>0$. (See \S 5 for detailed definitions.) Then $I$ is a  $\bigvee$-definable ideal.  It is not S1; but 
 the counterexamples are always families contained in a definable set of dimension
 $< \e \log(\a)$ for each $\e>0$.
    \>{example}

\ssec{Wide global types}  

We now note the existence of useful global types relative to an ideal $I$,
in three slightly different situations.   The combinatorial applications of the present
paper can be deduced from either   \lemref{claim0} or \lemref{sym}; the former
has a shorter, more general but much more impredicative proof.

   \<{lem} \lbl{sym-d} Let $I=I(x)$ be a $\bigvee$-definable ideal, defined over a model
 $M$.     
 Then there exists a global  type $p$, finitely satisfiable in $M$, such that if
 $b \models p|M$, $a \models p |M(b)$, then $tp(b/Ma)$ is $I$-wide.  (In fact, whenever $p$ is finitely satisfiable in $M$  and $p|M$ is wide, then $p$ has this property.)
 \>{lem}
 
\prf  Let $p_0$ be any wide type over $M$, and let $p$ be any extension to $\Uu$, finitely
satisfiable in $M$.  
 Let $b \models p|M$, $a \models p |M(b)$.  If $tp(b/Ma)$ is not wide, then for some $\phi(x,y)$
 we have $\phi(a,b)$ and $\phi(a,y) \in I$;  by $\bigvee$-definability,   for some $\theta \in tp(a/M)$, for all $a'$ with  $ \theta(a')$,  $\phi(a',y) \in I$.   Since $tp(a/Mb)$ is finitely satisfiable in $M$, there exists $a' \in M$ with $\theta(a')$ and $\phi(a',b)$.  It follows that $p_0=tp(b/M)$ is not wide, a contradiction.
 \eprf

 \<{lem}\lbl{claim0} Let $I=I(x)$ be an $A$-invariant ideal.  There exists a model $M \geq A$, a global $M$-invariant type 
$q$, finitely satisfiable in $M$, such that if  $a \models q|M$ and $b \models q | M(a)$ then $tp(a/M(b))$ is wide.   
\>{lem}

\prf  
 Let $T_{sk}$ be a Skolemization of the theory, in a expansion $L_{sk}$ of 
 the language $L$; so the $L_{sk}$-substructure $M(X)$ generated by a set $X$
 is an elementary submodel.
Define a sequence of elements $a_i$ ($i <\beth_{\omega_1}$), 
and sets $A_i = M(\{a_{j}: j<i \})$, with $tp_L(a_i/A_i)$ wide.  By Morley's theorem
\cite{morley}, there exists an indiscernible sequence 
$(c_i: i < \omega+2)$ such that for any $n$, for some $i_1< \ldots < i_n$,
$tp(c_1,\ldots,c_n) = tp(a_{i_1},\ldots,a_{i_n})$.     In particular, 
$tp(c_i / \{c_j: j <i \})$ is wide.  Let $U$ be an ultrafilter on $\Nn$, and
let $q$ be the set of   formulas $\phi(x)$ of $L(\Uu)$ such that $\{i: \phi(c_i) \} \in U$.
Let $M=A_\omega$.  Then $q$ is finitely
satisfiable in $M$.
Let $a=c_{\omega+1}$, $b=c_{\omega}$.  Then   $a \models q|M$ and $b \models q | M(a)$, and $tp(a/M(b))$ is wide.

\eprf

 \<{lem} \lbl{sym} Let $I=I(x)$ be an $\inft$-definable ideal, defined over a model
 $M$ with $L(M)$ countable.     Assume (``Fubini")  there exists an   ideal $I^2(x,y)$  
 on $L_{x,y}(M)$ 
 such that:   
 (i) if $\phi(a,y) \in I(y)$ whenever $tp(a/M)$ is $I$-wide, then $\phi \in I^2$; (ii) if
 $\phi(x,b) \in I(x)$ whenever $tp(b/M)$ is $I$-wide, then $\phi \in I^2$; 
 (iii) if   $\phi(x) \wedge \phi(y) \in I^2$ then $\phi \in I$.   
  
 Then there exists a global  type $p$, finitely satisfiable in $M$, such that if
 $b \models p|M$, $a \models p |M(b)$, then $tp(a/Mb)$ and $tp(b/Ma)$ are $I$-wide.

 \>{lem}

 \prf  Let $B$ be the Boolean algebra of formulas of $M$ modulo $I$.  We show that a generic ultrafilter $p_0$ on $B$
 (in the sense of Baire category) can be extended to a type satisfying the lemma.

\claim{}  Let    $\phi_i(x,y)$  ($i=1,2,3$) be a triple of formulas, and let   $P(x)  \in B \m I$.
Assume 
$$P(x) \wedge P(y) \vdash   \bigvee_{i=1}^3 \phi_i(x,y)$$
       Then for some $P' \in B \m I$ implying $P$,   for any $a,b \in P'$, we have (*): \, 
   $\phi_1(a,y) \notin I$ or $\phi_2(x,b) \notin I$ or $\phi_3(c,b)$ for some $c \in M$.

\prf    
If  $( P(x) \wedge \phi_3(c,x)) \notin I$ for some $c \in M$, we can let $P'(x) =  P(x) \wedge \phi_3(c,x)$;  
  then the third option in (*) is met.   
  Otherwise, $(P(x) \wedge \phi_3(c,x)) \in I$ for all $c \in M$. 
   It follows from the $M$-$\inft$-definability   of $I$ that $(P(x) \wedge \phi_3(c,x)) \in I$ for all $c$.  So $P(y) \wedge \phi_3(x,y) \in I^2$.
 
 If for some $P' \in B \m I$ implying $P$   we have:  $P'(a)$ implies
 $ \phi_1(a,x) \notin I$, then the first disjunct of (*) holds.   Otherwise, using the  $M$-$\inft$-definability of $I$,
 we see that for   all $a \in P$ with $tp(a/M)$ $I$-wide,  $ \phi_1(a,y) \in I$.  By the Fubini assumption (i),
 $ P(x) \wedge \phi_1(x,y)  \in I^2$.
 
 Similarly, if for some such $P'$, $P'(b)$ implies $\phi_2(x,b)>0$, then the second
 disjunct holds.  Otherwise, by Fubini (ii), $ (\phi_2(x,y) \wedge P(y)) \in I^2$. 
 
 Since $P(x) \wedge P(y)$ implies the disjunction of the $\phi_i$, we have
 $  (P(x) \wedge P(y)) \in I^2$; so $P \in I$;  this contradicts the choice of $P$, and proves
 the claim.    \eprf
 
 It is now easy to construct  a type $p_0$ over $M$ such that, for any  $\phi_1(x,y), \phi_2(x,y), \phi_3(x,y)$,  If $p_0(x) \union p_0(y) \vdash \bigvee_{i=1}^3 \phi_i(x,y)$, then (*) of the Claim holds for any $a,b \models p_0$.  Namely, 
 we let $p_0= \{P_n\}$, where $P_n \in B \m I$ is constructed recursively.   If $n$ is even, we choose $P_{n+1}$ so
 as to imply $\psi$ or $\neg \psi$, where $\psi$ is the $n/2$-nd element of
 some enumeration of the formulas $\psi(x)$.  If $n=2m+1$ is odd,
 consider the $m$'th triple $(\phi_1,\phi_2,\phi_3)$ in some (infinitely repetitive) enumeration of all triples of formulas over $M$.  If $P(x) \union P(y) \vdash \bigvee_{i=1}^3 \phi_i$, let $P'$ be as in the Claim, and let $P_{n+1} =P_n \wedge P'$.

  Let $b\models p_0$, and let 
$\G(x,b) =  p_0(x) \union \{\neg \phi_1 (x,b):   \phi_1(x,b) \in I \}
 \union \{ \neg \phi_2(x,b):    \phi_2(b,x) \in I \} \union 
 \{ \neg \phi_3(x,b):   (\forall c' \in M)(\phi_3(c',x) \notin p_0 )\}$.  If $\G(x,b)$
 is inconsistent, then $p_0(x) \union p_0(y) \vdash \phi_1(x,b) \vee \phi_2(x,b) \vee \phi_3(x,b)$ for some $\phi_1,\phi_2,\phi_3$ with $\phi_1(x,b) \in I, \phi_2(b,x) \in I,
 \phi_3$ such that $(\forall c' \in M)(\phi_3(c',x) \notin p_0 )\}$.  But this contradicts
 the construction of $p_0$.   Thus $\G(x,b)$ is consistent, and in view of the formulas $\neg \phi_3$, finitely satisfiable in $M$.   Let  $p$ be any   extension of $\G(x,b)$ to a global type finitely 
 satisfiable in $M$.  Let $b \models p | M, a \models p |M(b)$.  Then $tp(a/Mb)$ 
 is wide because of the formulas $\neg \phi_1$, and $tp(b/Ma)$ is wide because of the
 formulas $\neg \phi_2$.
   \eprf

We now come to the 3-amalgamation statement.  It says roughly that given a triangle of types, an arbitrary replacement of
one edge by another with the same vertices will not affect the wideness of the opposite vertex over
the edge.   To simplify notation we work over $A=\emptyset$, so ``divides" means ``divides over $\emptyset$."

 \<{thm} \lbl{3replacement} Let $\mu=\mu_z$ be an invariant S1- ideal.  
  Assume $tp(c/a,b)$ is $\mu_{z}$-wide,    $tp(b/a)$ and   $tp(b'/a)$ 
  do not divide,  
    $tp(a)$ extends to an invariant global type, 
and    $tp(b)=tp(b')$.  Then there
 exists $c'$ with $tp(c'/a,b')$ wide, and $tp(c'b')=tp(cb), tp(c'a)=tp(ca)$.   \>{thm}
 
 \prf   Let $Q \in tp(cb), P \in tp(ca)$.  By compactness, it suffices, for any such pair of formulas, 
  to find $c'$ with   $tp(c'/a,b')$ wide, and $Q(c',b'),P(c',a)$.  In other words
  it suffices to show that $\mu_{z} (Q(z,b') \wedge P(z,a)) >0$.
  
 Consider the relation $R(x,y)$ such that $R(d,e)$ holds iff $\mu_z(P(z,d) \wedge Q(z,e)) =0$.    By \lemref{stable}, it is a stable relation.  

    By assumption, $tp(b'/a)$ and $tp(b/a)$ do not divide.  By \lemref{stat}, 
 since $R(a,b)$ fails, $R(a,b')$ must fail too.  Thus $\mu_z(P(a,z) \wedge Q(bz)) >0$. 
 \eprf

 \<{remark} \lbl{3replacementr}    \<{enumerate}
 
 \item  The hypothesis that $tp(b/a)$ and   $tp(b'/a)$ do not divide can be replaced by: 
 $tp(b/a)$ and $tp(a/b')$ do not divide, using \lemref{stat} (3).

\item  Over a model, the hypothesis that $tp(a)$  extends to an invariant global type holds automatically.  

\item  If $E$ is an $\inft$-definable equivalence relation over $A$ with boundedly many classes, and the class
of $a$ is {\em not } the unique wide class within $tp(a)$, then $3$-amalgamation can fail;  one cannot amalgamate
a type $p(x,y)$ implying $\neg (xEy)$ with any types implying $xEy,yEz$.  It is possible that this is the only obstruction,
so that as in \cite{kim-pillay},  \thmref{3replacement} holds over  any set $A$ which is boundedly closed.  
\>{enumerate}
 \>{remark}
  
\ssec{Complements}
 In the remainder of this section we  mention a variant   of   \thmref{3replacement} in a measured setting, bringing out the 3-amalgamation aspect, and discuss connections to NIP and to probability theory.
  None of this will be needed for the combinatorial applications of \S 2-5.

An arbitrary triangle of 2-types cannot be expected to give a consistent 3-type,   for instance since a definable linear ordering may be present;
 types including $x<y,  y<z, z<x$ are obviously not consistent together.
But in a measured setting, contrary to  initial appearances, this obstruction has effect only on a measure zero set.

  Below, $i$ ranges over elements of $\Ups := \{1,2,3\}$, while $u$ ranges over subsets of $\Ups$ of size $2$.    Let $x_i$ be a sort, and
$X_i$ the space of types in this sort, over a fixed base set $M$.   We assume every type in $X_i$ extends to an invariant type (as is the case over an elementary submodel.)  We also assume, for simplicity's sake,
that $L(M)$ is countable.   
For $i \in \Ups$ let $\mu_i$ be an $M$-definable measure on $X_i = X_{x_i}$.  In fact it suffices to assume that $\mu_i$ is Borel-definable over $M$,
meaning that $\mu_i(\phi(x,b))$ is a Borel function of $tp(b/M)$.

Assume the $\mu_i$ commute, in the sense that  for any $i \neq j \in \Ups$,
for any formula $\phi(x_i,x_j)$ over $M$, 
$$\int  \mu_j(\phi(x_i,x_j)) d \mu_i =  \int \mu_i(\phi(x_i,x_j)) d \mu_j$$
see \cite{nip3}.  Any measures obtained as ultraproducts of counting measures will certainly have this property.
 
The common value
is denoted $\mu_{ij}(\phi)$; this defines a measure with variables $(x_i,x_j)$,   referred to as the tensor product of $\mu_i,\mu_j$.  Similarly, for 
$u \subseteq \Ups$, let $\mu_u$ be   the tensor product measure  on $X_u$.
  In particular we have $\mu=\mu_{123}$ on  $X_{123} = X(\Ups)$.  

We will occasionally refer to random elements; this can be given precise set-theoretic foundations, but we will not do this here.  Instead we will understand by this an  element 
of a type space, or a product of type spaces, avoiding a certain countable collection of measure-zero Borel sets, that
can be explicitly specified by inspecting the proof.    We will
also omit the foundational details of the notion of conditional measures, noting only  that in the context of separable totally disconnected spaces
we have a canonical countable Boolean algebra, namely the clopen subsets, making things easier.  

Consider the natural maps $X(\Ups) \to X(\{23\})$, $X(\Ups) \to X(2) \times X(3)$, etc.  For  any such map, with target $Y$ carrying measure $\mu_Y$, and given a random 
(for the pushforward measure of $\mu$) element $y \in Y$, we let $X_{123}(y)$ denote $X_{123}$ with the measure $\mu_{123/y}$ conditioned on $y$.  These conditional measures concentrate on the fiber over $y$, and satisfy:  $\mu(B) = \int (\mu_{123/y}(B)) d \mu_Y$ for any clopen $B$.  This formula defines $\mu_{123/y}$ uniquely for random $y$, in the sense that any two choices will agree for almost all $y$.   Again we refrain from giving the  foundational details, noting only that they are much easier in the present  context of separable totally disconnected spaces; this is due to the availability  of a  canonical countable Boolean algebra generating the measure algebra,  namely the clopen subsets.   See \cite{halmos}.

We will consider formulas $\theta_u$  in variables $(x_i: i \in u)$, and let $\theta =  \bigwedge_{|u|=2} \theta_u$. 
We interpret $\theta$ on the one hand as a clopen subset of $X_\Ups$, on the other hand as 
   a clopen subset of $\Pi_u X_u$, namely $\Pi_u \theta_u$.

\<{lem}  \lbl{a1.5}  Let $(q_1,q_2,q_3) \in \Pi_i X_i$ be a random triple.  Let $q_{23}=tp(a_2a_3/M)$
where $tp(a_3/M(a_2))$ does not divide over $M$, and $tp(a_i/M)=q_i$ ($i=1,2$).    Let $\theta_{1j}(x_1,x_j)$
be a formula of positive measure for $X_{1j}(q_1,q_j)$ (the space $X_{1j}$ with measure $\mu$ conditioned
on $(q_1,q_j)$.)   Then $\theta_{12}(x_1,x_2) \wedge \theta_{13}(x_1,x_3) \union q_{23}$ is consistent.  
In fact for $(a_2,a_3) \models q_{23}$, $\theta_{12}(x_1,a_2) \wedge \theta_{13}(x_1,a_3)$ has positive $\mu_1$-measure.
\>{lem}

\prf  Choose $p_{12} \in \theta_{12}$, random in $X_{12}(q_1,q_2)$ over $(q_1,q_2,q_3)$.  
Note that $p_{12}$ extends $q_1,q_2$.
Since $q_2$ is random over $(q_1,q_3)$, $p_{12}$ is random in $X_{12}(q_1)$ over $(q_1,q_3)$,
and in $X_{12}$ over $(q_3)$.  Hence   $(q_3,p_{12})$ are random in $X_3 \times X_{12}$.

Choose $p_{13} \in \theta_{13}$, random in $X_{13}(q_1,q_3)$ over $(p_{12},q_3)$.  Again $p_{13}$
extends $q_1,q_3$.   And (as $q_3$ is random over $(p_{12})$ in $X_3$),
 $p_{13}$ is random in $X_{13}(q_1)$ over $(p_{12})$, so 
 $(p_{12},p_{13})$ is random in $X_{12}(q_1) \times X_{13}(q_1)$ over $(q_1)$.  Now the product measure
 on  $X_{12}(q_1) \times X_{13}(q_1)$ coincides with the pushforward measure from $X_{123}(q_1)$.  (This is best seen
``over $q_1$".)  
 So by choosing $p_{123}$ at random in $X_{123}(p_{12},p_{13})$ (with the conditional measure),
 we find $p_{123}$ containing $p_{12},p_{13}$ and random.   Let $p_{23} $ be the restriction of $p_{123}$
 to the $2,3$-variables.   Let $(b_2,b_3) \models p_{23}$.  Note that $p_{23}$ is random in $X_{23}$,    so  
 $tp(b_3/M(b_2))$ does not divide over $M$.
 
Now $\theta_{12}(x_1,b_2) \wedge \theta_{13}(x_1,b_3)$ has positive $\mu_1$-measure (otherwise $p_{123}$ could not be random.)
By \thmref{3replacement}, $\theta_{12}(x_1,a_2) \wedge \theta_{13}(x_1,a_3)$  has positive $\mu_1$-measure too.   
\eprf

   \<{thm}    \lbl{ind}  Assume $L(M)$ is countable.  Let $\Ups=\{1,2,3\}$.  For $i \in \Ups$ let $\mu_i$ be an $M$-definable measure on $X_i = X_{x_i}$, and assume the $\mu_i$ commute.  For $u \subseteq \Ups$, $|u|=2$, let $\mu_u$ be   the tensor product measure  on $X_u$.    Then there
 exist measure-one Borel subsets $\Omega_u \subset X_u$ and $\Omega \subset X_1 \times X_2 \times X_3$ 
  with  the following amalgamation property.   Assume $q_u \in \Omega_u$, $(q_1,q_2,q_3) \in \Omega$, 
 $q_{u} | i = q_i $ for $i \in u$.    Then there exists $q \in X_{\Ups}$, $q | u = q_{u}$.   
 
 In fact,    we can take $\Omega_{23}$ to be the set of all $tp(bc)$ such that $tp(b/c)$ does not divide over $M$. \>{thm}
 
    \prf    It suffices to show that if $(q_1,q_2,q_3)$ is random,   in $X_1 \times X_2 \times X_3$, 
   $q_u$ is random in $X_u$ for $|u|=2$, and $q_i \subset q_u$ for $i \in u$, then there exists $q \in X_{\Ups}$, $q | u = q_{u}$.  Fix  such $q_i,q_u$.  By compactness, it suffices to show for  any given triple of formulas $\theta_u \in q_u$ that $\theta = \bigwedge_u \theta_u$ is consistent.  Fix such $\theta_u$.  Since $q_{1j}$ is random in $X_{1j}$, it is random in $X_{1j}(q_1,q_j)$
 over $(q_1,q_j)$.  Hence $\theta_{1j}$ has positive measure in $X_{1j}(q_1,q_j)$.  By \lemref{a1.5}, even
 $\theta_{12}(x_1,x_2) \wedge \theta_{13}(x_1,x_3) \union q_{23}$ is consistent.  \eprf

Note that since $\Omega_w$ has measure $1$, for a random choice of $q_{i} \in X_{i}$ $(i=1,2,3$), one expects the existence of $q_w \in S_w$ ($w \subset \{1,2,3\}$ with $|w|=2$)
with $q_i \subseteq q_w$ when $i \in w$.  The (obviously necessary) hypothesis of compatibility on the $q_w$ is therefore frequently attained. 
 
   Thanks to Pierre Simon for his comments on this.   
This result admits a more precise numerical  version, or alternatively a formulation using ideals, and a higher dimensional generalization;  this and related issues will be taken up elsewhere.

\ssec{NIP and de Finetti} \lbl{NIPdF}
 
   \<{example} \lbl{nip}  \rm Let $\mu$ be an $A$-definable Keisler measure in a NIP theory, cf. \cite{nip2}.   Let $\phi(x,y), \phi'(x,y)$ be formulas.
   For any real $\a$, let $R_\a(a,b)$ denote the relation:  $\mu( \phi(x,a) \meet \phi'(x,b)) < \alpha$.  Then
   $R_\a$ is  equational.  This uses the fact that for an indiscernible sequence $(c_j)$  over $A$ we have $\mu(\psi(x,c_i) \meet \psi(x,c_j) ) = \mu(\psi(x,c_i))$, applied to $c=(a,b), \psi(x,c)=\phi(x,a) \wedge \phi'(x,b)$.   \>{example}   
 
   When  $\alpha>0$, the relation 
$\mu( \phi(x,a) \meet \phi'(x,b)) = \alpha$ need not be equational, as one sees for instance by taking $\phi=\phi'$ and an indiscernible sequence $(a_i,b_i)$ with $a_i=b_i$.  

However, in  {\em any} theory, we have:

\<{prop}  \lbl{stable2}For any invariant measure $\nu$,  the relation $\nu( \phi(x,a) \meet \psi(x,b)) = \alpha$ is stable.  In other words, when $(a_i,b_i)$ is an indiscernible
sequence of pairs, the function $(i,j) \mapsto \nu( \phi(x,a_i) \meet \psi(x,b_j))$ is symmetric in $i,j$.    \>{prop}

 It follows that for any subset $Y$ of $[0,1]$, the relation:    $\nu( \phi(x,a) \meet \psi(x,b)) \in Y$ is stable.    
 
 The proof is related to a classical theorem of de Finetti, classifying the so called {\em exchangeable sequences} of random variables, i.e. sequences such that the action of the symmetric group does not change joint distributions.    This was subsequently
generalized by \cite{hewitt-savage}, \cite{krauss-scott}-\cite{krauss}, and in a different direction by Aldous and Hoover, see \cite{kallenberg}.  Thanks to Benjy Weiss  for telling me about this theory.  Though the assumption is classically stated as symmetry, indiscernibility suffices for the arguments;
the proof below is essentially a subset of the one in \cite{krauss} (in turn a modification of \cite{hewitt-savage}).   The higher dimensional case will be
considered elsewhere.

\prf of \propref{stable2}.  We show more generally that if $(a_i: i \in \Nn)$ is an indiscernible sequence, and $\psi_1,\ldots,\psi_k$ any
formulas, then 
$\nu(\psi_1(x,a_1) \wedge \cdots \wedge \psi_n(x,a_k))$ is invariant under the action of the symmetric group on $\{a_1,\ldots,a_k\}$,
i.e. 
$$\mu( \psi_1(x,a_1) \meet  \cdots \meet \psi_k(x,a_k)) =  \mu( \psi_1(x,a_{\si 1}) \meet  \cdots \meet \psi_k(x,a_{\si k}))$$ for any $\si \in Sym(k) $

 Let $B(\Nn)$ be the Boolean algebra generated by the formulas $\psi_i(x,a_j)$ for $i\leq k, j  \in \Nn$.   Let $S=S(\Nn)$ be the Stone space of $B(\Nn)$.  
Let $\mathcal{M}$ be the space of  countably additive regular Borel probability measures on $S(\Nn)$.   For a finite $J \subset \Nn$, let $B(J)$ be the subalgebra
generated by the $\psi_i(x,a_j)$ with $j \in J$, $S(J)$ the Stone space, and for $\mu \in \mathcal{M}$, let $\mu|J$ be the induced measure,
i.e. the pushforward of $\mu$ under the restriction map.  
Let $ \mathcal{M}_{ind}$ be
the subset of {\em indiscernible} measures, i.e. measures $\mu$ on $S$ such that for any finite $J_1,J_2 \subset \Nn$ of the same size,
with  order preserving bijection $j: J_1 \to J_2$, the induced map $j:  B(J_1) \to B(J_2)$ is measure-preserving, i.e. $j_*(\mu|J_1) = \mu|J_2$.

Let $\mathcal{M}_{sym}$ be the apparently smaller subset of {\em symmetric} (or {\em exchangeable}) measures, where we demand
that $j_*(\mu|J_1) = \mu|J_2$ for {\em any} bijection $j: J_1 \to J_2$.   

\claim{1}    $\mathcal{M}_{sym} = \mathcal{M}_{ind}$

 To prove the claim, note that both sets are convex and weak-* closed  subsets of the unit  ball of  $\mathcal{M}$.  Hence
 by Krein-Milman (cf. e.g. \cite{zimmer}), to show equality it suffices to prove that any extreme point of $\mathcal{M}_{ind}$ is in $\mathcal{M}_{sym}$.
 So assume $\mu$ is an extreme point of  $\mathcal{M}_{ind}$.  Now Claim 1 follows from:

\claim{2}  When $\mu \in \mathcal{M}_{ind}$ is extreme, we have independence:   $\mu(\phi_1(x,a_1) \wedge \cdots \wedge \phi_n(x,a_n)) = \Pi_{i=1}^n \mu \phi_i(x,a_i)$, for any  $\phi_i(x,a_i) \in B(\{a_i\})$.   

Let $\a= \mu(\phi_1(x,a_1))$.   
If $\a=0$, then $\mu(\phi_1(x,a_j))=0$ for any $j$, and both sides of the equation vanish.  If $\a=1$, then $\phi_1$ may be deleted on both sides,
and the claim follows by induction on $k$.   Assume therefore that $0<\a<1$.
Let $\mu'$ be obtained from $\mu$ by conditioning on $\phi_1(x,a_1)$, and shifting indices:
$$\mu'(\theta(x,a_1,\ldots,a_m)) = \mu(\theta(x,a_2,\ldots,a_{m+1}) \wedge \phi_1(x,a_1)) /\a $$
Similarly, let  $\mu''$ be obtained from $\mu$ by  conditioning on $\neg \phi_1(x,a_1)$, and shifting indices.  Then $\mu= \a \mu'+(1-\a) \mu''$;
and $\mu',\mu'' \in {\mathcal M}_{ind}$.  As $\mu$ is extreme, we have $\mu=\mu'$.  This means:
 $$\mu(\phi_1(x,a_1) \wedge \theta(x,a_2,\ldots,a_m)) = \mu(\phi_1(x,a_1)) \mu(\theta(x,a_2,\ldots,a_m))$$
Here $m,\theta$ are arbitrary.    Claim (2)  follows by induction on $m$, letting
 $ \theta(x,a_2,\ldots,a_n)  =  \phi_2(x,a_2)    \wedge \cdots \wedge \phi_n(x,a_n)  $.
 
Claim (1) follows easily:  the right hand side of the formula of Claim 2 is clearly symmetric. 
Any formula in $B(\Nn)$ is a disjoint union  of set-theoretic differences of conjunctions as considered in Claim (2)  The measure of the 
  difference  of two such expressions can be computed using the 
inclusion-exclusion formula, and of disjoint unions by additivity.

 Finally note  that if $\nu$ is an invariant measure, indiscernibility of the $(a_i)$ implies indiscernibility of $\nu | B(\Nn)$; hence the proposition
 follows from   Claim 1.
 
    \eprf
 
\>{section} 

 \<{section}{The stabilizer}

 Let $G$ be a group, $X$ a subset, defined over some model $M_0$.   
    Let $\tG$ be the subgroup of $G$ generated
 by $X$ (cf. \cite{hpp}, \S 7.)   
 By a definable subset of $\tG$, we mean a definable subset of $(X \union X \inv)^{\leq n }$
 for some $n$.   
 A subset $Y$ of   $\tG$ is {\em locally definable} if $Y \meet D$
is definable for every definable subset $D$ of $\tG$.

\<{remark} \rm  In sections 3 and 4 
we will never use $G$, only $\tG$.   It is thus natural to use a many-sorted reduct, whose   universes consist of the
sets $(X \union X \inv)^{\leq n}$, with the inclusion maps and multiplication maps between them, and a distinguished
predicate for $X$.   We will  speak of the inclusion maps as if they were actual inclusions.  

Going further, we can note that we actually use only a bounded number of multiplications.  In this section we will use only elements of 
  $(XX \inv)^3$, and will only use associativity for products of at most twelve elements of $X$ and their inverses.  
 
 Hence the results of this section are valid for structures $(X,X',G)$ with $X \subset X' \subset G$, with a binary map $m: (X')^2 \to G$ and
 an inversion map $\inv: X' \to X'$,  such that products of up to twelve elements
 of $X \union X \inv$ are defined, and independent of order.     We will refer to this as a ``local group" situation (cf.   \cite{goldbring}).       In this case $\tG$-translation invariance for a measure is replaced by the condition that $\mu$ measures $X$, and  $\mu(Y) =\mu(Ya)$
    for $Y \subseteq X\inv X$ and $a \in X \inv X$.    To avoid too technical a language
    we will state the results using the Ind-definable group $\tG = \union_n (X \inv X)^n$,
   indicating occasionally how to restrict to $(X \inv X)^3$.   The reader is welcome to ignore
   these refinements at a first reading.
         \>{remark}  
  
      An $\inft$-definable subset of $(X \inv X)^3$  closed under $m$ and $\inv$ will be called an $\inft$-definable subgroup of $\tG$ (though in the local setting there is a priori no group of which it is a substructure).  The main case is that of countable intersections;
 in this case one can write $H=\meet_{n \in \Nn} H_n$,  with $H_n$ definable,   $H_n=H_n \inv$,  and $H_nH_n \subseteq H_{n+1}$.  It is easy to see that
 any $\inft$-definable subgroup is an intersection of such countably- $\inft$-definable
 subgroups.     $\tG/H$ is {\em bounded} if for any definable subset $Y$ of $\tG$, the $\inft$-definable equivalence relation: $y \inv y' \in H$ has boundedly many classes   in the sense of \secref{basic}.   
  (equivalently, if $\tG,H$ are defined over $M_0$, the cardinality of  $\tG(N)/H(N)$ remains bounded when $N$ runs over all elementary extensions of $M_{0}$.)

Let $\tG$ be generated by the definable set $X$.  
Let $\mu$ be an ideal on $\tG$, invariant under right translations by elements of $X$ (i.e. $Z \in \mu $ iff $Z b \in \mu$, $b \in X$);
equivalently, $\mu$ is invariant under right translations by elements of $\tG$.  Assume   $\mu(X)>0$.
 Recall that a partial type $Q$ is called {\em wide}   if $\mu(Q')>0$ for any definable $Q' \supseteq Q$.
 
A definable subset $Z$ of $\tG$ is called {\em right generic} if finitely many
right translates of $Z$ cover any given definable subset of $\tG$.  If $Z$ is right generic
then clearly $\mu(Z)>0$.  In the converse direction we have  the observation, due to Ruzsa in the combinatorics literature, and Newelski
in the model theory literature,  that if  
$\mu(Z)>0$ then  $Z \inv Z$ is right generic.  We state this as a lemma for later reference.

\<{lem}\lbl{ruzsa} Let $\mu$ be an ideal on $\tG=<X>$, invariant under right translations by elements of $X$, and with  $\mu(X)>0$.
If $Q$ is a wide partial type, then so is $Q \inv Q$.  If $Z$ is a definable set with $\mu(Z)>0$, then   $Z \inv Z$ is right generic.
\>{lem}  

\prf  The statement for partial types follows by definition from the same statement for definable sets; so consider a definable
set $Z$ with $\mu(Z)>0$.  We have to show that $Z \inv Z$ is right generic, and wide.  

Let $X_n=(XX\inv)^n$;  say $Z \subseteq (XX\inv)^n$, and let $\{Za_i: i \in I\}$  be a maximal collection of 
pairwise disjoint subsets of $Z$, with $a_i \in X_n$.  We claim that $I$ is finite.  Otherwise, by the usual Ramsey/compactness argument on existence of indiscernibles, one can find indiscernible $(a_i: i \in \Nn)$ with $a_i \in X_n$ and  $Za_m \meet Za_{m'} = \emptyset$ for $m \neq m'$; 
by the S1 property, since $\mu(Za_i)>0$ for
each $i$ by right invariance, while $\mu(Za_i \meet Za_j)=0$ for $i \neq j$, $I$ must be finite.  
 If $a \in X_n$ then $Za \meet Za_i \neq \emptyset $ for some $i$;
so $a \in Z \inv Z a_i$.  This shows that $Z \inv Z$ is right-generic. 

 In particular,    $X \subseteq \union_{b \in B} X \inv X_b$, for some finite $B$; since $\mu(X)>0$ it follows that $\mu(X \inv X b)>0$ for some $b \in B$, so $\mu(X \inv X)>0$.

\eprf

In the local case, we say $Z \subseteq (X \inv X)^2$ is right-generic if finitely many
translates $Zb$ ($b \in X$) cover $X \inv X$.  Again if $Z \subseteq X \inv X$
has positive $I$-measure, then $Z \inv Z$ is right-generic.

\<{lem} \lbl{gen-equiv}[cf. \cite{nip2}]     Let $H$ be an $\inft$-definable subgroup  of $\tG$.  Then $\tG/H$ is bounded
 iff  every
definable set containing $H$ is right generic.  For any right invariant S1-ideal $\mu$ on $\tG$
this is also equivalent to:  $H$ is wide.

\>{lem}

 \prf  Consider $H = \meet H_n$ as above.  If each $H_n$ is generic, since $\tG$ is a countable union of 
 definable sets, there exists a countable set $C_n$ such that  $H_nC_n = \tG$.  Let $C=\union_n C_n$.
 Let $\pi: \tG \to \tG / H$ be the natural map.  Say that a sequence $u_n$ of elements of $C$ converges
 to $uH \in \tG/H$ if for each $m$, for all sufficiently large $n$, we have  $H_m u_n = H_m u$.  Then each
 sequence has at most one limit, and each point of $\tG/H$ is the limit of some sequence from $C$.  Hence
 the cardinality of $\tG/H$ is at most continuum.  (We will later define the "logic topology" on $\tG/H$; in this language we have just shown it is separable.)   
 
  Conversely if $\tG/H$ is bounded, let $X$ be a definable subset of $\tG$.   
The condition:  $H_{k+1}x \meet H_{k+1}y = \emptyset$ is a definable relation on $(x,y)$, since $H_{k+1}$ is definable.
Say $\tG/H$ is bounded by $\lambda$; then a fortiori there cannot be more than $\lambda$ distinct   $(a_i)$ with $H_{k+1} a_i$ disjoint.  
Compactness applies, so any such family is finite.  Let $(a_i)$ be   a maximal family $H_{k+1} a_i$ of disjoint cosets of $H_{k+1}$, with $a_i \in X$.  Then there are finitely many elements $a_i$ in the family,  and it follows that $X \subseteq \union_i H_{k+1} \inv H_{k+1} \inv a_i = \union_i H_ka_i$,  i.e. $H_k$ is right-generic.  

 Given a right invariant S1-ideal  $\mu$, if $H$ is wide then there can be no infinite family of disjoint
cosets of $H_{k+1}$, so as above $H_{k}$ is generic.  Conversely if  $H_{k}$ is generic
then $\mu(\union_j H_k b_j) > 0$ for some finite set $b_1,\ldots,b_l$, so
$\mu(H_k b_j) > 0$, and by right invariance $\mu(H_k)>0$.   \eprf
 
 If an $\inft$-$A$-definable   subgroup of bounded index exists, then there is a   {\em minimal}  one;   it is denoted $\tG^{00}_A$.
 For a discussion of the dependence on   $A$, see \cite{hpp}. 
 
 \<{lem}  \lbl{normal0}  $\tG^{00}_A$ is normal in $\tG$.  \>{lem}
 
 \prf   Let $H=\tG^{00}_A$.  Then $H$  has boundedly many
 $\tG$-conjugates; their intersection is an $\inft$-definable normal subgroup $N$
 of $\tG$.  On the face of it the definition of $N$ requires additional parameters; but
 $N$ is $Aut(\Uu/A)$-invariant, and in 
 general if an $\inft$-definable set is invariant under $Aut(\Uu/A)$
 then it is an infinite intersection of $A$-definable sets.  \eprf

 \<{thm}\lbl{stabilizer}  \lbl{normal} Let $M$ be a model, $\mu$ an $M$-invariant S1-ideal on definable subsets of $\tG$, 
 invariant under (left or right)  translations by elements of $\tG$.    Let $q$ be a wide type over $M$ (contained in $\tG$.)  Assume:
 
 (F)  There exist two realizations $a,b$ of $q$ such that $tp(b/Ma)$ does not fork over $M$ and
  $tp(a/Mb)$ does not fork over $M$.
 
     Then  there exists a wide, $\inft$-definable over $M$ subgroup $S$ of $G$.   We have 
     $S=(q \inv q)^2$; the set $q q \inv q$ is a coset of $S$.  Moreover, $S$ is normal in $\tG$, 
      and $S \m q \inv q$ is contained in a union of non-wide $M$-definable sets. 
               \>{thm}

Some remarks before turning to the proof.

 \<{enumerate} 
 
 \item It follows from the statement of the theorem that $S$ can have no proper $M$-$\inft$-definable subgroups of bounded index.
For suppose such a subgroup $T$ exists.   Then $q$ is contained in a bounded union of cosets of $T$.  Being a complete type over a model, it is contained in a single coset.  But then $q \inv q q\inv$, a coset of $S$, is contained in a coset of $T$; so $S=T$.
 
 \item  
  The statement about $S\m q \inv q$  can be read to say that a random element of $S$ lies in $q \inv q$; for instance when $M$ is countable, and $\mu$ is the  ideal of
 definable measure zero sets for some finitely additive measure $\mu$ on the Boolean
 algebra of $M$-definable sets, $\mu$ extends to a Borel measure on the space of types,
 and almost all types of elements of $S$ lie in $q \inv q$. 
 
 \item When $\mu$ is the zero-ideal of a measure, note that translation invariance is assumed of
 the ideal, not of the measure.  In particular, regardless of unimodularity, this assumption is true
 for  Haar measures on a locally compact group.
 
 \item (Weakening of left invariance.)   Most of the proof is devoted to showing that
 $S=(q \inv q)^2$ is a   subgroup of $\tG$, and $q q \inv q$ is a coset of $S$.  For this,
 left-translation invariance can be replaced with existence of an $f$-generic extension of  $q$, 
 in the sense of
 \cite{nip2}, i.e. the existence of an  $M$-invariant ideal $J$ containing  
 the forking ideal, and with $q$ wide for $J$.  We will use such a $J$ in Claims 3' and 5' (without assuming
 that $\mu=J$.)  The statement is essentially that left generics do not fork, and involves $\mu$ but not $J$.
       
 The word "wide" will refer to $\mu$ unless explicitly qualified.

Normality of $S$  will also follow under these assumptions, but we do not obtain the final statement about
$S \m q\inv q$ in this case.

 \item  In place of any form of  left translation invariance,  we could use  a stronger Fubini-type assumption on $\mu$ itself. 
(In Claim 3' of the theorem, we need to find $(c_1,c_2,a)$ with $tp(c_i/M)$ specified, $c_i \in q \inv q$,
and with $tp(a/M(c_1,c_2))$ wide.  Given a version of Fubini we can achieve this by choosing $a$ first, then 
$c_1,c_2$.)

 \item (Locality). 
   Inspection of the proof will show that for all assertions except the normality of $S$,  we only use $\mu$ (as an S1 ideal)  on definable  subsets of $X X \inv X$.  To show normality $S$,  we also require $X  a X \inv$, where 
   $a \in X$ or $a \in X \inv$.       
         Moreover the group structure is used only up to $(X \inv X)^3 $. 
  This is explicitly so everywhere except in Claim 5.  There, note that $qc \subseteq  X X \inv X$.  
  Hence $qc \meet Y \subseteq  X X \inv X$ for any set $Y$, and it makes sense to say
   that this intersection is wide.     In the proof, by the time we use $qab_1$, we know that $ab_1$ is in $q \inv q$.

\noindent It is also possible to combine (4) and (6); see
  \exref{locality+}.  

 \item The theorem implies that $S \subseteq X \inv X X \inv X$; or that for a  appropriate translate $Y=a \inv X$, we have $S = Y  Y \inv Y$.   Example 6.1.10 of \cite{CH} shows
 that this cannot be improved to $S \subseteq X \inv X$. 
 
 \item An easy L\"owenheim-Skolem argument shows that the theorem reduces to the case where   the language is countable, and $M$ is countable.
 
 \item  We show in fact that  $S \m St_0(q)$ is contained in a union of non-wide $M$-definable sets, where $St_0(q)   = \{s: qs \meet q \hbox{ is wide } \}$.
If $s \in S$ is arbitrary now, and $tp(s'/M(s))$ is wide, then $tp(s's/M)$ is wide, 
so $s', s's \in St_0(q)$.  Hence $s = (s') \inv (s's) \in St_0(q) \inv St_0(q)  =St_0(q)^2$.

\item The assumption that $M$ is a model, rather than just a substructure of the universal domain, is used via the consequence that any type extends to an invariant type; thus  \thmref{3replacement}  applies to any type $tp(a)$.    See \remref{3replacementr}.

\item  The  proof uses both the nonforking ideal and the ideal of wide sets with respect to $\mu$ (or $J$).   The former allows \thmref{3replacement} to be used
for an arbitrary type, since any type has a nonforking extension.  On the other hand the ideal of wide sets, unlike the nonforking ideal, enjoys translation invariance.

  \>{enumerate}
      
 \prf[Proof of \thmref{stabilizer}]   
  We also write $q$ to denote $\{a: tp(a/M) = q\}$;
 and $q \inv = \{a \inv: tp(a/M) = q\}$.
  
 Given two subsets  $X,Y$ of $\tG$, let  
 $$X \times_{nf} Y = \{(a,b) \in X \times Y: tp(b/M(a)) \hbox{ does not fork over } M\}$$
  
 Let $Q=  \{a \inv b:  (a,b) \in q \times_{nf} q \}$.   Let $J$ be as in Remark (4) (or just set $J=\mu$ for the basic statement of the theorem),   and set $Q' = \{a \inv b:  a,b \in q, tp(b/Ma) \hbox{ is $J$-wide} \}$.

  Note $q q\inv$ is obviously wide by right-invariance,  and similarly $q \inv q$ is wide assuming left-invariance.
 If we wish to avoid  the left invariance assumption, but are willing to use $\mu$ on $X^2$ instead, then wideness of $q \inv q$
 follows from \lemref{ruzsa}.   
 
Throughout this proof, we will use the fact (\lemref{stable}) that wideness of   $q x \meet q y \inv$   is a stable relation between $x$ and $y$.   By \lemref{stat}, or \thmref{3replacement}, for any two types $p_1,p_2$, this relation holds for one pair $(a_1,a_2) \in p_1 \times_{nf} p_2$
iff it holds for all pairs iff it holds for one or all pairs $(a_2,a_1)$ in $p_2 \times_{nf} p_1$.

\claim{1}  $q\inv q \subseteq QQ$.

\prf  Let  $a,b \in q$.      Using (F), find $c \models q$ be such that
$tp(a/Mc)$ does not fork over $M$, and $tp(c/Ma)$ does not fork over $M$.  By extending
$tp(c/Ma)$ to a type over $M(a,b)$ and realizing this type, we may assume $tp(c/Mab)$ does not fork over $M$.  
So we have $(b,c) \in q \times_{nf} q$, and    $(c,a) \in q \times_{nf} q$.   So $b \inv c, c \inv a \in Q$,
hence $b \inv a \in QQ$.  \eprf

\claim{2}   For  all $(a,b) \in q \times_{nf} q$, $qa \inv \meet q b \inv$ is wide.

 \prf  
 By \thmref{3replacement}, it suffices to show that for {\em some} $(a,b) \in q \times_{nf} q$, $qa \inv \meet qb \inv$ is wide.
  Let $a_{1},a_{2},\ldots$ be an $M$- indiscernible sequence of   elements of $q$,    such that $tp(a_i/ A \union \{a_j: j<i \})$ does not fork over $M$.
  Then $(a_{i},a_{j}) \in q \times_{nf}  q$ for any $i < j$. It suffices to show that $q a_{1} \inv \meet q a_{2} \inv$ is wide;
  by compactness, for any definable set $D$ containing $q$, it suffices to show that $\mu(Da_{1} \inv \meet Da_{2} \inv) >0$.  This is clear since $\mu$ is an S1-ideal, and by right-invariance, $\mu(Da_i \inv)>0$.
   \eprf

\claim{3'}  For all $(c_{1},c_{2}) \in (q \inv q) \times_{nf} Q'$, 
$qc_{1} \inv \meet qc_{2} \inv$ is wide.
 
\prf Let $p_{i} =tp(c_{i}/M)$.  As in Claim 2,  it suffices to see that $q c_{1}\inv \meet q c_{2}\inv$ is wide for some $(c_{1},c_{2}) \in p_{1} \times_{nf} p_{2}$.  Let $a_0 \models q$.  Then
there exists $a_1 \in q$ with $tp(a_0 \inv a_1 /M) = p_1$.  Since $c_2 \in Q'$, there
exists  $a_{2}'$ such that $r=tp(a_2'/M(a_0))$  is $J$-wide and   $tp(a_{0} \inv a_{2}' /M)=p_{2}$;     extend $r$ to a $J$-wide type $r'$ over $M(a_{0},a_{1})$, and let $a_2 \models r'$.    We thus have  
$(a_{0},a_{1},a_{2}) \in (q \times q) \times_{nf} q$,  with $tp(a_{0} \inv a_{i} /M) = p_{i}$ for $i=1,2$.  Note also, using left invariance of $J$,  that    $tp(a_{0} \inv a_{2}/M(a_{0},a_{1}))$ is $J$-wide,
hence so is   $tp(a_{0} \inv a_{2}/M(a_{0} \inv a_{1}))$, so it  does not fork over $M$.
 
 By Claim 2  
 we have $q a_{1} \inv \meet q a_{2} \inv$ wide.
  By the right invariance of $\mu$,  $q a_{1}\inv a_{0} \meet q a_{2}\inv a_{0}$ is wide.
\eprf

 \claim{3}  For all $(c,d) \in (q \inv q) \times_{nf} Q$, 
$q{c} \inv \meet q{d} \inv$ is wide.

\prf Let ${d} = a \inv b$, with $tp(b/M(a))$ wide for the forking ideal over $M$.
We have to show that $q{c} \inv \meet q b \inv a$ is wide.  By \thmref{3replacement}, it suffices
to show this for {\em one} instance $({c},b,a)$ with $tp(b,a)$ specified and such that $tp(b,a / M({c}))$
does not divide over $M$.  We may thus take $tp(a/M(c))$ to be a nonforking extension of $q=tp(a/M)$,
and $tp(b/M(a,c))$ to be a non-forking over $M$ extension of $tp(b/M(a))$.  The latter is possible
using the assumption that $tp(b/M(a))$ does not fork over $M$.  
 
By right-invariance, we need to show that  $q {c} \inv a \inv \meet q b \inv$ is wide.  We apply 
\thmref{3replacement} to the pair $(a,b)$ (viewed as a single tuple) and $c$.  So it suffices to show that $q{c} \inv a \inv \meet q (b') \inv$ is wide, where
$tp(b/M)=tp(b'/M)$ and $tp(b'/M(a,{c}))$ is $J$-wide.  By left-invariance of $J$, the type
$tp(a \inv b' /M(a,{c}))$ is $J$-wide, and hence $tp(a \inv b' / M({c}))$ is $J$-wide; so $tp(a \inv b' /M({c}))$
does not fork over $M$.  Also $tp(b' / M(a))$ is $J$-wide, so $a \inv b' \in Q'$.  By Claim 3', $q {c} \inv \meet q (a \inv b') \inv$
is wide.  By right invariance, $q {c} \inv a \inv \meet q (b') \inv$ is wide, as required.     
\eprf

\claim{4} Let $(b,a) \in Q \times_{nf} q \inv q$.   
Then $ab \in q \inv q$.    In fact $q a \meet q b \inv$ is wide.

\prf We have $a \inv \in q \inv q$.   Since $M$ is a model, $tp(a \inv/M)$ extends to a
global type $r$  finitely satisfiable type in $M$; so $r$ is   $M$-invariant.   Use \lemref{stat} (1), and Claim (3) to conclude that $q c \inv \meet q b \inv $
is wide if $c \models r |M(b)$.  Now $tp(c/M(b))$ does not divide over $M$,
so by \thmref{3replacement}, since $tp(a \inv/M(b))$ does not divide over $M$ either,
$q a \meet q b \inv$ is wide.   In particular, for some $d,e \in q$ we have $da = e b\inv$.
So $ab = d \inv e   \in q \inv q$. 
\eprf

\claim{5}  Let $a \in q \inv q$, $b_1,\ldots,b_n \in Q$ and assume
$tp(a/M(b_1,\ldots,b_n))$ is wide.  Then $ab_1\cdots b_n \in q \inv q$.   In fact
$q a \meet q (b_1 \cdots b_n) \inv$ is wide.

\prf  Since $tp(a/Mb_1)$ is wide, it does not fork over $M$ (\lemref{fork}).  Hence
by Claim 4 we have $ab_1 \in q \inv q$.  By right-invariance of $\mu$,
$tp(ab_1/ M(b_1,\ldots,b_n))$ is wide, and in particular
$tp(ab_1/ M(b_2,\ldots,b_n))$  is wide.  By induction,
$q ab_1 \meet q (b_2 \cdots b_n) \inv$ is wide.  Multiplying by $b_1 \inv$ on the right,
$q a \meet q (b_1 b_2 \cdots b_n) \inv$ is wide. Hence as in Claim 4, $ab_1\cdots b_n \in q \inv q$.
\eprf

In view of \thmref{3replacement}, Claim 5 is also valid assuming $tp(a/M)$ is wide, and $tp(a/M(b_1,\ldots,b_n))$ does
not fork over $M$.    To show that $q q \inv q$ is a coset, we will later need a variant of Claim 5, proved in the same way:

\claim{5'}   Let $a \in q \inv q$, $b_1,\ldots,b_n \in Q$ and assume
$tp(a \inv/ M(b_1,\ldots,b_n))$ is $J$-wide.  Then $ab_1\cdots b_n \in q \inv q$.   In fact
$q a \meet q (b_1 \cdots b_n) \inv$ is wide.

\prf  Since $tp(a \inv /Mb_1)$ is $J$-wide, it does not fork over $M$, and so $tp(a/Mb_1)$ does not fork over $M$.    Hence
by Claim 4 we have $ab_1 \in q \inv q$.  By left-invariance of $J$,
$tp((ab_1) \inv / M(b_1,\ldots,b_n))$ is $J$-wide, and in particular
$tp((ab_1) \inv/ M(b_2,\ldots,b_n))$  is $J$-wide.  By induction,
$q ab_1 \meet q (b_2 \cdots b_n) \inv$ is wide.  Multiplying by $b_1 \inv$ on the right,
$q a \meet q (b_1 b_2 \cdots b_n) \inv$ is wide. Hence as in Claim 4, $ab_1\cdots b_n \in q \inv q$.
\eprf

\claim{6}  $Q^n \subset q\inv q q \inv q $.

\prf  Let $b_1,\ldots,b_n \in Q$.  Let $a \in q \inv q$ with 
$tp(a/M(b_1,\ldots,b_n))$  wide.  Then $ab_1\cdots b_n \in q \inv q$,
so $b_1\cdots b_n = a \inv (ab_1\cdots b_n) \in q\inv q q \inv q $.  
\eprf

 It follows from Claim 1 that 
  $Q$ and $q \inv q$ generate the same subsemigroup, which is hence a group $S$.
   By Claim (6), this  group is in fact equal to the $\inft$-definable set
$ q\inv q q \inv q $.

 Since $q\inv q \subseteq S$,
we have $q \subseteq bS$ for any $b \in q$, and so 
 $q q \inv q \subseteq bS$.   Conversely, choose $b \in q$.  
Any element $x$ of $bS $   can be written $x=ba_1 \cdots a_4 $ with $a_i \in Q$.
Let $d \in q$ be such that $tp(d/M(a_1,\ldots,a_4,b))$ is $J$-wide.  Let $e=d \inv b$.  
Then $tp(e \inv /M(a_1,\cdots,a_4,b))$ and hence $tp(e \inv /M(a_1,\cdots,a_4))$ are $J$-wide.
By Claim 5' we have $e a_1  \cdots a_4   \in q \inv q$.    
So  $x = b a_1  \cdots a_4     \in d q \inv q \subset q q\inv q$.
Thus  $q q \inv q = bS$

We know that $S$ is an $\inft$-definable group over $M$.  I claim any 
  $\inft$-definable over $M$ subgroup  of $S$ of bounded index  must be equal to $S$.   For let $T$ be such a subgroup.  We have $q \inv q \subseteq S$, so $q \subseteq a S$ for any $a \in q$.  Thus $q$ 
 is contained in a left translate $R$ of $S$; we have $R=qS$ so $R$ is defined over $M$.   Now $T$ acts
on $R$ on the right; the equivalence relation induced is $\inft$-definable over $M$ with boundedly many orbits.
Since $q$ is a complete type over $M$, it has an $Aut(\Uu/M)$- invariant extension to $\Uu$;
this extension must pick a specific $T$-orbit $cT$, which is hence $\inft$-definable over $M$; by completeness again, 
as the realizations of $q$ intersect $cT$, 
$q$ is contained in $cT$.     But then $q \inv q \subseteq T$; so $S \subseteq T$.

We know at this point that
$S$ has no proper $\inft$-definable over $M$ subgroups of bounded index. 
 Let $r$ be a type of elements of $X \union X \inv$ over $M$.   There cannot exist an unbounded
 family of cosets $a_iS$ with $a_i \in r$, for then the sets $a_ibq$ would also be disjoint for any $b \in q \inv$,
 so for some definable $X'$ with $q \subset X' \subset X$ the sets $a_ib X'$ can be taken disjoint,
 contradicting the S1 property for $\mu$ within $r b X \subseteq (X \union X \inv)^3$.
  Thus $r$ is contained in boundedly many left cosets of $S$, hence (being a complete type over a model) in one; call it $C_r$.  
 So $C_r$ is $M$-definable, and hence the conjugate group   $S^r=C_r \inv S C_r$ is $M$-definable. 
 
  For any $c \in X \union X \inv \union \{1\}$, $r=tp(c)$,
the image of $qc$ in $G/S$ is bounded.  Otherwise there is a large collection of disjoint sets
of the form $a_i c S$, with $a_i \in q$.  Pick $b_0 \in q$; then $q \inv b_0 \subseteq S$; 
the sets $a_i c S b_0 \inv$ are also disjoint, hence so are the $a_i c q \inv$.  Thus there exists a 
definable $X' \subset X$ with $a_i c (X') \inv$ disjoint.   So the sets $X c \inv a_i \inv$ are disjoint,
and wide.  But this contradicts the S1 property within $Xc X \inv$.    Thus $qc/S$ is bounded.
It follows that $q$ is contained in boundedly many cosets of $cS c \inv=S^r$.  
So $q$ is contained in a single coset   $g S^r$.  It follows that $q \inv q \subseteq S^r$, so $S \subseteq S^r$.
Similarly $S \subseteq S^{r \inv}$, so $S^r \subseteq S$ and $S^r = S$.  This shows that $X \union X \inv$
normalizes $S$, i.e. $S$ is normal in $\tG$.

 At this point we begin using left invariance freely.

We argued above that $q \inv q$ is wide; in particular $S$ is wide.  
     $Q$ is also wide:  suppose otherwise.  So $Q \subseteq D$ for some definable $D$ with $\mu(D)=0$.
Let $a \in q$.  Then $a \inv q$ is wide.  So $a \inv q \m D$ is wide.  However $ q \m aD$ forks over $M$,
 since if $b \in q \m aD$ then $a \inv b \notin Q$ so   $tp(b/M(a))$ forks over $M$.  Thus
$D' \m aD$ lies in the forking ideal, for   some definable $D'$ containing $q$.  By \lemref{fork} we have
$\mu(D' \m aD)=0$; so $\mu(a \inv D' \m D)=0$.  It follows that $\mu(a \inv D')=0$ and $\mu(D')=0$,
contradicting the wideness of $q$. 

We finally show that $S$ is contained in $q \inv q$ up to a union of
non-wide definable sets.  Let $r$ be a wide type over $M$ extending $S$; we have to show that
$r \subseteq q \inv q$. Pick  $a_0 \in r$ and $c \in Q$ with $tp(c/M(a_0))$ wide.     As $a_0 \inv \in S$, we may
write $a_0=b_1 \cdot \cdots \cdot b_n$ with $b_i \in Q$; also as $c \in Q$ we have $c \in q \inv q$.  Thus Claim 5 applies
(with $c$ playing the role of $a$ of Claim 5); and we obtain that   	 $qc \meet q a_0$ is wide.   
Choose $b_0 \in r$ with  $tp(b_0 /M(c))$ wide.  In particular
 $tp (b_0/ M(c))$ does not fork over $M$.  
By stability of the relation and   \lemref{stat} (3), $q c \meet q b_0$ is wide too.
Thus $ b_0 c \inv \in  q \inv q$.  Now $tp({b_0}c\inv / M(c))$ is a right translate of $tp({b_0}/M(c))$, so it is wide.  
By Claim 5 (or 3), ${b_0}=({b_0} c\inv) c \in q \inv q$ (and $q b_0 c \inv  \meet qc \inv $ is wide;
so $qb_0 \meet q$ is wide.)
  So  $r \subseteq q \inv q$ as required.
In fact this shows that $r \subseteq St_0(q)$, in the notation of 
Remark 9.

     \eprf

 \<{cor} \lbl{stab1}  Let  $\mu$ be  an  invariant S1-ideal on definable subsets of $\tG$, 
 invariant under   translations by elements of $\tG$.    
     Then  there exists a model $M$ and a wide, $\inft$-definable over $M$ subgroup $S$ of $G$, with $\tG/S$ bounded.
     For an appropriate complete type $q$ over $M$ we have $S=(q \inv q)^2$, and the complement $S \m q \inv q$ is contained in a union of non-wide $M$-definable sets.     
     
  If $\mu$ satisfies the conditions of \lemref{sym} over a model $M_0$, or if 
 $\mu$ is $\bigvee$-definable over $M_0$,    then one can take $M=M_0$.          
     \>{cor}       
     
     \prf  \lemref{claim0} provides $M$ and an $M$-invariant global type $q^*$ such that if $q=q^*|M$,  $a \models q|M$ and $b \models q^* | M(a)$ then $tp(a/M(b))$ is $\mu$- wide.   This implies (F).      
     In case the assumptions of \lemref{sym} or \lemref{sym-d} hold, these lemmas provide   a type over $M_0$ with (F) and so \thmref{stabilizer}
     applies with $M=M_0$.  
     \eprf

    \<{example} \lbl{unbdd} \rm  Consider  the theory of divisible ordered Abelian groups $(G,+,<)$, or any   o-minimal expansion, and let $M$ be a model.  We have a two-valued definable measure $\mu$, assigning measure $0$ to any bounded definable set.  A two-valued invariant
 measure is always S1.  The measure $\mu$ is translation invariant.  Let $q_A$   be the set of all measure-one $M$-definable formulas over $A$, $q=q_M$.    If $a \models q_M$ and $b \models q_{M(a)}$, then $tp(a/Mb)$
 does not fork over $M$ since it is finitely satisfiable in $M$, and $tp(b/Ma)$ does not fork over $M$ since it extends
 to an $M$-invariant type.  Hence (F) of \thmref{stabilizer} holds.  We can take $\tG=G$, $X=\{x \in G: x > 0 \}$.
 The subgroup $S$ is then $G$.  Note that $q \inv$ is {\em not} wide in this example.     
 \>{example}

   Here is an example of the situation discussed in comments (4,6), 
 where the $\inft$-definable group $S$ is not normal.   

 \<{example}   \lbl{locality+} \rm  Consider the theory ACVF of algebraically closed valued fields, say of residue
 characteristic $0$; the field of Puiseux series over $\Cc$ is a model.  Alternatively,
 let $M$ be an ultraproduct of the $p$-adic fields $\Qq_p$. 
 Let $K$ denote the valued field, 
 $\Oo$ the valuation ring, $\Mm$ the maximal ideal.  
 Let $G$ be the semi-direct product
 of the additive group $G_a$ with the multiplicative group $G_m$.  So $G=TU$ where $T,U$ are Abelian 
 subgroups, $U =G_a$ normal, $T \cong G_m$.  Let $t \in K$ be an element of valuation $>0$, and let
 $g$ be the corresponding element of $G$, so that conjugation by $g$ acts on $U$ as multiplication by $t$.
  Let $U_0 = \{x \in K: \bigvee_{m \in \Nn} \val(x) \geq -m \val(g)  \}$.
 View $\Oo \leq U_0$ as subgroups of $U$.  
 Within $G$, let $X=g \Oo$.  The group $\tG$ generated by $X$ is $g^{\Zz} U_0$.  
 Let $p$ be a generic type of $\Oo$; it avoids any coset of $\Mm$ in $\Oo$.  
 Let $\mu$ be the right-invariant ideal generated by $g \Mm$, and $J$
 the left-invariant ideal generated by $g \Mm $.  These are not the same; 
 notably $\Oo$ is in $J$ but not in $\mu$.      $\mu$ is not S1, but it is so when restricted to $X=XX \inv X$.      Let $q= g p $. 
   As in Remarks 4 and 6, the proof of \thmref{stabilizer}  goes through to give a subgroup $S$, namely $\Oo$ (it is definable in this case).  But $\Oo$ is not a normal subgroup of $\tG$.  \>{example}

\<{defn}\lbl{near}
  We call $X$ a {\em near-subgroup} of $G$ if there exists an  invariant S1-ideal  $\mu$ on definable subsets of $(X \union X \inv)^3 $, with $\mu(X)>0$, and with $\mu(Y)=\mu(Y')$ whenever $Y,Y' \subseteq  X X \inv X$ and  $Y'=cY$ 
or $Y'=Yc$ for some $c$. 
   \>{defn}
   
 We will see in  \corref{3-to-n} that asking for $\mu$ defined on $\tG=\union_n (X \union X \inv)^n$ would result
 in the same definition; in later sections we will work with this stronger definition.

 \<{remark}  Lou Van den Dries has shown that  a weaker condition suffices:  $0 \in X$,
 and $\mu$ is defined on $X X \inv X$.  Moreover the element $c$  (which must by definition be in $(X \union X \inv)^6$)   can in fact be chosen so that all products are taken within $X X \inv X$.  This condition is essentially sharp, in view of \exref{locality+}.  See \cite{vddries}
 \>{remark} 
 
     When $X$ is finite, any right-invariant measure must be proportional to the counting measure.
     Asymptotically, when $(X,G)$ vary in some family, we have that every ultraproduct is a near-subgroup iff
     $|X X \inv X| / |X|$ is bounded in the family.

 The following   corollary of \thmref{stabilizer} is analogous to Lemma 3.4 of \cite{tao-nc};
 the point is that we do not assume a priori that $(X \inv X)^n$ has finite measure.    
 The Fubini-type assumption on the ideal is much weaker here, but the conclusion is purely qualitative.      We state the
 extension lemma for S1 ideals; a similar statement for measures is also valid, with a similar proof.
  
\<{cor} \lbl{3-to-n} Let $X$ be a near-subgroup of $G$.  Then for any $n$, $(X \inv X)^n$ is contained in a finite union of right translates of  $(X \inv X)^2$.  $\mu $ extends   to an invariant S1-ideal $\mu'$ on $\union_n (X \inv X)^n$;
$\mu'$ is the unique right-invariant ideal extending $\mu | (X \inv X)^2$.    
\>{cor} 

\prf  For this we may add parameters, and work over a model.  Let $\tG$ be the group generated by $X$.  
By \thmref{stabilizer} and Remark 6 to that theorem,  there exists a wide $\inft$-definable 
subgroup normal $S$ of $\tG$.  The proof also shows that $S \subseteq (X \inv X)^2$ and that the image of $X$ modulo $S$ has bounded  cardinality.   Hence $\tG/S$ is bounded, and in particular for any $n$,
$(X \inv X)^n$ is contained   in boundedly many   cosets of $S$, and hence in boundedly many right
translates of $(X \inv X)^2$.   By compactness, finitely many right cosets of 
$(X \inv X)^2$ suffice to cover $(X \inv X)^n$.  
 If $D$ is a definable subset of $\union_n (X \inv X)^n$, it follows that 
we can write $D = \union_i D_ib_i$ where $D_i \subseteq (X \inv X)^2$ and $b_i \in  (X \inv X)^{n+2}$. 

Define $\mu'$ to be the collection of all definable sets  $\union_i D_ib_i$, where $D_i$ is a definable subset of $(X \inv X)^2$,  $b_i \in  \union_n (X \inv X)^{n}$,
and $\mu(D_i)=0$.    This is clearly a right-invariant ideal whose restriction to definable subset of  $(X \inv X)^3$ is precisely $\mu$.  
        (if $Y$ is a definable subset of $(X \inv X)^3$ and $Y=\union_i D_ib_i$ as in the definition of $\mu'$, then by invariance we have  $\mu(D_ib_i)=0$        for each $i$, so $\mu(Y)=0$.)
        
         If $\mu''$ is any  right-invariant ideal $\mu'$ extending $\mu$, and $D$ a definable subset of $(X \inv X)^n$, write  $D = \union_i D_ib_i$ where $D_i \subseteq (X \inv X)^2$ and $b_i \in  (X \inv X)^{n+2}$. We have $\mu''(D)=0$ iff $\mu''(D_i)=0$ for each $i$ iff $\mu(D_i)=0$
for each $i$.   This shows that $\mu''=\mu'$.

To see that $\mu'$ is S1, it suffices to show for each $n$ that the restriction to $(X \inv X)^n$ is S1.  
 As above, write $(X \inv X)^n = \union _{j=1}^l D_j b_j$.  It is clear that any ideal on the union of finitely many sets must be S1, if
 the restriction to each of these sets is S1.  So it suffices to show that $\mu' | (D_jb_j)$ is S1 for each $j$.  But 
 $\mu' | (D_j b_j)$ is isomorphic, via translation by $b_j$, to $\mu | D_j$, which is S1.  
\eprf



 This kind of characterization incidentally makes some functorialities evident, that are not so directly from the definition
of a near-subgroup or an approximate subgroup; see \remref{towers} (0),(2).

Given  elements $a_1,\ldots,a_l$ and $b_1,\ldots,b_m$ of $G$, let
$A_i = \{x \inv a_i x : x \in X\}$ be the set of $X$-conjugates of $a_i$, 
and let $W_n(a_1,\ldots,a_l,b_1,\ldots,b_m)$ be the set of words of length $\leq n$
in $A_1 \union \ldots \union A_l \union \{b_1,\ldots,b_m\}$.  Let
$d(X;a_1,\ldots,a_l) $ be the smallest integer $n$ such that
  $X \subseteq W_n(a_1,\ldots,a_l; b_1,\ldots,b_l)  $ for some $b_1,\ldots,b_l \in X$; or $\infty$ of there is no such $n$.

\<{prop} \lbl{anti-nil2} For any $k,l,n \in \Nn$, for some $M,K \in \Nn$, the following holds:  

Let $G$ be a group, $X$ a finite subset.  Assume  $|X X \inv X| \leq k |X|$.
Also assume that   there exist $x_1,\ldots,x_M \in X$ such that:

(*)  for any
$1 \leq i_0< i_1< \cdots < i_l \leq M$,   
$d(X,  x_{i_0} \inv x_{i_1},\cdots,x_{i_0} \inv x_{i_l}) \leq n$.  

  Then  there exists a subgroup $S$ of $G$, $S \subseteq (X \inv X)^2$, such that $X$ is contained in $\leq K$ cosets of $S$.  \>{prop}

\prf   Fix $k,l,n$.  Suppose there are no such $M,K$; then there are groups $G_M$ 
and $X =X_M \subset G_M$
such that  there exist $x_1,\ldots,x_M$  with (*),   and there is no subgroup $S$ of $G$, $S \subseteq (X \inv X)^2$, such that $X$ is contained in $\leq M$ cosets of $S$.   
Consider $(G_M,X_M,\cdot)$ as a structure, and enrich it using the $Q_\a$ -quantifiers
for the normalized counting measure on $X_M$, as in \secref{measurequantifiers}.  
By compactness, there exists a  countably saturated group $G$ and a subset $X$ such that  there exists an infinite indiscernible sequence $x_1,x_2,\dots \in X$ such that

(i)  (*) holds for any $1 \leq i_0< i_1< \cdots < i_l < \infty$.

(ii)  For any {\em definable} subgroup $S$ of $G$ with
 $S \subseteq (X \inv X)^2$,  $X$ is not contained in finitely many right translates of $S$. 
 
 Let $\tG$ be the subgroup of $G$ generated by $X$.
By \thmref{normal} there exists $\inft$-definable normal subgroup $S$ of bounded index in $\tG$, with $S \subseteq (X \inv X)^2$.  Since the sequence $x_1,x_2,\ldots $ is indiscernible and $G/S$ is bounded, all $x_i$ lie in the same
coset of $S$.  So the elements $y_i = x_1 \inv x_i$ all lie in $S$.  
Now $d(X,y_1,\ldots,y_l) \leq n$; so  $X \subseteq W_n(y_1,\ldots,y_l; b_1,\ldots,b_l)  $ for some $b_1,\ldots,b_l \in X$.  Let $N$ be the normal subgroup  of $\tG$ generated by the $y_i$, and let $\bX$ be the image of $X$ modulo $N$.  Then
$\bX \subseteq W_n(1,\ldots,1;\bar{b}_1,\ldots,\bar{b_l})$, where $\bar{b}_i=b_iN$.  
Hence $\bX$ is finite.   As $S \subseteq (X \inv X)^2$,
it follows that the image of $S$ modulo $N$ is finite, i.e.
 $[S:N] < \infty$.   
   Since $N$ is $\bigvee$-definable, so is $S$.  But $S$ is $\inft$-definable; so it is a definable group.
 Now $N \subseteq S$, so  $X$ is contained in finitely many translates of $S$, in contradiction to (ii).   
\eprf

By  Ruzsa's argument,     the condition $|a_1^X \cdots a_l^X | \geq |X| /m$ implies that $X \union X^2$ is contained in the union of $\leq m$ translates
$c_j a_1^X \cdots a_l^X (a_l \inv)^X \cdots (a_1 \inv )^X$, with $c_j \in X$; so that 
$d(a_1,\ldots,a_l;G) \leq \max (m,2l)$.   We can now deduce \corref{anti-nil} from \propref{anti-nil2} using Ramsey's
theorem, but will give a direct argument.  We denote the $l$'th Cartesian power of $X$ by $X^{(l)}$.

\prf[Proof of \corref{anti-nil}]    Fix $k,l,m$ and suppose for contradiction that the conclusion fails.   Then for arbitrarily large $K$, letting $p = 1-1/K$,
there exists a group $G_K$ and a finite subset $(X_0)_K, X=(X_0)_K \inv(X_0)_K$, such that the situation of \corref{anti-nil} holds but no subgroup $S$
of $G_K$ with $S \subset X^2$ is contained in $\leq K$ cosets of $S$.   Let $\mu_K$ be the counting measure on $G$, normalized so that $\mu_K(X)=1$,
and let $\mu_K^l $ be the counting measure on $X^l$, similarly normalized.  Let $Q_K$ be the set of $l$-tuples $(a_1,\ldots,a_l)$  such that  
$\mu_K(a_1^X\cdots a_l^X|) \geq 1/m$; then $\mu_K^l(Q_K) \geq p=1-1/K$.  

By compactness (as in the proof of \propref{anti-nil2}, and in \secref{measurequantifiers}) there exists a 
structure including a 
group $G$,  a definable subset $X=X_0 \inv X_0 \subset G$,   a definable measure $\mu$ on definable sets,
as well as a definable measure $\mu^l$ on $l$-tuples, such that  Fubini holds between $\mu$ and $\mu^l$, and  $1=\mu(X) \leq k \mu(X_0) < \infty$.    Further, there exists a definable set $Q \subseteq  X^{(l)}$ such that if $(a_1,\ldots,a_l) \in Q$ then 
$\mu(a_1^X \cdots a_l^X)  \geq 1/m$, and $\mu^l(Q) \geq 1-1/K$ for any $K=1,2,\ldots$, i.e. $\mu^l(Q)=1$.  
Finally, for no definable group $S \subseteq XX$ is $X$ contained
in finitely many cosets of $S$.   

In fact  only one   instance of Fubini will be required:  
$\mu^l(Y^l) = \mu(Y)^l$.

We take $G$ to be countably saturated.  
Recall that countable saturation means that any countable family of definable
sets with the finite intersection property has nonempty intersection; we will actually need it for the family $R_{j}$ below.    

By \thmref{normal} there exists $\inft$-definable normal subgroup $S$ of bounded index in $\tG$.    Find a countable set of definable (with parameters) equivalence
relations $E^j$ on $X_0$, such that each $E^j$ has finitely many classes, $E^{j+1}$ refines $E^j$,
 and if $(a,b) \in E^j$ for each $j$ then $a \inv b \in S$.  (For instance, say $S = \meet S_j$, and let $C^j$ be a maximal
subset such that $S_jx \meet S_jy = \emptyset $ for $x \neq y \in C^j$;   define $E^j$ so
that $(x,y) \in E^j$ implies 
 $\{c \in C^j: x S_j \meet cS_j = \emptyset \} = \{c \in C^j: yS_j \meet cS_j = \emptyset \}$.   Alternatively note that   if $a,b$ have the same type over some countable
model then $a \inv b \in S$.)  

Some class $F_j$ of $E^j$ has measure $\e_j>0$; so $\mu(F_j \inv F_j ) \geq \e_j >0$; thus
 $(F_j \inv F_j )^{(l)} \geq \e_j^l $; and hence (as $\mu(Q)=\mu(X^l)=1$) we have
 $\mu(Q \meet (F_j \inv F_j )^{(l)})  \geq \e_j^l > 0$.  
Hence for each $j$  there exist $(a_1,\ldots,a_l) \in Q$ such that for each $i \leq l$, 
we have $a_i = b_i \inv c_i$ for some $(b_i,c_i) \in E^j$.  As we took the $E_j$ to refine each other,
this holds for any finite set of indices $j$ at once.  In other words, the   family of sets $\{R_j\}$ has
the finite intersection property:
$$R_j= \{(a_1,b_1,c_1,\ldots,a_l,b_l,c_l):  (a_1,\ldots,a_l) \in Q,   \bigwedge_{i \leq l}  (a_i,b_i) \in E^j, \hbox{ and }  a_i = b_i \inv c_i  \}  $$
By countable saturation, $\meet_{j} R_{j} \neq \emptyset$, i.e.  
 there exist $(a_1,\ldots,a_l) \in Q$ and $b_1,c_1,\ldots,b_l,c_l$  
 such that for each $i \leq l$  
we have $a_i = b_i \inv c_i$ and  $(b_i,c_i) \in \meet_j E^j$.  By the choice of $E^j$,
this implies $a_i \in S$.   

Now $S$ is normal in $\tG$, so $a_1^X \cdots a_l^X   \subseteq S$.  Since
$\mu(a_1^X \cdots a_l^X) > 1/(m+1)$, it follows that 
$S$ cannot have $\mu(X_0X)(m+1)$ disjoint cosets $x_iS$.   So $X_0 / S$ is finite; it follows
that $XX/S$ is finite, so $XX = S \union \union_{\nu=1}^k (XX \meet c_iS)$ for some $c_1,\ldots,c_\nu$.  
Since $S$ is $\inft$-definable, so is each $c_iS$, and we see that the 
 the complement of $S$ in $XX$  is also $\inft$-definable.  When a subset of a definable set and its complement are both $\inft$-definable, they are both definable.   Hence $S$ is a definable group.  But finitely many cosets of $S$ cover $X$.  This contradiction proves the corollary.  
\eprf

\medskip

Though we stated \propref{anti-nil2} for finite $X$, it holds with the same proof
if the hypothesis  $|X X \inv X| \leq k |X|$ is replaced by $\mu((X \union X\inv)^3) \leq k \mu(X)$,
with $\mu$ an arbitrary right-invariant finitely additive  measure on $\tG$, or
even in the above sense on $(X \union X\inv)^3$.
 
\bigskip

%
 
 \>{section}

\<{section}{Near subgroups and Lie groups}   \lbl{lie}

Let $X \subseteq \tG$ be a near-subgroup with respect to an $M$-invariant, right-invariant ideal $\mu$, as in the previous section.

 Any compact neighborhood $X$ in a Lie group $L$ is (obviously) an approximate subgroup, and a near-subgroup with respect to Haar measure.   We will show
 that all  near-subgroups are related to these classical ones.   
We will   use logical compactness to connect to the locally compact world, and then   the Gleason- Yamabe structure theory for locally compact groups in order to find  Lie groups.  

\ssec{Some preliminaries}  
We will require the following statement:  
every locally compact group $G$ has an open subgroup $G_{1}$ which is isomorphic to a  projective limit of Lie groups.     (Gleason defines a  topological group $G$ to be 
 a {\em generalized Lie group} if for   every neighborhood $U$ of the 
identity there is an open subgroup $H$ of $G$ and a compact normal subgroup $C$ 
of $H$ such that $C \subseteq U$ and $H/C$ is a Lie group.   (\cite{gleason}, Definition 4.1).
According to  \cite{yamabe}, Theorem 5', every locally compact group is a generalized Lie group.
By \cite{gleason}, Lemma 4.5, if $G$ is a generalized Lie group with connected component $G^{0}$ of the identity,
and $G/G^{0}$ is compact, then $G$ is a projective limit of Lie groups.  Now $G/G^{0}$ is totally disconnected.
So there exists an open subgroup $G_{1}$ of $G$ containing $G^{0}$, such that $G_{1}/G^{0}$ is compact.
Hence $G_{1}$ is an open subgroup of $G$, and a projective limit of Lie groups.)
 
We will also use the fact that in a connected Lie group $G$,  for any chain $C_1 \subset C_2 \subset \cdots $ of compact  normal   subgroups, $cl( \union_n C_n)$ is
also a compact normal subgroup.     Indeed the dimension of the Lie algebras of the $C_n$  must stabilze, so they are locally equal, and hence the connected components $C_n^0$ stabilize.   Factoring out the compact normal subgroup $\union_n C_n^0$, we may assume the $C_n$  are discrete, i.e. finite.   Since $G$ is connected,
the $C_n$ are contained in the center $Z$.   The connected component $Z^0$ of $Z$ has universal covering group $\Rr^n$, so $Z^0 \cong \Rr^k \oplus (\Rr/\Zz)^l$.
The  discrete group $Z/Z^0$ is a homomorphic  image of the fundamental group of $G/Z$,  hence is finitely generated; it has a finite torsion part $A/Z^0$; since $Z^0$ is divisible,
$A$ can be written as a direct sum $A_0 \oplus Z^0$.    It is clear that the torsion points 
of $Z$, and hence all the $C_n$, are contained in the compact central subgroup $A_0 \oplus  (\Rr/\Zz)^l$.  
 
   In particular any
closed subgroup contains a unique maximal compact  normal subgroup of $G$. 

 Further down (\lemref{sop1}), we will also need to know that a  compact Lie group has no infinite descending
sequences of closed subgroups; this follows easily along the same lines.

 The  results of this section will also be valid for local groups, using the following local version of Gleason-Yamabe 
  due to Goldbring:
    for a compact local group $G$  there exist a continuous map $h:D \to L$ into a Lie group $L$, whose domain $D=D ^{-1}$ is a smaller compact
neighborhood of $1$ in $G$,  and whose image $hD$ is a compact neighborhood of $1$ in $L$, such that $xy$ is defined for any $x,y \in D$,
and we have:
$xy \in D$ iff $h(x)h(y) \in hD$, in which case $h(xy)=h(x)h(y)$.   
\cite{goldbring} has generalized the ``no-small-subgroups" theory to the
local group setting; to apply it 
one needs to know that some neighborhood of $1 \in G$
contains a compact normal subgroup, such that the quotient has no small subgroups; this is  Lemma 9.3 of \cite{goldbring}.   

Recall that  we call two subsets  $X,X'$ of a group {\em commensurable} if each one is contained in finitely many right translates of the other.  
If $H,H'$ are subgroups, and $H$ is contained in finitely many cosets of $H'$, then it is contained in the same number of cosets of $H \meet H'$,
so $[H: H\meet H'] < \infty$; thus for groups this coincides with the usual notion.

\<{thm}\lbl{ap1}  Let $X$ be a near-subgroup of $G$, generating a group $\tG$.   Then there exists a $\bigvee$-definable subgroup
 $\breve{G}$    contained in $\tG$,   a $\inft$-definable subgroup $K \subseteq \breve{G}$, 
    a connected, finite-dimensional Lie group $L$,  with no nontrivial normal compact subgroups, 
 and a homomorphism $h: \breve{G} \to L$ with kernel $K$ and dense image, with the following property:
 
    If  $F \subseteq F' \subseteq L$ with $F$ compact and $F'$ open,   then there exists a definable $D$
with $h \inv (F) \subset D \subset h \inv(F')$.   Any such $D$ is  commensurable to $X \inv X$.
 
$\breve{G}$ and $K$ are defined without parameters.  The Lie group $L$ is uniquely determined. 
 \>{thm}  

Let us bring out some facts implicit in the statement of the theorem (and also visible directly in the proof.)  
\<{remark} \lbl{ap1-rems} 
 \<{itemize}
\item  If $(G',X')$ is a countably saturated elementary extension of $(G,X)$, then $h$ extends to
$h': \breve{G}' \to L$, and $h'$ is surjective.
\item  The Lie group $L$ is determined up to isomorphism by $(\tG,\cdot,X)$,
where $\tG$ is the subgroup of $G$ generated by $X$; 
 in fact by the theory of $(\tG,\cdot,X)$, with $\tG$ viewed as many-sorted.  We call it the associated Lie group.
\item Since any compact subset of a Lie group is a countable intersection of open sets, it follows
that if $W \subseteq L$ is compact, then $h \inv(W)$ is $\inft$-definable.
\item Similarly,  if $W \subseteq L$ is open, then $h \inv(W)$ is $\bigvee$-definable.
\item  If $W \subseteq L$ is a neighborhood of $1$, then $h \inv (W)$ contains a definable set 
 of the form $U \inv U$, with $U$ a definable subset of  $(G,X)$ contained in $\breve{G}$ and commensurable to $X \inv X$.  
 \item  Any definable set containing $K$ contains some $h \inv (W)$, with $W$ a neighborhood of
 $1$ in $L$.
 \item If $L$ is trivial, taking $F=F'=L$ in the statement of the theorem we see that $\breve{G}$ is a definable group, commensurable to $X \inv X$.
 \item We have $K \subseteq (X X\inv)^m$ for some $m$.    \thmref{stabilizer} provides
 an $\inft$-definable stabilizer contained in $(X \inv X)^2$, but converting it to a 0-definable one involves some
 (finite) enlargement.  
 
\>{itemize}
\>{remark}

We first show the main statement of \thmref{ap1} holds after saturation and base change

\<{lem} \lbl{ap1.01}    Let $X$ be a near-subgroup of $G$, generating a group $\tG$.   Assume the structure $(G,X, \ldots)$ is countably saturated.
Then over parameters there exists a $\bigvee$-definable subgroup
 $\breve{G}$    contained in $\tG$,   a $\inft$-definable subgroup $K \subseteq \breve{G}$, 
    a connected, finite-dimensional Lie group $L$  and a homomorphism $h: \breve{G} \to L$ with kernel $K$ and dense image, with the following property:
 
    If  $F \subseteq F' \subseteq L$ with $F$ compact and $F'$ open,   then there exists a definable $D$
with $h \inv (F) \subset D \subset h \inv(F')$.   Any such $D$ is  commensurable to $X \inv X$.
   
\>{lem}

\prf     Let $\tG$ be the subgroup of $G$ generated by $X$; let $S_0=XX \inv$.  
      \thmref{normal} (via  \corref{stab1})  provides definable subsets $S_n \subseteq (X X \inv)^2$  of $\tG$ such that $S=\meet_{n \in \Nn} S_n$ is  normal subgroup of $\tG$ bounded index; we may take $S_{n+1}=S_{n+1} \inv$ and $S_{n+1}S_{n+1} \subseteq S_n$.  Define $S_n$ for negative $n$ too by  $S_n=S_{n+1}S_{n+1}$.  So   $S_0$ is $0$-definable, and $\union S_n = \tG$.  

We define a topology on $\tG/S$ using the quotient map $h: \tG \to \tG/S$  by:

(*)  $W \subset \tG/S$ is closed iff $h \inv(W) \meet S_n$ is $\inft$-definable for each $n$.

See  \cite{hpp}, Section 7 for a more detailed description. Let   $L_{0}=\tG/S$.     This is easily seen to be a locally compact toplogical group.  Compactness is an immediate consequence
of   saturation and logical compactness:  an intersection of a small number  of $\inft$-definable subsets can never be empty, unless a finite sub-intersection is empty.  Continuity of the group operations follows from the definability
of the group structure on $G$.  The images of the sets $S_{n}$ form a neighborhood basis for the identity of $L_{0}$
as noted below, so that $L_{0}$ is Hausdorff.

Let $\pi_S: \tG \to \tG/S$ be the projection.  
Note
that   $\pi_S\inv \pi_S(S_n) = \meet_m S_nS_m \subseteq S_nS_n$.  In particular $\pi_S \inv \pi_S(S_n)$ is contained
in a definable subset of $\tG$.  In fact for any definable set $D \subseteq S_n$, 
$\pi_S \inv \pi_S(D) = SD$ is an $\inft$-definable subset of $S_{n-1}$.   More
generally for any locally definable\footnote{see definition in the first lines of \S 3} subset $D$ of $\tG$ , $\pi_S(D)$ is closed.
Indeed $\pi_S \inv \pi_S(D) \meet S_n = \pi_S \inv (\pi_S( D \meet S_nS_n))$.

In particular,  the image of $\tG \m  S_nS_n$ in $\tG / S$ is closed, and disjoint
from  $\pi_S(S_n)$; since $\pi_S(S_nS_n) \union \pi_S(\tG \m S_nS_n) = \tG/S$,
  $\pi_S(S_n)$ lies in the interior of $\pi_S(S_nS_n)$.  In particular, each $\pi_S(S_nS_n)$ is a neighborhood of $1$, as is theorefore $\pi_S(S_{n+1})$. 

By  Yamabe,  $L_{0}$ has an open subgroup $\breve{G}/S$, isomorphic to a projective limit of Lie groups.  
 $\breve{G}/S$ is also closed, so both $\breve{G} \meet D$ and
$D \m \breve{G}$ are $\inft$-definable, for any definable $D$ contained in $\tG$.  Thus $\breve{G}$ is locally definable in $\tG$, i.e. it has a definable intersection with any definable subset of $\tG$.

The topology  of a projective limit $\liminv L_i$  is generated by pullbacks of open subsets of individual factors $L_i$.  So there
exist   a Lie group $L$,   a neighborhood $U_1$ of the identity
in $L$, and a homomorphism  $h: \breve{G}/S \to L$, such that $h \inv(U_1) \subseteq \pi_S(S_1)$.  
  By shrinking $\breve{G}$ down further
to the pullback of the (open) connected component of $1$ in $L$, we can take $L$ to be connected.
Let $\pi: \breve{G} \to \breve{G}/S \to L$ be the composition.  

Now (*) holds for $L$ :  the morphism from a projective limit to one of the factors is closed; so  $Y \subseteq L$ is closed iff $h \inv(Y)= \pi \inv (Y) / S$ is closed
iff $\pi \inv(Y)$ meets every definable set in an $\inft$-definable set.

We also have:  (**) For any compact neighborhood $U$ of $1$ in $L$, $\pi \inv(U)$ is commensurable to 
$X \inv X$.    For any two compact  neighborhoods of $1$ in $L$ are commensurable, each one being 
contained in a union of translates of the other, which can be reduced by compactness to a finite 
union.  This comparability is  preserved
 by $\pi \inv$.   So it suffices to show that $\pi \inv(U)$ contains $X \inv X$ for some $U$, and that
 $\pi \inv(U')$ is contained in finitely many translates of $X \inv X$ for some $U'$.   On the other hand
 by the Ruzsa argument (above \lemref{gen-equiv}), any $S_n$ is commensurable to $X \inv X$.
 We saw that $\pi_S(X \inv X)$ is compact; hence $\pi(X \inv X)$ is compact, so it is contained in 
 some compact open neighborhood $U$, and thus $X \inv X \subseteq \pi \inv (U)$.  And
 by construction, $\pi \inv(U_1) = \pi_S \inv h \inv(U_1) \subseteq S_1$, giving the second direction.

 If  $F$ is a compact subset of $L$ and $F'$ an open subset, with $F \subset F'$, then there exists a definable $D$
with $h \inv (F) \subset D \subset h \inv(F')$.   
 Indeed $\pi \inv(F)$ is an $\inft$-definable set contained in the $\bigvee$-definable set $\pi \inv(F')$, so there exists a definable $D$ with $\pi \inv(F) \subseteq D \subseteq \pi \inv(F')$.  
 
\eprf

We now begin to address the issue of parameters.  

\<{lem}    \lbl{ap1.02}  With the assumptions and notation of \lemref{ap1.01}, there exists an $\inft$-definable
subgroup $S$ of $\tG$ without parameters, with $\tG/S$ bounded.  \>{lem}

\prf 
 We may work in a homogeneous elementary extension $\Uu$ of $(G,X,\cdot)$, so that  $\inft$-definable sets are 
$\inft$-definable without parameters   as soon as they are $Aut(\Uu)$-invariant.  

Let  $\a$ be the set of pairs $(H,\G)$ such that $\G \leq H \leq \tG$,  
and for some small base $A$, $\G$ is a normal subgroup of $H$, $\G$ is $A$-$\inft$-definable,  $H$ is an locally definable subgroup of $\tG$ over $A$, and $\tG/ \G$ is bounded. 
Let $\b$ be the set of pairs $(H, \G) \in \a$ such that if $(H',\G') \in \a$ and $\G \leq \G' \leq H' \leq H$ then 
$H=H'$ and $\G=\G'$.  Equivalently, the locally compact group $H/\G$ is connected, with no nontrivial compact normal subgroups.   (Hence by Yamabe, is a Lie group.) 

For$ (H,\G) \in \b$ it is clear that $H$ determines $\G$, since if $(H,\G') \in \b$ then $\G= \G \G' = \G'$.

\claim{1}  $\b$ is nonempty.  

\prf We saw above that there exists $(H,\G) \in \alpha$ with $H/\G$ a connected Lie group.  In the preliminaries
to this section we saw that $H/\G$ has a maximal compact normal subgroup; it has the form $H/\G'$ with
$\G'$ $\inft$-definable.  Then $(H,\G')$ is in $\b$.  \eprf

\claim{2}  
Let $(H,\G),(H', \G') \in \b$.  Then 
$(H \meet H', \G \meet H') \in \b$. 

\prf  Since $H'$ is locally definable, while $\G$ is contained in a definable set, it is clear that 
$H' \meet \G$ is $\inft$-definable.  Since $\tG/\G$ and $\tG/H'$ are bounded, so is $\tG/(\G \meet H')$.
Also $H \meet H'$ is locally definable.  
Thus $(H \meet H', \G \meet H') \in \alpha$.  

 Now $\G'/(\G' \meet H)$ is bounded (as it embeds into $\tG/H$).
  By \lemref{finiteindex}, $\G' \meet H$ has finite index in $\G'$.  
  
Similarly, $\G$ is contained in finitely many costs of $H'$, hence of $H' \meet H$.  So $\G(H' \meet H)$ is a finite union 
of cosets of $H' \meet H$, and hence is a locally definable subgroup of $H$.  We saw that
for any definable set $D$ containing $\G$, the image of $D \inv D$ contains an open neighborhood of the identity.
Hence the image of $\G(H' \meet H)$ in $H/\G$ is open.  

Now  the natural map
$(H \meet H') / (\G \meet H') \to H / \G$ is   injective.  
But it has open image and the group $H / \G$ is connected, so the map is surjective.  Thus  
$(H \meet H') / (\G \meet H') \cong H / \G$ and hence has no nontrivial compact normal subgroups.  \eprf

 Similarly $(H \meet H', \G' \meet H) \in \beta$.  So $\G' \meet H = \G \meet H' $ and thus 
 $\G' \meet H = \G \meet \G'$. 
  
We noted that $\G' \meet H$ has finite index in $\G'$; moreover since this holds for any pair from $\a$,
in particular it holds for $\G'$ and any $Aut(\Uu)$-conjugate $\si(H)$ of $H$, so $\G' \meet \si(H)$ has index
bounded independently of  $\si$. 
  
 So $\G \meet \si(\G')$ has finite index in $\G'$, bounded independently of $\si \in Aut(\Uu)$.  By
 symmetry, $\G,\G'$ are commensurable, and all conjugates of $\G$ are uniformly commensurable.

 Pick $\G_1 \in \b$.      By  
\cite{bergman-lenstra}, there exists an $Aut(\Uu)$-invariant group $S_1$ commensurable to each conjugate
of $\G_1$.  
The proof of \cite{bergman-lenstra} shows that $S_1$ contains a finite intersection of conjugates of $\G_1$ as a subgroup of finite index; so $S_1$ is $\inft$-definable, and of bounded index in $\tG$; being $Aut(\Uu)$-invariant, it is $\inft$-definable over $\emptyset$.  Let $S$ be the intersection of all $\tG$- conjugates of $S_1$; then $S$
is normal in $\tG$, $\inft$-definable over $\emptyset$, and of bounded index.  
\eprf

 \prf[Proof of \thmref{ap1}]  We may assume $(G,X)$ is countably saturated, since the statements descend from a saturated extension of $(G,X)$ to $(G,X)$ by restriction, using the same ($0$-definable) $\breve{G},h,L$.    Let $S$ be the $0$-$\inft$-definable group given by \lemref{ap1.01}.   Since $\b \neq \emptyset$
 in \lemref{ap1.01}, we know that $\breve{G},h$ exist over parameters, and it remains only to show that 
 $\breve{G}$ and $\ker(h)$ can be chosen to be $\bigvee$- definable and $\inft$-definable (respectively) without parameters.  

  We begin with $\breve{G}$.  We may replace $\breve{G}$ by the pre-image of any open subgroup of $\breve{G}/S$ 
(the``connected-by-compact" condition  will remain valid.)
Let $G_c$ be the group generated by $\breve{G} \meet S_1$; note that
$G_c$ is locally definable, and is generated by $Z_c = G_c \meet S_1$; so each of $G_c$, $Z_c$ can be used to define the other.    
Now $Z_c$ is a definable set, with parameter $c$ say.  Let $Q$ be the set of realizations of $tp(c)$.  If $c' \in Q$,  then $Z(c')$
generates a group $G_{c'}$, and $G_{c'} \meet S_{1} = Z(c')$.  Thus for $c'',c' \in Q$, $G_{c''}=G_{c'}$ iff
$Z_{c''}=Z_{c'}$; this is a definable equivalence relation.  
%


We   have $h: \breve{G} \to L$.   Note that if $C$ is a   compact normal subgroup
of $L$, the composition of $h$ with the quotient map $L \to L/C$ has the same properties
(1,2) as $h:\breve{G} \to L$.  
Replacing $L$ by $L/C$ for a maximal compact normal subgroup $C$ of $L$, we may assume $L$ has no compact normal subgroups. 

 Let  $K$ be the kernel of $h$.  Then $K$ is $Aut(\Uu)$-invariant.  
For if $K'$ is an  $Aut(\Uu)$-conjugate of $K$, then $K,K'$ are 
$\inft$-definable normal subgroups of $\breve{G}$;
$K'K /K$ is a compact normal subgroup of $\breve{G}/K$, hence it is trivial;
and similarly $K'K/K'$ is trivial; so $K=K'$.   Thus $K$ is    $\inft$-definable without parameters.

It remains to prove the uniqueness of $L$.  Let us compare $L$ to the locally compact group
 $\tH:=\tG/S$, where $S=\tG^{00}_{\emptyset}$ is the smallest $0$-$\inft$-definable subgroup of $\tG$ of bounded index.
   Let $H$ be the image of $\breve{G}$ in $\tH$.  Then $H$ is an open subgroup of $\tH$,
so the connected component of the identity $\tH^0$ is contained in $H$, and equals $H^0$.  Let $C$ be the image of $K$ in
$\tH$.  So $C$ is a normal subgroup of $H$.  Since $H/C$ is connected, we have $H/(C \tH^0)$ both connected and totally disconnected.  (Unlike the situation in the category of topological spaces, in the category of topological
groups the image of a totally disconnected group is still totally disconnected. Indeed it has a pro-finite open subgroup, and this remains the case for a quotient group.)  So $CH^0=H$.  Both $C$ and $H^0$ are normal in $H$, 
so letting ${C_0}=C \meet H^0$ we have $H/{C_0} \cong C/{C_0} \times H^0/{C_0}$.  Thus the action of $C$ by conjugation on $H^0$ is trivial modulo ${C_0}$.    Now $C$ is a maximal normal compact subgroup of $H$; ${C_0}$ is a   compact normal subgroup of $H_0$, maximal with respect to being normalized by $C$ too, but we have just shown that this last condition is trivial, so ${C_0}$ is a maximal compact normal subgroup of $H_0$.  We have $L = H/C \cong H^0/{C_0}$ canonically.   Now it is clear that $C_0$ is the unique maximal normal compact subgroup
of $H^0$.  (If $C_1$ where another, $C_0C_1$ would be still bigger.)  
This proves the uniqueness of $L$.   \eprf

If we expand  $(\tG,X,\cdot,\ldots)$ to a structure $\underline{\tG} = (\tG,X,\cdot,\ldots,R_{new},\ldots)$ with 
 more definable sets, the smallest  $0$-$\inft$-definable subgroup  of bounded index may become smaller: 
   $\underline{\tG}^{00}_0 \subset  {\tG}^{00}_0$.   Thus $H=\tG/S,H^0,C_0$ will change with the added structure.  Nevertheless
the isomorphism proved in the last paragraph of the proof remains valid; hence   {\em the associated Lie group $L =   H^0/{C_0}$ does not change if the structure is enriched.}

\<{defn}   \lbl{associatedlie}
Let $\tG$ be a $\bigvee$-definable group, $X$ a definable near-subgroup of $\tG$, generating $\tG$.  Let $M=(\tG,X,\cdot,\ldots)$,
where $\ldots$ indicates possible additional structure.    \<{itemize}

\item   $LC(M) = \tG/S$, where $S$ is the smallest $\inft$-definable subgroup of $\tG$, without parameters,
of bounded index.

\item  $L(X)$ is the Lie group associated to $X$; so $L(X)= LC(M)^0 / C_0$, with $C_0$ a maximal normal compact
subgroup of $LC(M)^0$.    Let $\widehat{L}(X) = LC(M)/C_0$; then $L(X) = \widehat{L}(X)^0$, the connected component.

\item  $l(X) = \dim L(X)$.      \>{itemize}  \>{defn}

We refer to $l(X)$ as the Lie rank of $X$, or of $\tG$.

 \<{example}\lbl{lie0}  If a near-subgroup $X$ has $l(X)=0$, then there exists a definable group $S$ with $X,S$ commensurable.  \rm  Indeed in this case kernel $S$ of the homomorphism $\breve{G} \to L$ is equal to $\breve{G}$; but $S$ is $\inft$-definable and $\breve{G}$ is $\bigvee$-definable,
so they are definable.  
  \>{example}

\<{lem} \lbl{lie-doubling} 
   In the situation of \thmref{ap1}, assume the S1-ideal arises from an invariant, translation invariant measure $\mu$.  Let $k_5 = \mu(X X \inv X X \inv X) / \mu(X)$.
  Extend $\mu$ to the $\si$-algebra generated by the $\infty$-definable subsets of $\tG$,
and let $\lambda$ be the pushforward of $\mu$ to $\widehat{L} = \widehat{L}(X)$, i.e. $\lambda(U) = \mu(\pi \inv(U)) \in \Rr_\infty$,
where $\pi$ is the quotient map.  Then $\lambda$ is  a Haar measure on $\widehat{L}$.
We have $\lambda( (\pi X)(\pi X) \inv(\pi X)) \leq k_5 \lambda( \pi X)$.  Moreover, there exists a compact subset $W$ of  $L=L(X)$
with $\lambda(W) >0$ and $\lambda(W W \inv W ) \leq k_5 \lambda(W)$.  We can take $1 \in W$.  
 \>{lem}  
 
 \prf    It is clear that $\lambda$ is a nonzero,  translation invariant measure, hence   a Haar measure. 
We have  $X X \inv X  \subseteq  \pi \inv( (\pi X) (\pi X) \inv (\pi X)) \subseteq  (X X \inv X X \inv X)$  , since $\pi \inv  (1)   \subseteq X \inv X$.   
This implies the first inequality, by definition of the pushforward measure.  
 Moving to $L$, recall that we have $h: \breve{G} \to L$ with kernel $K$ (\thmref{ap1}), with $\breve{G}$ a $\bigvee$-definable subgroup of $\tG$,
 and $\tG/\breve{G}$ bounded.  In particular $X / \breve{G}$ is bounded, so $X$ intersects finitely many cosets of $\breve{G}$;
 say $X= \union_{i=1}^r  X_i$, with $X_i \subseteq c_i \breve{G}$, and $c_i$ lying in distinct cosets of $\breve{G}$.  
 Let $k_3 = \mu(X  X \inv X ) / \mu(X)$.  Then, noting that $ X_i X_i  \inv X_i \subseteq c_i \breve{G}$, and the $c_i \breve{G}$ are disjoint, we have: 
$$\sum_i  \mu(X_i X_i \inv X_i)    \leq \mu(X  X \inv X ) \leq k_3 \mu(X) = \sum_i k_3 \mu(X_i)$$
The sum being extended over all $i \leq r$ such that $\mu(X_i)>0$.      It follows that for at least one $i$ with $\mu(X_i)>0$, we have
$\mu(X_iX_i \inv X_i) \leq k_3 \mu(X_i)$.   Similarly, for at least one $i$ with $\mu(X_i)>0$ we have
$\mu(X_i X_i \inv X_i X_i \inv X_i) \leq k_5 \mu(X_i)$.    Let $Y= c_i \inv X_i$.  Then $h(Y)$ is a compact subset of $L$;
$\lambda(h(Y)) = \mu(h\inv h(Y)) \geq \mu(Y)>0$; and $\lambda(Y Y \inv Y) \leq k_5 \lambda(Y)$ by the same argument as for $\widehat{L}$ above. 
By translating $W$, we can arrange $1 \in W$.
  \eprf

Can \lemref{lie-doubling} be used to bound $l(X) = \dim(L)$ in terms of doubling constants of $X$?  
When $G$ is nilpotent, we have:    $l(X) \leq \log_2(k_5)$.     
This follows  from  \lemref{lie-doubling}  and  \lemref{nil-doubling},
due   (with a different proof) to Tsachik Gelander; thanks for allowing me to include it here.
   Use $1 \in W$ to obtain $WW   \subseteq W W \inv W$ in order to apply the lemma.

  \<{lem}[Gelander]  \lbl{nil-doubling} Let $X$ be a compact subset of $\Rr^d$, or more generally of a connected, simply connected Lie group, and let  $\lambda$ be Haar measure. Then $\lambda(XX)\geq 2^d \lambda(X)$.  \>{lem}

\prf  In fact we have $\lambda(s(X)) \geq 2^d \lambda(X)$, where $s(x)=x^2$.  The ambient group $H$ is isomorphic to a subgroup of the strict upper
triangular matrices over $\Rr$, of some dimension; the map $s$ is hence injective.   Moreover $H$ is diffeomorphic to $\Rr^d$, and 
the differential $ds$ of $s$  at any point is a linear transformation of the form $2+M$, with $M$ nilpotent.  It follows that the Jacobian determinant 
has value $2^d$, so by the change of variable formula for integration, the diffeomorphism $s$ expands volume by exactly $2^d$.
\eprf

  \<{rem}  \lbl{towers} \rm  (Compare   \cite{tao-nc}, Lemma 7.7 and Theorem 7.12.)
  
Let     $\G$ be a $\bigwedge$- definable subgroup
of bounded index in the $\bigvee$-definable group $\tG$.  Let $\tN$ be a  locally definable  normal subgroup of $\tG$,
and let $\pi: \tG \to \tG/\tN$ be the quotient map.  The main case is that $\tN$ is the intersection with $\tG$ of a definable normal subgroup $N$ of $G$.

(0)    The image $\mbG$ of $\G$ has bounded index in the image $\bG$ of $\tG$ modulo $\tN$, and also $\G \meet \tN$ has bounded index in $\tN$.  Conversely in this situation the boundedness of $\tG/\G$ follows from that of $\bG/ \mbG$ and of 
$\G \meet \tN$ in $\tN$.  
 
(1)   View $\tG/\G$, $\bG / \mbG$ and $\tN/ (\G \meet \tN)$  as locally compact groups.  
Then $\tN/ (\G \meet \tN)$ with the logic topology is homeomorphic to the image
of $\tN$ in $\tG/\G$, with the subspace topology.     
 Indeed the natural map  $\tN/ (\G \meet \tN) \to \tG / \G$ is  a continuous injective homomorphism.    To see that it is also a closed map, since $ \tG / \G$ is covered by the   interiors of sets of the form $\pi(D)$, with $D$ definable, we may restrict
 attention to the inverse image of such a set.  But then we are looking at an injective
 continuous map between compact Hausdorff spaces, hence an isomorphism.  

Similarly,  $\bG / \mbG \cong (\tG/\G) / (\tN/ (\G \meet \tN))$ as topological groups. 

(2) 
 $\G$ is definable iff the topology on $\tG/\G$ is discrete.   This makes it plain that $\G$ is definable iff $\pi(\G)$ and $\G \meet \tN$ are.

(3)   If $\tG/\G$ is a Lie group then so are $\tN/(\G \meet \tN)$ and $\bG/{{\mathbf \G}}$, and we have
an exact sequence
$$1 \to \tN/(\G \meet \tN) \to \tG/\G \to \bG/{{\mathbf \G}} \to 1$$
This in turn induces an exact sequence of homomorphisms among the Lie algebras.   
It follows  that 
$\dim(G/\G)=\dim(\bG/{{\mathbf \G}}) + \dim(\tN / (\G \meet \tN))$.  

(4)    From (3) it follows that
  $$l(\tG) \geq l(\tG/\tN)+l(\tN)$$ 
  Indeed we may move from $\tG$ to $\bG$, changing none of the three numbers.  Then we may enlarge $\G$
 so that $\bG/\G$ has no nontrivial normal compact subgroups.   By (3) we obtain
 in this situation:  
  $l(\tG) = \dim(\bG/{{\mathbf \G}}) + \dim(\tN / (\G \meet \tN))$.   Now  $\tN/(\G \meet \tN)$ may have nontrivial compact subgroups,
  but we have at all events $l(\tN) \leq \dim(\tN / (\G \meet \tN))$ (the inequality may be strict.)  Similarly $l(\tG/\tN) \leq \dim(\tG/\tN)$, and
  (4) follows.  
  

See \S 7 for a an inductive use of this invariant, similar to Gromov's use of the growth rate in the case of his polynomial growth assumption.

 \>{rem}

\<{rem}  \rm    The canonicity of $L$ in \thmref{ap1} is achieved at a price.  We noted already that it requires moving from $(X \inv X)^2$
to $(X \inv X)^m$ where $m$ is difficult to control.  In addition, factoring out the maximal compact normal subgroup can lead to substantial
loss of information.  

In some cases there will exist  a largest   $\inft$-definable normal
subgroup $\Delta$ of $\tG$ with $\Delta \subseteq X\inv X$.  By Yamabe,   
  $\tG/\Delta$ is a Lie group $\widetilde{L}$.  In this case $\widetilde{L}$ too is an invariant
of $(\tG,X)$, and is superior in both respects.  When it exists, we may call  $\widetilde{L}$ the directly associated Lie group.  More generally we may need to look at a number of $\tG/\Delta$, differing by compact isogenies.  

For example,  let $\a > 10$ be an irrational real number, and let 
 $$X[n] = X[n,\a]= \{ [m \a]:  m \in \Zz, -n \leq m \leq n \}$$
 where $[m \a]$ is the integer part of $m \a$.
 $X[n]$ is symmetric, and satisfies $|X[n]X[n]| / |X[n]| \leq 4$.      Let $(G,X,n^*)$ be a nonprincipal 
 ultraproduct of $(\Zz,X[n],n)$.  Then the directly associated Lie group is the product of the circle 
 $\Rr/\a \Zz$ with $\Rr$.  The map $\tG \to \Rr$ takes $x$ to  the standard part of $x/n^*$.  The map
 $\tG \to \Rr/ \a \Zz$ takes $x$ to the standard part of the image of $x$ in the nonstandard circle $\Rr^* / \a$.  
 The image of $X$ in the cylinder $\Rr \times \Rr/ \a \Zz$ is  the image of the square $[-\a,\a] \times [-1,0]$.    The image of the element $[\a]$ is $(0,m + \a \Rr)$ for
 some nonzero integer $m$; it follows that $(0) \times \Rr/\a \Zz$ is contained
 in the image of $\tG$, so that $\tG \to \Rr \times \Rr/ \a \Zz$ is surjective.  
  The doubling of this square within the cylinder is similar to the doubling of the $X[n]$ within $\Zz$.  By contrast
the  associated Lie group without compact subgroups is $\Rr$, which does  not account  
 for the doubling of $X$ or $XX$ very well, and only begins to work around the $[\a]$'th set power of $X$.   
 
 It is also interesting to note here that if one takes $X[n]' = \{ [m \a_n]:  m \in \Zz, -n \leq m \leq n \}$ where
 $\a_n$   approaches $\infty$, the associated Lie group will be $\Rr^2$; this limit is natural
  for the directly associated Lie group but not for the reduced one.  
  \>{rem}

For the record we state a version of \thmref{ap1} 
waiving canonicity but  
gaining more control of the location of the kernel.

\<{lem} \lbl{ap1-2}  Let $X$ generate a $\bigvee$-definable group $\tG$, and assume an ideal on $\tG$ exists satsifying
the assumption of \lemref{sym}.  
   Then there exists a $\bigvee$-definable subgroup
 $\breve{G}$ contained in the group generated by $X$,   a $\inft$-definable subgroup $K \subseteq \breve{G}$, 
    a connected, finite-dimensional Lie group $L$  and a homomorphism $h: \breve{G} \to L$ with kernel $K \subseteq (X \inv X)^2$ and dense image, such that:

    If  $F \subseteq F' \subseteq L$ with $F$ compact and $F'$ open,   then there exists a definable $D$
with $h \inv (F) \subset D \subset h \inv(F')$.   Any such $D$ is  commensurable to $X \inv X$.
 
$\breve{G}, K$ may be defined with parameters in any given model.  \>{lem}  

\prf  By \lemref{sym} and \thmref{stabilizer}, one obtains an $\inft$-definable stabilizer $S$ defined over a given  model, and with $S \subseteq (X \inv X)^2$.  It follows that the image $U$ of $(X \inv X)^2$ in $\tG/S$ contains the identity in its interior.  
We follow the proof of \thmref{ap1}, taking care to factor out only by a   compact subgroup contained in the given neighborhood $U$.  \eprf

In the local group setting, the  conclusion reads:  there exists a homomorphism $h: W \to L$
of local groups, $W$ a subset of $X$ commensurable to $X \inv X$, such that 
  $\bX=h(X)$ is a compact neighborhood of $1 \in L$;   and if   $F \subseteq F' \subseteq \bX$ with $F$ compact and $F'$ open,   then there exists a definable $D$
with $h \inv (F) \subset D \subset h \inv(F')$.   Any such $D$ is again commensurable to $X \inv X$.
I have not checked the question of parameters for local groups.

\<{cor}  \lbl{ap1.1}  Let $X$ be a near-subgroup of a group $G_0$, generating $\tG$.   Then there exist 0-definable
subsets $X_{1},X_{2},\ldots$ of $\tG$,  commensurable to $X \inv X$,  and $c \in \Nn$, with:  

\<{enumerate} 
\item  $1 \in X_{n}=X_{n} \inv$
\item   $X_{n+1}X_{n+1} \subseteq X_{n}$ 
\item $X_{n}$ is contained in $\leq c$ translates of $X_{n+1}$.
\item   $aX_{n+1}a \inv \subseteq X_{n}$ for $a \in X_{1}$. 
\item  $[X_{n},X_{m}] =\{x  y x \inv  y \inv: x \in X_n, y \in X_m\} \subseteq X_{k}$ whenever  $k < n+m$.
In particular each $X_n$ is closed under the commutator bracket.
\item  $X_{n+1}=\{x \in X_{1}: x^4 \in X_n \}$   
\item   Let $x,y \in X_{m}, m\geq 2$ and suppose $ x^2=y^2$.  Then $x y\inv \in \meet_n X_n$.
\>{enumerate}
\>{cor}

\prf     We may assume $(G_0,X)$ is $\aleph_0$-saturated.  Let $h,L$ be as in \thmref{ap1}.
 We first show that $L$ has a system $U_{n}$ of compact neighborhoods of the identity with
  properties (1-3).
\def\mL{\mathfrak{L}}
Let $\mL$ be the Lie algebra of $L$, $exp: \mL \to L$ the exponential map, and
fix a Euclidean inner product on $\mL$.  Let
 $V$ be a simply connected open neighborhood 
of $0 \in \mL$ such that $exp$ is a diffeomorphism $V \to exp(V)=U$, and such that
the image of $(X \inv X)^2$ in $L$ contains $U$ in its interior.
  Let $V_n$ be the ball of 
radius $r_02^{-n}$ around $0$.  Here $r_0>0$ is chosen small enough so
that  $V_{0}$ is contained in $V$; some further constraints  on $r_0$ will be specified later.   Viewing $\mL$ as the tangent space
at $1$ of $L$, fix on $L$ the unique left-invariant Riemannian metric extending the given
inner product at $1$.    
Let $U_{n}=exp(V_{n})$.  Note that $U_n$ is    the set of points of $U$
 at Riemannian distance $\leq  r_02^{-n}$ from the identity element (cf. e.g. \cite{lee}, Prop. 6.10).  It follows that
(1-2) hold:  $1 \in U_n = U_n \inv$  and $U_{n+1}U_{n+1} \subset U_n$.  

Fix an invariant volume form $\omega$ on $L$. 
We claim that for some constant $c' > 0$, we have  $vol(U_{n+1}) \geq c' \vol(U_n)$ for large enough $n$.    We have $vol(U_n) = \int_{V_n} exp^* \omega$ where $exp^* \omega$ is the pullback.    Now $V_{n}$ has volume proportional to $2^{-nd}$,
with respect to the standard Euclidean volume form $\omega_1$.
We have $exp^* \omega = f \omega_1$ for some non-vanishing smooth function $f$,
that we can take to be positive.  On $V$ we have $(c'') \inv \leq f \leq c''$ for some $c''>0$, so 
$vol(U_{n}) \leq c'' \vol(V_{n}) \leq 2^{d}c'' \vol(V_{n+1}) \leq 2^{d}(c'')^{2} \vol(U_{n+1})$.
 
Now $U_{n-1}$ contains at most $\vol(U_{n-1})/\vol(U_{n+2})$ disjoint $U_n$- translates of $U_{n+2}$; hence $U_{n}$ is contained in that many translates of $U_{n+2} \inv U_{n+2} \subseteq U_{n+1}$.   This gives the analogue of (3).  

To obtain (4), we may begin with $r_1$ small enough so that for $x \in U_{1}$,   $1-ad_{x}$ has operator norm $<1/2$.  Then $ad_{x}(V_{n+1}) \subseteq V_{n}$, so $x \inv U_{n+1}x \subseteq U_{n}$.

(5) Let $c(x,y)= log (exp(x)exp(y)exp(-x)exp(-y))$.  We have to show that
$c(V_n,V_m) \subseteq V_k$ when $k \leq N, k < n+m$.  Now if $|u|< 2^{-n}$ and $|v| < 2^{-m}$ then
$| c(u,v) - [u,v] |O(2^{-m-n-\min(m,n)})$, where $[u,v]$ is the Lie algebra bracket.   
This can be seen by looking at the power series expansion of $c$; it begins with $[u,v]$,
followed by higher order terms.  
 So the statement
holds for large enough $m,n$; by renormalizing (replacing $V_n$ by $V_{n+k}$)
we obtain the result.

Finally note that $U_{n+1} = \{u \in U_1: u^2 \in U_n \}$; since for $u=exp(v)$ we have 
$u^2=exp(2v)$ and
$u \in U_{n+1}$ iff $v \in V_{n+1}$ iff $2v \in V_n$ iff $u^2 \in U_n$.

Since $h \inv(U_2)$ is an $\inft$-definable set contained in the definable set $h \inv(U_1)$,
there exists a definable set $Y_1$ with $h \inv(U_2) \subseteq Y_1 \subseteq h \inv(U_1)$. 
Define   $Y_n$ inductively by: $Y_{n+1} = \{y \in Y_1: y^2 \in Y_n \} $. 
 It follows that $h \inv(U_{n+1}) \subseteq Y_n \subseteq h \inv(U_n)$.    (If $h(x) \in U_{n+1}$ then $h(x^2) \in U_n$; by induction $x^2 \in Y_{n-1}$; so $x \in Y_n$.
If $x \in Y_n$ then $x^2 \in Y_{n-1} \subseteq h \inv U_{n-1}$ so $h(x)^2  \in U_{n-1}$ and $h(x) \in U_n$.)  Clearly $Y_n=Y_n \inv$.  It follows from the intertwining of the $Y_n$ in the $h \inv U_n$ that 
$Y_{n+2}Y_{n+2} \subseteq Y_n$,   that $Y_n$ is contained in at most $c^2$ translates of $Y_{n+1}$, 
$a Y_{n+2} a \inv \subseteq Y_{n}$, and $[Y_n,Y_m] \subseteq Y_k$
whenever  $k+1 < n+m$.

Let $X_n=Y_{2n}$.  Then it is clear that (1-6) hold.  (7) follows from the fact that squaring is injective on $U_1$ (if
one chooses $U_1$ small enough.)
\eprf

\<{rem} \lbl{ap2}     (1)   Let $(X_n)$ be as in \lemref{ap1.1}.
Let   $(G,X)$ be a non-principal ultraproduct of $(G_0,X_n)$, and let 
   $\tG$ be the subgroup of $G$ generated by $X$. 
Then there exists a locally definable subgroup $\breve{G} \leq \tG$,   a connected, finite-dimensional 
  Lie group $U$, 
 and  
a homomorphism $h: \breve{G} \to U$ as in \thmref{ap1}, such that  in addition, 
$U$ is Abelian.   
\>{rem}  

\prf      Let $h,L,X_{n}$ be as in \ref{ap1.1}.  
  For $k \in \Zz$, for all $n \geq k$, define  $X_n[k]= X_{n+k}$.
  This carries over to the ultraproduct, so $X[k]$
is defined for all $k \in \Zz$, and has similar properties:  $X[k]X[k] \subseteq X[k+1]$.   Since the ultrafilter is nonprincipal, it concentrates on $n>k$, so by (5) of \ref{ap1.1} we have $[X,X[k]] \subseteq X[k+n-1] \subseteq X[k']$
for all $k' \in \Nn$.   Factoring out 
$\meet _{k}  hX[k]$ we obtain a commutative locally compact group.  As in \thmref{ap1}
we may replace it with a commutative Lie group.   
\eprf

\bigskip

We now deduce a version in the asymptotic setting.  Here we do not obtain an infinite chain, but the function $f$ serves to say that the length of the chain is arbitrarily large compared to $e,c,k$.   Say two sets are $e$-commensurable
if each is contained in the union of $\leq e$ cosets of the other.   Taking $\nu$ to be the counting measure, we obtain 
(a strengthening of)  \thmref{gc0}.

 \<{cor}\lbl{gc} Let $f: \Nn^2 \to \Nn$ be any function, and fix $k \in \Nn$.   Then there exist
$e^*,c^*,N \in \Nn$  such that the following holds.  

 Let $G$ be any  group, $X$ a subset, and assume there exists
a translation - invariant finitely additive real-valued measure $\nu$ on the definable subsets of 
$G$ contained in some power of $X X \inv$, with $\nu(X X \inv X) \leq k \nu(X)$.   

 Then there exist $e \leq e^*, c \leq c^*$ and  $0$-definable subsets 
$X_N \subseteq X_{N-1} \subseteq \cdots \subseteq X_1 $,    $N > f(e,c)$ 
  such that $X\inv X$ and $X_1$ are $e$-commensurable 
   and for $1 \leq m,n < N$ we have
\<{enumerate} 
\item  $X_{n}=X_{n} \inv$
\item   $X_{n+1}X_{n+1} \subseteq X_{n}$ 
\item $X_{n}$ is contained in $\leq c$ translates of $X_{n+1}$.
\item   $aX_{n+1}a \inv \subseteq X_{n}$ for $a \in X_{1}$. 
\item  $[X_{n},X_{m}] \subseteq X_{k}$ whenever $k \leq N$ and $k < n+m$.
In particular each $X_n$ is closed under the commutator bracket.
\item  $X_{n+1}=\{x \in X_{1}: x^4 \in X_n \}$ 
\>{enumerate}
\>{cor}  

\prf   Fix $f,k$.  We consider groups $G$ and subsets $X$ admitting a measure 
as above, with $\nu(X)=1, \nu(X X \inv X) \leq k$ (as we can always arrange by renormalizing $\nu$.)  
Consider integers $c$, and formulas $\phi$ of one free variable.  
Given $c,\phi$ and $X$, let
$X_1$ be the subset defined by $\phi$.   Let $ N=f(e,c)+1$, and
define $X_n$ using (6) for $2 \leq n \leq N$.  
Let us say that $(c,e,\phi)$   {\em works} for $X$ if properties (1-5) hold for
the sets $X_n$ defined in this way.

We will show that for some finite set  
$(c_1,e_1,\phi_1),\ldots,(c_n,e_n,\phi_n)$, for any $G$ and any  $k$-near-subgroup $X$ of $G$, some $(c_i,e_i,\phi_i)$ works for $X$.  Suppose this is false.  Then by the compactness theorem there exists $G$, a measure $\mu$
on the definable subsets of the group $\tG$ generated by $X$ with $\mu(X)=1$, and a definable subset $X$ of $G$ with $\mu(X X \inv X) \leq k \mu(X)$
such that no $(c,\phi)$ works.  But let $X_1$ be the definable set provided by \corref{ap1.1}.  Let $c$ be the integer of \corref{ap1.1} (3), and $e$ the number
of translates of $X_1$ needed to cover $X$.  Then 
   $(c,e,\phi)$ works for $X$ (indeed the $X_n$ have the required properties beyond
   any bound.)   This contradiction proves the the statement, and the theorem.
\eprf

\<{remark} \lbl{gc-rems} \rm \<{enumerate}
\item
Again we can also add (7):  if $x,y \in X_2$ and $x^2=y^2$ then $xy \inv \in X_N$.

\item  We could add that $X_1 \subseteq (X \inv X)^2$ (as we do in the statement of \thmref{gc0}), if we waive
the 0-definability of the $X_n$.  They remain definable over parameters from the given structure.  See \lemref{ap1-2}.

\item  This type of proof is always effective in the sense of G\"odel.  
Note that if   $f$ is recursive, then $e,c,N$ are automatically given by a recursive function (it suffices
to search for $e,c,N,X_1$ such that (1-6) hold.)  
\>{enumerate} 
\>{remark}
The sequence of subsets $X_n$ in \corref{gc} is recursively determined by $X_1$, via (6).  
The question is thus how to describe $X_1$.  Studying the proof of the fundamental 
theorems on locally compact groups should provide  detailed information; for now
we state what is clear a posteriori when they are treated as a black box.

Regarding $X_1$, we have:

\<{cor} \lbl{poly}  Fix $k \in \Nn$, and $f$ as above.   Then there exists $m$ and an algorithm that  
accepts as input   the multiplication table of a finite near-subgroup $X$ up to $(X \inv X)^{m}$, and yields the set $X_1$   in polynomial time.  \>{cor}

\prf  A formula in a logic with measure quantifiers $Q_\e$ can be computed in polynomial time.
\eprf

Note that we do not assume that $G$ itself is finite; and even  if finite, it is not 
available to the algorithm, beyond $(X \inv X)^m$.   Indeed the algorithm can be made to work for local groups.   

We can improve the $m$ to $3$ if, as in \remref{gc-rems} (2),
we use a formula with a parameter from $X$.  In this case the algorithm will first search for a parameter satisfying an
appropriate auxiliary formula, then compute $X_1$ using this parameter.  

\medskip

To illustrate \thmref{ap1-2}  we recover an easy version of a theorem of Freiman's (see \cite{t-v})  (generalized  to the non-commutative case).

\<{cor}  \lbl{freiman}  Fix $m,k$.  Then there exists $e=e(k,m)$ with the following property.
Let $G$ be a group of exponent $m$, i.e. $x^m=1$ for any  $x \in G$.  Let $X$ be a finite $k$-approximate
subgroup of $G$.  Then there exists a subgroup $S$ of $G$ such that $S,X$ are $e$-commensurable.  \>{cor}

\prf  By compactness it suffices to show that if $X$ is a near-subgroup of a group $G$ of bounded exponent,  then there exists a definable subgroup $S$ of $G$  such that $S,X$ are commensurable. 

By a theorem of Schur's (\cite{curtis-reiner} 36.14), a periodic subgroup $P$ of $GL_n(\Cc)$ has an abelian normal subgroup of finite index.
When the period is bounded, the abelian subgroup and hence $P$  must be finite.  If $L$ is a connected Lie group with center $Z$,
by considering the action of $L$ on its Lie algebra we see that $L/Z$ is linear.  Hence a periodic subgroup $P$ of $L$ of bounded period 
must be contained in $Z$ up to finite index, and again it follows that $P$ is finite.  

Thus the image of $\breve{G}$ in the Lie group $L$ associated to $X$ is finite.  Since this image is dense in $L$, and $L$ is connected, it follows that
$L$ is trivial.   The conclusion follows from \remref{lie0}.  
\eprf

Actually it is easy to see in the same way (using the fact that $L$ has no compact normal subgroups) that if $G$ is a periodic group and $X$ is a near subgroup,
then there exists a subgroup $S$ of $G$ such that $S,X$ are commensurable.  This does not extend to  families of  finite approximate groups, since without a uniform bound taking an ultraproduct will not preserve periodicity.

\>{section}

 \<{section}{Linear groups}  
 
 Up to this point the hypotheses in this paper were purely measure-theoretic,  at the top dimension as it were.  We will now look at lower dimensions as well.  Numerically this
 means that if a subset $Y$ of $X$ has about $c |X|^\alpha$ elements, we pay attention to $\alpha < 1$
 and not only to $c$ when $\alpha=1$.   Our main tool  is a  cardinality estimate due in its original  form to Larsen and Pink; in \cite{hw} it was presented as a dimension comparison lemma and slightly generalized in a number of directions; one of these will be needed here.   We will first define quasi-finite dimension in general, and specifically
 for ultraproducts of finite approximate subgroups.
 Then, assuming the group is linear (or indeed densely embedded in a group with a nice dimension theory) we show that the ambient dimension $\dim$ constrains strongly the quasi-finite dimension.
 Finally, knowing that the group looks sufficiently non-commutative by certain measures using quasi-finite dimension, we can conclude using the stabilizer that it is in fact definable.   
  
\ssec{The semi-group of dimensions}
Let $K$ be an ultraproduct of structures $K_i$ for some language $L$.  For each $i$ we consider, 
along with $K_i$, the counting measure on definable sets, as a map from the class of definable sets into $\Rr$.
Taking the ultraproduct of these maps as well, we obtain a map from the class of definable sets of
$K$ into the ultrapower $\Rr^*$  of $\Rr$.  This is a countably saturated real closed field.   For nonempty definable $X$
(represented by a sequence $X(K_i)$), we have a nonstandard real number $\log |X|$ (represented by the sequence $\log |X(K_i)|$.)  

Let $C$ be a  convex subgroup of $\Rr^*$.  We assume $C$ is  a countable union or a countable intersection of  definable subsets of $\Rr^*$.  
   Then $\Rr^*/C$ is an ordered  $\Qq$-vector space.
We define $\d(X)$   to be the image of $ \log |X| $   in   $\Rr^*/C$.  We view $\d$ as a (non-integral) dimension.  

{\bf Subadditivity:}  Let 
 $f: X \to Y$ be a definable map; assume $\d( f\inv(y)) \leq \a=a+C$ for each $y \in Y$, and $\d(Y) \leq \beta=b+C$.
 Then $\d(X) \leq \a + \beta$.   To see this, if $C= \union_n C_n$ is a countable union, we may take $C_1 \subset C_2 \subset \ldots$.  We have
 $\log |Y| \leq b+c$ with $c \in C_n$ for some $n$; and by compactness,
 $\log | f\inv(y)|  \leq a+ c'$, with $ c' \in C_{n'}$ for some $n'$.    It follows that $\log |X|  \leq a+b +c+c'$,
 so $\d(X) \leq \a+ \beta$.  If $C = \meet _n C_n$, then $\log |Y| -b  \in C_n$, and
  $\log | f\inv(y)| -  a \in C_n$ for each $n$; hence $\log|X| -a-b \in C_n$ for each $n$.
 
 We would like to extend the dimension to $\inft$-definable sets.  Fix $\delta_0 \in \Rr^*, \delta_0 > C$.  Let $V_0=V_0(\delta_0)$
be the group of elements $a \in \Rr^*/C$ such that $-n \d_0 +C \leq a \leq n \d_0 +C$ for some $n \in \Nn$.  
Let $V=V(\delta_0)$ be the set of  cuts of $V_0$, i.e. subsets $s \subset V_0$ that are nonempty, bounded above, and closed downwards.   This is a semi-group 
under set addition, linearly ordered by set inclusion.  $V_0$ embeds into $V$, by $a \mapsto \{v: v \leq a \}$.  
We identify $V_0$ with its image in $V$.    
Any subset of $V$ that is bounded below has a greatest lower  bound, namely the intersection.   We note that $V_0$ consists of invertible 
elements of $V$, and that it is semi-dense in $V$, in the sense that if $u< v \in V$ then there
exists $z \in V_0$ with $u <z \leq v$.  

  It will suffice for our purposes to use the intermediate subsemigroup $V_1$
consisting of infima of bounded countable subsets of $V_0$.

We could also form the  linearly ordered semi-group $V'$ of 
cuts in $V_0':  = \{a \in \Rr^*: \hbox{ for some }n \in \Nn, \ \ -n \d_0 \leq a \leq n \d_0  \}$.   The natural map $V_0' \to V_0$ maps cuts to cuts, and respects addition and $\leq$.  
$V$ can be identified with the ordered subsemigroup  of cuts $I \in V'$ with $C+I=I$.   Note that for a subset of $V$,
the infimum (in the sense of $V'$) lies in $V$ and agrees with the infimum in $V$. 

When $\a_n, \b_n$ are descending sequences of cuts,
 $\inf_n (\a_n + \b_n) = \inf_n \a_n + \inf_n \b_n$ holds in $V$.   By the above remark, it suffices to check this in $V'$.   The inequality $\geq$ is clear,
 since $\a_n + \b_n \geq \inf_n \a_n + \inf_n \b_n$ for each $n$.  For the other inequality, 
   let $\a_n' \in \a_n \m \a_{n+1}$,
 $\b_n' \in \b_n \m \b_{n+1}$.  Suppose $c \in \Rr^*$ and $c  \leq \inf_n (\a_n + \b_n)$.
 Then by countable saturation there exist $(\a,\b) \in \Rr^*$ with $\a \leq \a_n', \b \leq \b_n'$
 for each $n$, and $c \leq \a+\b$.  Hence $c \leq \inf_n \a_n + \inf_n \b_n$. 
 
 Let us also point out that if $\a < \a'$ and $\b < \b'$ are cuts, then $\a+\b< \a'+\b'$.  This holds
 for any semigroup of cuts in a dense linearly ordered group; to prove it we may consider the semigroup
 of all cuts.    Let $a \in \a' \m \a, b \in \b' \m b$.  Let $a^- = \{x: x<a \}$, and similarly $b^-$.  We have $a^-+b^- \leq a+b$,
 and $a^- + b^- \neq a+b$ since the cut  $a^- + b^-$ has no maximal point.  Thus
 $\a+\b \leq a^- + b^- < a+b \leq \a'+\b'$.  (One strict inequality and one weak inequality would not suffice for the same.)
 
  The multiplicative group $\Qq^{>0}$ acts on
$V_0$, and hence on $V$.  
 
\ssec{Quasi-finite dimension}

For an $\inft$-definable set $X$ define:
 $$\d(X) = \inf \d(D)$$
 where $D$ ranges over all definable sets containing $X$.   Note a  continuity property of the dimension:  
  If $X= \meet X_n$ with $X_1 \supset X_2 \supset \ldots$  $\inft$-definable,
  then $\d(X) = \inf_n \d(X_n)$.     
  
 The subadditivity  property holds  for $\inft$-definable sets $X$:  let $f$ be a definable map, let $\g \in V_1$, and assume
  $\d(f \inv(a) \meet X) \leq \g$ for all $a$.   Then  
 $\d(X) \leq \d(f(X)) + \g $.
   Indeed if $X = \meet X_n$ with $X_n$ a descending sequence of definable sets,
 then $f(X) = \meet_n f(X_n)$ by compactness (saturation); say $\g=\inf \g_k$; then   for each $k$, 
 for some $n(k) $, we have $\d( f \inv (a) \meet X_{n(k)}) \leq \g_k$, again by compactness.  So 
 $\d(X_{n(k)}) \leq \d(f(X_{n(k)})) + \g_k$.  Thus 
 $\inf_n \d(X_n) \leq  \inf_k \d(X_{n(k)}) \leq \inf_k \d(f(X_{n(k)})) + \inf_k \g_k  = \d(f(X)) + \g$.   
 
As a very special case of subadditivity, noting that $\d(F)=0$
for finite $F$, we have $\d(D_1 \union D_2) = \max  ( \d(D_1),\d(D_2) )$.

Also, $\d(D_1 \times D_2) = \d(D_1)+\d(D_2)$:  let $\mE_i$ be the family 
of definable sets containing $D_i$.  For any definable $E$ with 
$D_1 \times D_2 \subseteq E$ there exist (by compactness) $E_i \in \mE_i$
with $D_i \subseteq E_i$ ($i=1,2$) and $E_1 \times E_2 \subseteq E$.  Thus
$\d(D_1 \times D_2) = \inf_{E_1 \in \mE_1,E_2 \in \mE_2} \d(E_1 \times E_2)=
\inf _{E_1,E_2}\d(E_1)+\d(E_2) = \inf_{E_1} \d(E_1) + \inf_{E_2} \d(E_2) = \d(D_1)+\d(D_2)$.

If $X$ is $\inft$-definable over a set $A$, there exists a complete type
$P$ over $A$ containing $X$ with $\d(X)=\d(P)$.  To see this it suffices to
check that $X \m \union \{D: \d(D) < \d(X) \}$ is nonempty, since any type extending this will do.  By compactness it suffices to see that $X$ is not contained in a finite union 
of sets $D$ with $\d(D) < \d(X)$.  This is clear using $\d(D_1 \union D_2) = \max  \d(D_1),\d(D_2)$.

If $\d(X) \in V_0$, we say that $X$ has strict quasi-finite dimension.  Note in this case, by saturation of $\Rr^*$, that if $\d(X) = \inf_{n \in \Nn} \a_n$ then $\d(X)=\a_n$ for large enough $n$. 
 
\ssec{Examples}  The best behaved case is of  totally categorical theories, i.e. theories $T$ with a unique model in each power; a basic example is the theory of vector spaces
over a fixed finite field.  In this case, any finite subset $T_0$ of $T$ has finite models, described by a single integer parameter; the   cardinality is precisely given by a polynomial $P$ in this dimension
parameter.   
In this case, regardless of the choice of convex subgroup, the quasi-finite dimension equals the degree of the Zilber polynomial (times a fixed scalar), and recovers the Morley dimension.    
This was an essential ingredient of Zilber's theory of totally categorical structures; see \cite{CH} for generalizations and converses. 

In general, the two most natural choices for a convex subgroup $C$ are the smallest nontrivial convex subgroup, the convex hull $C_{min}$ of $\Zz$;
and the largest convex subgroup $C_{max}$ with $\delta_0 \notin C_{max}$.  

$C_{max}$ is a countable {\em intersection} of definable subsets of $\Rr^*$.   The corresponding group of dimensions is canonically isomorphic to $\Rr$
(with $\delta_0$ mapping to $1$.)    
Each dimension $\alpha$ induces an ideal $I_\alpha$ (\exref{dim}), which is not S1.  Asymptotically, a set $X$ represented by a sequence of finite sets $X_i$
has the same dimension as $Y \subset X$ (represented by $Y_i$) if for any $\epsilon$, for almost all $i$, $|Y_i| \geq |X_i|^{1-\epsilon}$.  

$C_{min}$ is a countable {\em union} of definable subsets of $\Rr^*$.  The group of dimensions is more complicated; but when $\d(X)=\a$< 
the ideal of lower-dimensional subsets of $X$ is an S1-ideal.   Here $Y \subset X$ has the same dimension as $X$ if for some $k$, 
$|X_i| \leq k |Y_i|$ for almost all $i$.  

\ssec{Minimality}

Now assume each $K_i$ is a field, possibly with additional structure.  There will
be no loss of generality in assuming that $K_i$ is algebraically closed.   Let $K$ be an ultraproduct of the structures $K_i$.  
 Constructible
sets and varieties will be assumed to be defined over $K$.  Here the words ``constructible"
means:  definable in $K=K^{alg}$ as a field, whereas ``definable" means: definable in $(K,\cdots)$ as an $L$-structure.
For a constructible set $S$, asides from the pseudo-finite dimension $\d$ constructed above, we have the  dimension in the sense of
algebraic varieties.  This can be defined as the Morley rank of the set, viewed as definable in $(K,+,\cdot)$.  Or it can be defined
as the dimension of the Zariski closure of $S$, as in \cite{weil}; see \cite{poizat2}.  

Let $G$ be a simple algebraic group over the ultraproduct $K$.   $G$ can be viewed as a group subvariety of the group
$GL_n$ of invertible matrices.  We write $G$ when we think of the defining equations, and $G(K)$ when we think of the set of points of $K$.

Let $\G_0$ be a Zariski dense subset of $G(K)$.  
Consider the functions $F_c(x,y) = c x \inv c \inv y$, $c \in \G_0$.
Any subvariety $H$  of $G(K)$ closed under all the $F_c$ must be a group subvariety of $G$, normalized by $\G_0$, hence by the Zariski closure
of this group, i.e. by $G$.  Since $G$ is simple, we must have $H=1$ or  $H=G$. 

It follows that if $Y,Z$ are   constructible subsets of $G$, defined 
over a subfield $A$ of $K$, and 
$0 < \dim(Y) \leq \dim(Z) < \dim(G)$, then $\dim(F_c(Y \times Z)) > \dim(Z)$ for some
$c \in \G_0$.  Moreover, let $Y \times' Z = Y \times Z \m \union_j W_j$, where $W_j$
ranges over all $A$-definable constructible subsets $W$ of $Y \times Z$ with $\dim(W)<\dim(Y)+\dim(Z)$.  (This is the same, for the theory ACF of algebraically closed fields,
as the product $\times_{nf}$ encountered in the proof of \thmref{stabilizer}.  Note that this product depends on the base set $A$, but
as $A$ will be fixed we will omit it from the notation.)  
Then $\dim(F_c(Y\times' Z))> \dim(Z)$ for some $c$.  This is a typical application of Zilber's stabilizer, and in itself an instance of the ``sum-product"   phenomenon  in a constructible setting: we may assume $Z$ is irreducible.   If $\dim(F_c(Y\times' Z)) = \dim(Z)$, we find that $Y$ 
and all $\G_0$- conjugates of $Y$ are contained in finitely many cosets of the Zilber stabilizer $H = \{y:  \dim(yZ \triangle Z)< \dim(Z) \}$,
 up to smaller dimension.  But then $\meet_{x \in \G_0} x \inv H x$ (a finite intersection of conjugates of $H$)
is closed under all the $F_c$, so by the first paragraph it equals $G$, i.e. 
contradicts the previous paragraph.   

This property of $(G,F_c)_{c \in \G_0}$  is referred to as   {\em minimality}.  See    \cite{hw}, Example 2 for details and generalizations.

\ssec{The dimension inequality}

Let $\G \subseteq X$ be an $\inft$-definable subgroup, of strict quasi-finite
dimension $\d(\G)=\d(X)$; and with $\G_0 \leq \G$.   Let  $\delta_0=\d(\G)$, $V=V(\delta_0)$. For $n \in \Qq$ and $v \in V$, $n v$ is defined; we write ${\gamma_0} n $ for $n {\gamma_0}$.    Let ${\gamma_0} = \d(\G) / \dim(G)$.  
  
 For a constructible $Z \subseteq G(K)^n$, define $\delta_\G(Z) = \delta( Z \meet \G^n)$.
For any $W \subset G(K)^n$, let $\dim(W)$
denote the dimension of the Zariski closure of $W$.

\<{prop} \lbl{hw1}  For any constructible $Z \subseteq G^n$, we have $\delta_\G(Z)  \leq   {\gamma_0} \dim(Z)   $.  \>{prop}  

\prf   This is a special case of Corollary 1.12 of \cite{hw},  as generalized in Remark 1.11.
We give the proof in the present case, for $Z \subseteq G$.   For $G^n$ see the remark below.

  Let
$W$ be an $\inft$-definable subset (over $A$)  of $\G^n$.  Call $W$ {\em unbalanced} if 
$ \delta(W) > {\gamma_0}  \dim(W)$.   There exists a complete
type $W' \subset W$ defined over $A$ with $\d_\G(W')=\d_\G(W)$.   As $\dim(W')\leq \dim(W)$,
$W'$ is unbalanced if $W$ is. 
 
We must show that no unbalanced sets exist.  Otherwise, let $Y,Z$ be unbalanced $\inft$-definable sets  with $\dim(Y)$ minimal, and $\dim(Z)$ maximal possible.    Clearly
$0 < \dim(Y) \leq \dim(Z) < \dim(G)$.   
Say $Y,Z,c$
are defined over the countable  $A \leq K$.  By the above, we may take $Y,Z$ to be complete types over $A$.  
Form $Y \times' Z$.  By minimality of $(G,F_c)_{c\in X}$
there exists $c \in X$ with $\dim( F_c(Y \times' Z))   > \dim(Z)$. 

   We note first that $Y \times' Z$ is balanced:
 Let $f$ be the restriction of $F_c$ to $Y \times' Z$.  Then since
$Y \times' Z$ implies a complete quantifier-free type over $A$ in the language of fields, the fiber dimension $\dim f \inv(a)$ is constant (=${b}$) for $a \in  F_c(Y \times' Z)$, and 
from $\dim(f(Y \times' Z))> \dim(Z)$ it follows that 
 ${b}  < \dim(Y)$.   So the fibers are not unbalanced, i.e.
 $\d_\G( f\inv(a)) \leq {b} {\gamma_0}$.   On the other hand since $\dim(f(Y \times' Z))>\dim(Z)$, $f(Y \times' Z)$ is not unbalanced either, so $\d_\G(f(Y \times' Z)) \leq \dim( f(Y \times' Z)) {\gamma_0}$.   By subadditivity we obtain $  \d_\G(  Y \times' Z) \leq( {b} + \dim(f(Y \times' Z))) {\gamma_0} = \dim(Y \times' Z) {\gamma_0}$.  
  
Now there exists a complete type $Q$ over $A$ with $Q \subseteq Y \times Z$ and
$\d_\G(Q) = \d_\G(Y)+\d_\G(Z)$.
I claim that $Q=  Y \times' Z$ (formed over $A$).  For if $Q$ is any
other type,  the fibers $Q_a=\{w: (w,a) \in Q\}$ have dimension $\dim(Q_a) = b' < \dim(Y)$ for $a \in Z$
 (the dimension is constant on $Z$   since $Z$ is a complete type). 
  By minimality of $\dim(Y)$ we have $\d_\G(Q_a) \leq b' {\g_0}$.  By subadditivity 
  it follows that 
  $$\d_\G(Y)+\d_\G(Z) = \d_\G(Q) \leq b' {\g_0} + \d_\G(Z) < \dim(Y) {\g_0} + \d_\G(Z) \leq \d_\G(Y)+\d_\G(Z)$$
  So $Q=Y \times' Z$ .

Since $\g_0 \dim(Y) < \d_\G(Y)$ and $\g_0 \dim(Z) < \d_\G(Z)$, we have
$\g_0 \dim(Y \times' Z) = \g_0 (\dim(Y) + \dim(Z)) < \d_\G(Y)+ \d_\G(Z) = \d_\G(Y' \times Z)$.  So $Y \times' Z$ is unbalanced.  A contradiction. \eprf

\<{remark} \lbl{hw1.1}  \rm  Let  $f: X \to X'$ be a constructible map,  $\G$ an $\inft$-definable subset of $X$,
$\G' = f(\G)$.  
 If the inequality of \propref{hw1} holds for $\G' \subseteq X'$ and for each fiber $f \inv(a) \subseteq X$, $a \in \G'$, all with the same value of  ${\g_0}$, then it holds for $\G \subseteq X$.  This is an easy consequence of subadditivity and definability of Zariski dimension, cf. \cite{hw}.  
 \>{remark}

Recall that a morphism $f: U \to V$ of algebraic varieties is {\em dominant} if there exists no proper subvariety
$V'$ of $V$ such that the image of $U$, over any field, is contained in $V'$.
   
\<{lem} \lbl{hw2} Let $U$ be a Zariski open subset of $G^m$.  Let $f: U \to W \subseteq G^n$ be a dominant morphism
of varieties.  
 Then $\d( f(U \meet \G^m)) =  \dim(W) {\gamma_0}$.
\>{lem}  

\prf We first show that if $U$ is Zariski open in $G^m$, then $\d_\G(U) = m \dim(G) {\gamma_0}$. 
We have $\d(\G^m)=m \d_\G(G)$,   On the other hand
if $V$ is a proper Zariski closed subset of $\G^m$, then $\dim(V) \leq m \dim(G) -1$,
and by \propref{hw1} $\d(V \meet \G^n) \leq \dim(V) {\gamma_0} < m \d_\G(G)$.  It follows that
$\d(\G^m \m V) = m \d_\G(G) = m \dim(G) {\gamma_0}$.

 There   exists a relatively Zariski open $W' \subseteq W$, $\dim(W')=\dim(W)$,
such that $\dim f \inv(b)$ is constant for $b \in W'$.  Replacing $W$ by $W'$
and $U$ by $f \inv(W')$, we may assume  $\dim f \inv(b)=d$ is constant for $b \in W$.
 So $\dim(U) = d+\dim(W)$.
By \propref{hw1}, for any $b \in W \meet \G^n$, we have $\d_\G( f \inv(b)) \leq d {\gamma_0}$.
Hence if $\d(f(U \meet \G^m)) =\gamma < \dim(W) {\gamma_0}$, then by subadditivity
of $\d$ we have  $\d(U \meet \G^m) \leq \gamma+ d {\gamma_0} < (\dim(W)+d){\gamma_0} =  \dim(U) {\gamma_0}$;
this contradicts the first paragraph.  Note that adding the invertible element $d {\gamma_0}$ preserves strict inequalities.
    \eprf

In case $f( U \meet \G^m) \subseteq \G^n$, it follows that $\d_\G(W) =  \dim(W) {\gamma_0}$.   The proof shows more generally that the class of subvarieties $U$ of $G^m$ satisfying 
$\d_\G(U') = \dim(U) {\gamma_0}$ for all Zariski dense open $U'$, is closed under forward images of such morphisms.

Note that an $\inft$-definable subgroup $\G$ of $G(K)$ has strict quasi-finite dimension iff for
some $\bigvee$-definable $\tG$ containing $\G$, $\tG/\G$ is bounded. 

\ssec{\em From now on we assume $C$ is the convex hull of $\Rr$ in $\Rr^*$, a $\bigvee$-definable convex subgroup.}

For $Y \subseteq X$, let $\mu(Y)$ be the unique real number $r$ such that for any rational $\a$, 
$\a |X|> |Y|$ if $\a >r$ and $\a |X| < |Y|$ if $\a < r$.   Then $\mu$   is a definable measure on definable subsets of $X$, and we have $\mu(Y) > 0$ iff $\d(Y)= \d(\Gamma)=\d(X)$.

\<{prop} \lbl{1.3p}    Let $\G$ be a Zariski dense $\inft$-definable subgroup of  $G(K)$, $G$ a semisimple algebraic group over a $K$.  Assume $\G$ has strict quasi-finite dimension.  
Then $\G$ is definable.    
\>{prop}  

\prf   Assume first that $G$ is simple.  Let $K^a$ be the algebraic closure of $K$.  
Let $\G_0 \leq \G$ be any Zariski dense set of points, so that the previous lemmas apply.
   Since $G$ is a simple group,  any non-central conjugacy class $C$ of $G(K^{a})$ generates $G$ in a finite number $d \leq 2 \dim(G)$ of steps.  Thus for any noncentral $b$, the  morphism of varieties $f_b: G^d \to G$,
$f(x_1,\ldots,x_d,b) = x_1 \inv bx_1 x_2 \inv b x_2 \cdots x_d \inv b x_d$
is surjective on $K^a$-points.  By \lemref{hw2}, $\d(f_b(\G^d)) \geq \dim(G) {\gamma_0} = \d(G)$.  
 Let $X$ be a definable set containing $\G$
with $\d(X)=\d(\G)$, and let $\tG$ be the group generated by $X$.  Then $\tG/\G$ is bounded.
Let $S \subseteq \G$ be an $\inft$-definable normal subgroup of $\tG$ with $\tG/S$ bounded (\lemref{normal0}), and choose
a noncentral $b \in S$.   Let $Y$ be the definable set $Y=f_b(X^d)$.   Since $S$ is normal, $Y \subseteq S$.  We have
  $\d(Y) \geq \d(f_b(\G^d))=\d(G)$, so $\mu(Y) >0$.   Hence $S$ contains a bounded finite
  number of disjoint translates $s_i Y$ of $Y$, so any $s \in S$ lies in $s_i Y Y \inv$ for some $i$
  (Ruzsa's argument.)  Hence  $S = \union_i YY \inv$ is definable.     Since $\G/S$ is bounded and $\inft$-definable,
it must be finite, so $\G$ is definable too. 

When $G$ is semisimple, we proceed by induction on $\dim(G)$.  Let $N$ be a normal algebraic subgroup, 
$\pi: G \to G/N$ the natural homomorphism.  
Since  $\G$
has  strict quasi-finite dimension, for some $\bigvee$-definable $\tG$ we have $\tG/\G$ bounded.
It follows that $(N \meet \tG) / (N \meet \G)$ and $\pi(\tG)/\pi(\G)$ are bounded, so $N \meet \G$
and $\pi(\G)$ have strict quasi-finite dimension; by induction they are definable.  
By \remref{towers},  $\G$ is definable. \eprf

\prf[Proof of \thmref{linear}, and Corollary.]   Suppose not.  Then there exists an ultraproduct $(K,X)$ of $(K_i,X_i)$
such that for no definable subgroup $H$ of $G(K)$ do we have $H \subseteq (X \inv X)^2$ and
$X$ contained in finitely many cosets of $H$.  However $X$ is a near-subgroup of $G(K)$.  Let $\tG$
be the subgroup generated by $X$.  By \thmref{stabilizer}  there
exists an $\inft$-definable group $\G \subseteq (X \inv X)^2$, normal in $\tG$, with $X/\G$ bounded.  By \propref{1.3p},
$\G$ is definable.  By compactness, $X/ \G$ is finite.  This contradiction proves the theorem.

The corollary easily follows, and can also be quickly proved directly in the same way:   if it fails,  we obtain an ultraproduct  $(K,X)$ with $X$ an infinite near-subgroup, generating a $\bigvee$-definable group $\tG$ strictly bigger than $(X \inv X)^2$, and
  such that no infinite definable proper subgroup of $\tG$ is normalized by $X$.  Let $\G$ be as above.  Again
  $\G$ is definable, hence (being normalized by $X$) either $\G=\tG$ or $\G$ is finite.  If $\G$ is finite then since $X/ \G$  is bounded it is finite, contradicting the assumption that $X$ is infinite.  If $\G=\tG$ then since $\G \subseteq (X \inv X)^2$
  we must have $\tG=  (X \inv X)^2$, again a contradiction.

\eprf 

Similarly we can obtain $|S|/|X \inv X| \geq .9$ in \corref{1.4}; otherwise we obtain $(K,X)$ as above and also a measure
$\mu$ on $\tG$ with no infinite definable subgroup $H$, contained in $\tG$ and normalized by $X$,
 satisfying $\mu(H) / \mu(X \inv X) \geq .99$.   But again $\G$ is definable, and by \thmref{stabilizer},
 $\G \m X \inv X $ is contained in a union of non-$\mu$-wide sets; by saturation and definability
 of $\G$ is is contained in finitely many such sets, so $\mu(\G \m X \inv X) =0$; a contradiction.  One can
 also get $X_i X_i \inv X_i=S_i$ from the fact that  $q q \inv q $ is a coset of $S$ in \thmref{stabilizer}, and that $S_i$ has no subgroups of bounded index.   I noted this stronger statement after Laci Pyber pointed out that the statement of \corref{1.4} implies  $X_i X_i \inv X_i = S_i$, using \cite{n-p}.   
 
One can immediately deduce a version for arbitrary linear groups:  

 \<{cor} \lbl{linear-solv}   Let  $k,n\in \Nn$.  Then there exist $k' \in \Nn$,  such that if $X$ is a 
$k$-approximate subgroup of $GL_n(K)$ for some field $K$, then 
there exist algebraic subgroups $H \leq G$ of $GL_n$ with $H$ solvable and normal in $G$, and a subgroup $\Delta$ of $G$ (normalized by $X$) with $\Delta \subseteq (X \inv X)^2 H$
and such that $X$ is contained in $\leq k'$   cosets of $\Delta$.  \>{cor}

The groups $H,N$ are defined by polynomial equations in the matrix entries;
these equations can be taken to have degree bounded by a function of $k,n$ 
alone.

Jordan has shown that finite subgroups of linear groups are bounded, up to an Abelian part, provided they   contain no nontrivial unipotent elements.    (Jordan's beautiful proof occupies some 13 pages of   \cite{jordan}.  
\cite{curtis-reiner} contains a different proof in characteristic $0$, due to Frobenius.)  We may now extend this
to say that approximate subgroups of connected Lie groups are bounded, up to a (connected, closed) solvable subgroup. 

\<{cor} \lbl{lie-solvable}  Let  $k \in \Nn$, and let $L$ be a connected Lie group of dimension $d$.  Then there exist $k'' \in \Nn$,  such that if $X$ is a 
(finite) $k$-approximate subgroup of  $L$,   then 
there exist  a $d+ 2$- solvable subgroup $S$ of $L$ 
 such that $X$ is contained in $\leq k''$ cosets of $S$.   
\>{cor}  

\prf     
Let $X$ be a $k$-approximate subgroup of  $L$
Assume first that  $L$ embeds into $GL_d(\Rr)$.  In this case, let $H,G$ be the subgroups provided by \corref{linear-solv}; so $H$ is $d$-solvable.
So   $X$ is contained in boundedly many cosets of  a   subgroup $\Delta$ of $G/H$, with $\Delta/H$
finite (as $X$ is finite.)    By \cite{jordan}, $\Delta/H$ contains a normal Abelian subgroup $S/H$ of bounded index.  Then $S$
is $d+1$- solvable, and $X$ is contained in boundedly many cosets of $S$.

 In general, let $Z$ be the center of $L$.   Then $L/Z$  acts faithfully on the Lie algebra of $L$ by conjugation, so it embeds into $GL_d(\Rr)$.  By
 the linear case, 
   If $X$ is a $k$-approximate subgroup of  $L$, then the image of $X$ in $L/Z$
is a $k$-approximate subgroup of $L/Z$, so by the linear case it is contained in boundedly many cosets of a solvable subgroup $S/Z$.
The pullback $S$ of $S/Z$ to $L$ is $d+2$-solvable, and finitely many cosets of $S$ cover $X$.    \eprf

\noindent {\bf Remarks.}

\<{enumerate} 
\item Once the definability of $\G$ is established, it is known to be definable in the field
 language, possibly expanded by an automorphism, and indeed to be
  a simple  group of  (possibly twisted) Lie type; see \cite{larsen-pink}.  

\item  We could also deduce \thmref{linear} from \corref{anti-nil}; the proof of \propref{1.3p}, together
with saturation, shows that   for some $m \in \Nn$ we have $\mu(C_b^d) \geq 1/m$ for all non-central $b$.
It follows that  with probability very close to $1$ (in $b$),   $\mu(C_b^d) \geq 1/m$; so the hypothesis of \thmref{anti-nil} holds.

\item  Using another direction of generalization taken in \cite{hw},  results of this section are valid for near-subgroups of   groups $G$ of finite Morley rank, in place of algebraic groups.  

\>{enumerate}

We further remark that  Theorem 1.1 of \cite{bkt} in the sum-product setting, as well as the theorem of \cite{helfgott} for 
subsets of $SL_2(\Ff_p)$, can be put in the framework of \propref{1.3p}  if one takes the    
$C=C_{max}$ to be the largest convex subgroup of $\Rr^*$ not containing $\delta_0$ (in place of the smallest
nonzero convex subgroup, as we took it to be.)   \footnote{In fact, after these lines were written, Breuillard, Green and Tao essentially took this route; 
using a beautiful analysis of  the geometry of tori, 
continuing a line started in \cite{jordan}, they obtain an {\em effective, polynomial} version of \thmref{linear}.  See \cite{bgt}.}

\>{section}

\<{section}{Uniform definability of the topology} \lbl{simple}
 
We prove a stronger version of the stabilizer theorem for arbitrary S1-ideals on 
Ind-definable group, with more uniform
control of the topology of the Lie group.  It follows that the Lie group associated to
a near-subgroup is always associated already to the reduct to a finite sublanguage.
Stronger uniformity statements in this direction may give
a more powerful means   for finitization of results about saturated
models.

\<{remark}  Let $T$ be a simple theory, or a NIP theory.  Then the forking ideal is an S1-ideal.  \>{remark}

\prf  Let $(a_i)$ be an $A$-indiscernible sequence, and suppose $\phi(x,a_i)$ does not fork over $A$.  We have to show that $\phi(x,a_i) \wedge \phi(x,a_j)$ does not fork over $A$, for some $i \neq j$.  

Simple case:
 the $a_i$ are independent over some $M$.  Let $c_i$ be such that   $\phi(c_i,a_i)$ with $c_i,a_i$ independent; choose $c_i$ so that    $c_i,M$ are independent over $a_i$.  Then
$c_i,Ma_i$ are independent over $A$.  The sequence $(a_i)$ could be taken to be long; 
by refining it we can assume that $tp(a_i/M)$ is constant.  By 3-amalgamation we   can find $c$ independent over $M$ from
$(a_i)_i$, with $tp(c,a_i/M)=tp(c_i,a_i/M)$.  Since $tp(c_i/M)=tp(c/M)$, $c,M$ are independent over $A$, so $c$ is independent from $a_1,a_2$ over $A$.
Hence $tp(c'/a_1a_2)$ does not fork  $A$.

NIP case:  Let $q_i$ be a global type with $\phi(x,a_i) \in q_i$, such that $q_i$ does not fork over 
$A$.  Let $M$ be a model containing $A$.  Then $q_i$ does not fork over $M$.  
So $q_i$ is $M$-invariant.  There are few choices for $M$-invariant types, so $q_i=q_j$
for some $i \neq j$.  Since $q_i$ does not fork over $A$, $\phi(x,a_i) \wedge \phi(x,a_j)$ 
does not fork over $A$.
\eprf

 Let $X$ be a topological space, $p \in X$.  
 We say a collection $C$ of sets strongly generates the topology at $p$ if $p$ is in the interior of each set in $C$, and any open neighborhood of $p$ contains some element of $C$.    
 
 Let $M$ be a Riemannian manifold.  Let $\rho(p,q)$ denote  the Riemannian distance, and $B(p,r)$ (respectively
 $\bar{B}(p,r)$)   the open (resp. closed) ball of radius $r$.   A {\em geodesic ball} around $p$ is 
 the image under the exponential map $exp_p$ of a ball $b$ of radius $r$ around $0$ in the tangent space to $p$, where
 $r$ is small enough that   $exp_p$ is a diffeomorphism.  We have $\rho(p,exp_p(v)) = |v|$ if $v \in b$
 (\cite{lee}, Proposition 6.10, p. 105.)
  A subset $U$ is called {\em convex} if for each $p,q \in U$ there is a unique
 geodesic $x$ from $p$ to $q$ of length $\rho(p,q)$, contained entirely in $U$.  Any point has a convex neighborhood
 (\cite{lee}, 6-4, p.112).

\<{lem} \lbl{unif2}  Let $M$ be a Riemannian manifold, $G$ a topological group  acting isometrically and transitively on $M$ (the action $G \times M \to M$ is assumed continuous.)    Let $B(p_0,r)$ be a geodesic ball of $M$, contained in a convex set $W$.   Assume there exists a compact $Y \subseteq G$ such that if 
$x,x' \in B(p_0,r)$ then for some $g \in Y$ we have $gx=x'$  and $\rho(x,g^2x ) = 2 \rho(x,x')$.
  
  Let $U$ be any open set  of diameter $< r$.   Let $C$ be   
    the collection of neighborhoods of $p_0$ 
of the form $cl(g_1U \meet g_2U)$.  Then $C$ strongly generates the topology at $p_0$.  
 \>{lem}

\prf   It suffices to show that there are nonempty sets of the form $g_1U \meet g_2U$,
of arbitrarily small diameter.  For  then by translation we may take these sets to contain $p_0$, and their closures will
still have small diameter, and will strongly generate   the topology at $p_0$.  

Let  $\cU$ be the closure of $U$, and let $\d<r$ be the diameter of $U$.   

Find $p_n ,q_n \in U$ with $\rho(p_n,q_n) \geq \d-1/n$; and find $g_n \in G$
with $g_n p_n = q_n$ and   $\rho(p_n,g_n^2 p_n)= 2 \rho(p_n,q_n)$.  By 
 assumption, we may choose $g_n$ in a compact set; and all $p_n,q_n$ lie within
a compact set (a closed ball of radius $r$).
    Refining the sequence $(p_n,q_n,g_n)$, we may thus assume it converges to a point $(p,q,g) \in \cU^2 \times G$; and we have $\rho(p,q)=\d$, $\rho(p,g^2p)=2 \rho(p,q) = 2 \d$.  
 It follows from uniqueness of the minimizing geodesic between $p$ and $g^2p$ that
$\bar{B}(p,\d) \meet \bar{B}(g^2p,\d) = \{q\}$.  By definition of $\d$ we have $\cU \subseteq \bar{B}(x,\d)$
for any $x \in \cU$.  In particular, 
$\cU \subseteq \bar{B}(p,\d)$, and $\cU \subseteq \bar{B}(q,\d)$.  From
the latter we obtain: $g \cU \subseteq \bar{B}(g^2p, \d)$.  So
$\cU \meet g \cU = \{q\}$.  

The set $U \meet g_nU$ is nonempty, since $q_n \in U \meet g_nU$.  It remains only
to show that the diameter of $U \meet g_nU$ approaches $0$ as $n \to \infty$.    

Suppose otherwise; then  there exist $\g>0$,  and $a_n,b_n \in U \meet g_nU$
such that $\rho(a_n,b_n) \geq \g$ for infinitely many $n$.    We can refine the sequence again to assume
$a_n \to a, b_n \to b$; we have $\rho(a,b) \geq \g$ so $a \neq b$, and
$a,b \in \cU \meet g\cU$.  But we have seen that $\cU \meet g \cU$ consists
of a single point; a contradiction.    
\eprf

The hypothesis of \lemref{unif2} are satisfied when $G$ is a Lie group, acting on itself
by left translation,   $M$ is $G$ with a left invariant Riemannian metric, and $B(p_0,2r)$ is a geodesic ball.
For then $Y = \bar{B}(p_0,r)\bar{B}(p_0,r) \inv$ is compact.  For $x,x' \in  B(p_0,r)$,
let $g=x' x \inv$, $h=x \inv g x = x \inv x'$, and let $|u|=\rho(1,u)$.
    Then $gx=x'$.  We have $\rho(x,x')=\rho(x',g^2x)=\rho(1,x \inv g x)=|h|$, $\rho(x,g^2x) = \rho(1, x \inv g^2 x)=|h^2|$, so we have to show that $|h^2| = 2|h|$.  We have $h=exp(v)$ for some $v$, where $exp$ is the 
  the exponential map at $1$, $h^2=exp(2v)$, and $|h^2|=|2v|=2|v| = 2|h|$.    

\<{cor}[Stabilizer theorem] \lbl{stabilizer-full}  Let $X$ be a near-subgoup of $G$. 
\<{itemize}
\item 
There exist  a $\bigvee$-definable $\breve{G}$ and an $\inft$-definable normal subgroup $\G \subseteq \breve{G}$, both defined
without parameters,
such that 
$\breve{G}/ \G$ is bounded;   and any definable $D$ with $\G \leq D \leq \breve{G}$ is 
commensurable to $X \inv X$.
\item  There exist a connected Lie group $L$ and a homomorphism $\pi: \breve{G} \to L$ with dense image,
and kernel $\G$. 
 If $D$ is a definable subset of $G$, write $\pi D$ for the closure of $\pi(D)$.  
$\pi$ intertwines   the definable sets containing $\G$, 
contained in $\breve{G}$ with the   compact neighborhoods of $L$.
\item There exist a uniformly definable family of definable sets $D_a$, and a definable set $E$, 
with $\bdr \pi(D_a) \meet \pi E \subset int(\pi(E)) $ 
such that the neighborhoods of $1$ of the form    $\pi E  \m \pi D_a$
 generate the topology of $L$  at $1$.
\>{itemize} 
\>{cor}

\prf   The first two parts follow from \thmref{ap1}.   

   There remains to prove the uniform generation of the topology of $L$.   Fix a left-invariant Riemannian metric on $L$, and view $M=L$ as a Riemannian manifold.  Let $p_0=1$ and let $r$ be as in \lemref{unif2};
   renormalizing, we may assume $r=4$.  
      Write $B_s$ for $B(1,s)$.   By the above remark there exists a definable $E$ with 
 $\bar{B}_5 \subseteq \pi E \subseteq B_6$.  Similarly
there exists a definable $D$ such that $\pi D$ contains $\bar{B}_9 \m {B}_{2}$ (so $\bdr(\pi(D)) \meet E \subset  int(E)$) and is  disjoint from $\bar{B}_{1}$.
 Then $U=\pi(E) \m \pi(D) = B_{2} \m \pi(D)$ is an open neighborhood of $1$.    By \lemref{unif2}, there exists $g,g' \in B_3$ with 
$g U \meet g'U$ of arbitrarily small diameter, and containing $1$.  We compute 
$U \meet gU = (\pi(E) \meet \pi(gE)) \m \pi(D \union gD) = B_2 \m  \pi(D \union gD) $ = $\pi(E) \m \pi(D \union gD)$.  Similarly for $gU \meet g'U$.  
 The uniformly definable family is   the family of unions $D \union gD$.
\eprf

 \ssec{The locally compact Lascar group}
 Let $T$ be a theory, $\Uu$ a universal domain, $\tE$ a $\bigvee$-definable equivalence relation,  $\Sigma$ an $\inft$-definable equivalence relation, such that $\Sigma$ implies $\tE$.   Let $P$ be a complete type.  
  Let $\ta$ be a class of $\tE$ restricted to $P$,  such that $\tau=\ta / \Sigma$ is bounded.    Let $\pi: \ta \to \ta / \Sigma$
  be the quotient map. 
  Then $\ta/\Sigma$ admits a natural  locally compact topology, generated by the complements of the images $\pi(D)$ of definable sets.  $G=Aut(\Uu/\ta)$ acts on $\tau$.  Let $K$ be the kernel of this action, and $L=G/K$.  Then $L$ admits a natural locally compact group structure; we call it the locally compact Lascar group of $(\ta,\Sigma)$.  
 
 We have transposed from definable groups (as in  \thmref{stabilizer-full}) to automorphism groups.  In both cases,
 the set of conjugates of a definable set lie in a uniformly definable family.  We will use this in \lemref{sop1}
 below.

 \ssec{The compact Lascar group}  
 So far, the case where $\ta$ is a definable set and $L$ is compact has been useful.  For simplicity, 
 we too will restrict to this case in the statement below. {\em  For the rest of this section we assume $\tE$ is the indiscrete 
 equivalence relation, so $\ta  = P$ and $ \tau=P/ \Sigma$ is compact.  } 
  We do not expect any trouble in generalizing
 to the locally compact case.

  \<{lem} \lbl{sop1}    Let $L'=L/N$ be a finite dimensional quotient of $L$, so $N$ is a compact normal subgroup and $L'$ is a compact Lie group.  For large enough $k$, $L'$ has a regular orbit on  $ \tau^k /N=P^k/N$.   Let $\tau'$ be such an orbit.
   There exists a uniformly definable family of definable sets $D_a$, such that the sets
  $\tau' \m  \pi(D_a)$ strongly generate the topology on $\tau'$ at every point.
  \>{lem}  
  
  \prf    For $x \in \tau$, let $S_x$ be the stabilizer of $x$.   Let $\Xi$ be the set of finite subsets of $\tau$.
  For $u \in \Xi$, let $S_u = \meet_{x \in u} S_x$.   We have
  $\meet_{x \in \tau} S_x = K$.   
 Since $N$ is compact, $\meet_{u \in \Xi} S_u N = KN =N$.  (Let $a \in \meet_{u \in \Xi} S_u N$.  Pick an ultrafilter on $\Xi$ including all sets of the form $\{u: x \in u \}$.   Write $a=s_u n_u$ with $n_u \in N, s_u \in S_u$.  Then $s_u \to s$ and $n_u \to n$ for some $s,n$.  We have $s \in \meet_x S_x = K$ and $n \in N$.) 

  Now $L'$ is a compact Lie group, so it has no infinite descending
sequences of closed subgroups.  Thus for  some finite tuple $u=(x_1,\ldots,x_k)$
we have $S_u N = N$.  It follows that $\tau' = L'x$ is a regular orbit in $\tau^k /N$.     The uniformity statement  follows from 
\lemref{unif2}, as in the proof of \thmref{stabilizer-full}; since compactness is assumed, we can take $E$ to
be the entire ambient sort.
\eprf
 
 An earlier version of this section attempted an application to SOP theories, but in this   Krzysztof Krupinski found a gap.

\>{section}

 \def\brG{\breve{G}}
\def\brH{\breve{H}}

\<{section}{Groups with large approximate subgroups}

  In this section we aim to prove:

\<{thm} \lbl{gromov1}
Let $G_0$ be a finitely generated group, $k \in \Nn$.  Assume $G_0$ has a cofinal family of $k$-approximate
subgroups (i.e.   any   finite $F_0 \subset G_0$ is contained in one.)  Then $G_0$   is nilpotent-by-finite.
\>{thm}

 This generalizes Gromov's theorem \cite{gromov}, asserting the same conclusion if $G_0$ has polynomial growth.
 There is by now a small family of proofs of Gromov's theorem and extensions, descending from either Gromov's original proof
 or Kleiner's; the first may have been \cite{vddries-wilkie}, and the most recent,  \cite{shalom-tao}.     I believe all view the group as a metric space, via the Cayley graph,  and analyze it either geometrically or analytically.
 
We will consider an arbitrary sequence of approximate subgroups, rather than balls in the Cayley graph.    A Lie group $L$ lies at the heart of the proof, as in the case of Gromov's.     While Gromov's group arises is the automorphism group  of  the Cayley graph "viewed from afar", we find $L$ and a homomorphism $h: G_0 \to L$ using   the model theoretic/measure-theoretic construction \thmref{stabilizer}, which has no metric aspect.    

Beyond this point,  our proof will adhere  very closely to the outline of  Gromov's.     If the homomorphism into $L$ is trivial, we conjugate it to a nontrivial one in exactly  the   way taken by Gromov, 
 succeeding unless $G_0$ is already virtually abelian \footnote{We say a group $G_0$ is {\em virtually $P$} if some finite index subgroup is $P$.}
  (in which case we are already done).    We now use the earlier \thmref{lie-solvable} covering the linear case
 to show that the image is essentially solvable, and hence a nontrivial homomorphism into an Abelian group can be obtained.  Gromov used the  Tits alternative at the parallel point.      We show that the kernel satisfies the same assumptions as $G_0$; here we make some  further use of Lie theory.  
Induction is carried out on the Lie dimension,  rather than the growth rate exponent which is not available to us; we conclude
  that the group is polycyclic-by-finite, and in particular virtually solvable.    To pass from the polycyclic solvable to the nilpotent case, we quote Tao   \cite{tao-solvable} or Breuillard-Green \cite{bg}  where Gromov cited Milnor-Wolf.

 We will see along the way that $G_0$ is polycyclic-by-finite with $d$ infinite cyclic factors, where $d$ is the dimension of the associated Lie group.   

An alternative statement is that when $G_0$ is {\em not} nilpotent-by-finite, then for some finite $F_0 \subset G_0$, $G_0$ has {\em no} $k$-approximate subgroups containing $F_0$.  If one wishes to seriously use the ambient group $G_0$, some hypothesis on containing sets of generators is necessary (e.g. since any countable family of finitely generated groups embeds jointly in a  single one.)    

The strongest possible general conjecture on the structure of $k$-approximate subgroups would be this:  for some $k',k''$, any $k$-approximate subgroup of a group $G$ is
$k'$-commensurable with one induced by a map into a $k''$-nilpotent group.  Here we say that $X$ is induced by $h$ if $h$ is a homomorphism on some subgroup
$A$ of $G$ into a group $N$, and $X = h \inv h(X)$.  Statements in this vein, possibly restricted to approximate subgroups of a fixed group, have been suggested by Helffgott, E. Lindenstrauss, Breuillard and Tao.    

The conjugation method used in the present section would be powerless in the following scenario:  $X_n$ is a $k$-approximate subgroup of the alternating group $A_n$,
{\em and $X_n$ is conjugation-invariant.}      

Towards the proof of \thmref{gromov1},  we will study the following situation $\diamond$:

\<{itemize}
\item    A language with two sorts  $G,\Phi$;  $G$ carries a group structure; a relation on $G \times \Phi$ defines a family of definable subsets of $G$, $(X_c: c \in \Phi)$.  Additional structure is allowed.  
\item  $M^*$ is a saturated structure, with an elementary submodel $M$.  
\item   $G_0=G(M)$ is  finitely generated.  
\item $X=X_{c^*}$ is a  $c^*$-definable subset $X$ with $G(M) \subset X$  ($c^*$ is an element of $\Phi(M^*)$.)
\item  For all $c \in \Phi(M)$, $X_c$ is finite.
\item  There is an $\inft$-definable subgroup $\G$ of $G$, and a $\bigvee$-definable
subgroup $\tG$,  with $\G \subseteq X \subseteq \tG$, and $\tG/\G$ bounded. 
 $\G,\tG$ are defined over some small subset of $M^*$.  
\item Any subgroup of $G(M)$ has the form $S(M)$ for some 0-definable subgroup $S$ of $G$.
\>{itemize}

In this situation, note:  \<{enumerate} \label{notes}
\item We may replace $\tG$ by the group generated by $X$, without disturbing the hypotheses.
\item Let $G'$ be a 0-definable  subgroup of $G$; $X'_c = X_c \meet G'$; $\G'=\G \meet G', \tG'=\tG \meet G'$.
Then ($\diamond$) holds of the new data, except possibly for the finite generation of $G'(M)$.  When $G'$ has finite index in $G$, this too holds.
\item  There exists a $\bigvee$-definable $\brG \leq \tG$ and a normal $\inft$-definable subgroup $\G'$ of $\brG$
containing $\brG \meet \G$, such that $\brG / \G'$ is a connected Lie group.   (This is \thmref{ap1}.  We have $\brG \meet \G \subseteq \G'$
since the image of $\brG \meet \G$ in $\brG/ \G'$ is a compact normal subgroup, hence trivial.)
\item $X$ is contained in finitely many
cosets of $\brG$ (the image of $X$ modulo $\brG$ is a compact subset of the discrete space $\tG / \brG$.)
\item  $G_0 \meet \brG$ has finite index in $G_0$ (since $G_0 \subseteq X$,   by (4), $G_0$ is contained in finitely many cosets of $\brG$,
equivalently of $G_0 \meet \brG$.) 
\item   
Let  $H_0= G_0 \meet \brG$; let $H$ be a 0-definable group, with $H_0=H(M)$.
So $H$ has finite index in $G$.  Let $\widetilde{H} = \tG \meet H$, let $Y$ be a definable subset of $\widetilde{H}$ commensurable
with $X$, and containing $X \meet \widetilde{H}$, with corresponding family $(Y_c: c \in \Phi')$.  We choose $\Phi' \subset \Phi$ so that
$Y_c$ is commensurable to $X_c$ for $c \in \Phi'$: in particular, $Y_c$ is finite for $c \in \Phi'(M)$.    So $H(M) \subset Y$.
  Now the hypotheses $\diamond$ hold of $(H,Y,\widetilde{H},\G \meet H)$.  Let $\breve{H} = \breve{G} \meet H$,
and $\G'' = \G' \meet H$.  Then $\breve{H}/\G''$ has finite index in $\breve{G} / \G'$, but the latter is connected so they are equal.  Hence $\breve{H}/ \G''$ is
connected; and  $H_0 =  G_0 \meet \brG \subseteq \brG \meet H = \breve{H}$.
Now we are in the same situation $\diamond$, but have in addition $H_0 \leq \breve{H}$.    
\>{enumerate}

Before entering the proof proper, 
we can clarify the meaning of this setup by looking at the Lie rank zero case.

\<{lem}\lbl{fg4}  Assume $\diamond$, and further assume that $\tG/\G$ is totally disconnected.  Then $G=G_0$ is finite.  \>{lem}
 
\prf     Being totally disconnected, $\tG/\G$  contains a compact open subgroup $C$.  
  If $\psi: \brG \to L$ is the canonical map, then $H=\psi \inv(C)$ is a definable group by compactness of $C$, and is commensurable with $X$
  by openness.  
Since $X$ contains $G_0=G(M)$, $H$ is covered by finitely many cosets of $G_0$, so $H \meet G_0$ has finite index in $G_0$.   In particular
it is finitely generated.  Let $F_1$ be a finite set of generators for $H \meet G_0$.
Since $M \prec M^*$, there exists a definable group $H_c$ containing $F_1$ and commensurable with $X_c$, for some $c \in \Phi(M)$; so $H_c$ is finite.  It follows that the group generated by $F_1$ is finite, i.e. $H \meet G_0$ is finite; and thus $G_0$ is finite.  
\eprf

We will need some lemmas on finite generation.  First, if $E$ is a finitely generated group, $N$ a normal subgroup with $E/N$ 
finitely presented, then $N$ is finitely generated as a normal subgroup.  (In particular when $E/N$ is finite, this implies the finite generation of $N$,
a well-known statement used above.)     This in fact  valid  for any equational class:
{\em If $E$ is finitely generated and $N$ is a congruence with $E/N$ finitely presented,
then $N$ is finitely generated as a congruence}.   Indeed let $F$ be a finitely generated free algebra in this equational class, and $h: F \to E$ a surjective homomorphism.  Let $g: F \to E/N$ be the composition $F \to E \to E/N$.  Since $E/N$ is finitely presented, $g$ has a finitely generated kernel $K$.
Thus $N=h(K)$ is finitely generated.    
 
A $\bigvee$-definable subgroup is called {\em definably generated} if it is generated by
a definable subset.  If $G$ is a topological group, let $G^0$ denote the connected component of $1$; it is a closed normal subgroup of $G$.

  Let $H$ be a sufficiently saturated group (with possible additional structure),
 $\brH$ be a $\bigvee$-definable subgroup, $\G$ a $\bigwedge$-definable subgroup, with $\G \normal \brH$.    
 Let $\pi: \brH \to \brH/\G$ be the quotient map.  Recall  the logic topology on  $\brH/\G$  from \S 4.  In particular, a subset $Z$ of
 of the quotient is compact iff $\pi \inv (Z)$ is contained in a definable set.

\<{lem}\lbl{gr1}  Let $H,\brH,\G,\pi$ be as above, and assume $A=\brH/\G$ is locally compact.   Assume $A/A^0$  is finitely generated. 
Then   $\brH$ is definably generated.  \>{lem}

\prf        Let $U$ be a compact   neighborhood of $1$ in $A$.  Then $\pi \inv (U)$
is contained in a definable subset $D$ of $\brH$.  $U$ generates an open subgroup of $A$; this open subgroup is also closed, and must contain $A^0$.  On the other hand $A /A^0$ is generated by finitely many elements $\pi(h_1),\ldots,\pi(h_r)$.
Let $D' = D \union \{h_1,\ldots,h_r \}$.  
Since $D$ contains $\ker \pi$  and $\pi  \brH$ is generated by $\pi (D')$, it follows that $\brH$ is generated by $D'$.  

\eprf 

 In fact if $A/A^0$ is $m$-generated, the proof shows that $\brH$ is generated by $Y$ along with $m$ additional elements, whenever $Y$ is a definable set containing $\G$.

\<{lem}\lbl{gr1.1}  Let $G$ be a sufficiently saturated group (with additional structure), $H$ a definable normal subgroup with $G/H$ Abelian.  Let $\brG$ be a $\bigvee$-definable subgroup of $G$, $\G$ a $\bigwedge$-definable subgroup, with $\G \normal \brG$ and
 $E=\brG / \G$  a connected Lie group. 
Then $H \meet \brG$ is definably generated.  \>{lem}

\prf   Let $\brH =  H \meet \brG$, and $\pi: \brG \to E$ be the canonical map.  By \lemref{towers},    $\pi|\brH$ induces an isomorphism of topological groups
$ \brH / (\brH \meet \G) \cong \pi(\brH)$.  

 Since $H$ is normal in $G$, $\brH$ is normal in $\brG$, so $\pi(\brH)$ is normal in $E$, and hence so is $\pi(\brH)^0$,  
The commutator subgroup $[E,E]$ is contained in $\pi(\brH)^0$, so the quotient 
$E / \pi(\brH)^0$ is isomorphic to $\Rr^n \times \Rr^m/\Zz^m$.     The  image of
 $\pi(\brH)$ in $E/ \pi(\brH)^0$ contains no nontrivial connected groups, so it is discrete.  
  Now it is well-known that a discrete subgroup of $\Rr^n \times \Rr^m/\Zz^m$ is finitely generated, indeed admits a generating 
  set with  at most  $n+m$ elements.
    By \lemref{gr1}, $H \meet \brG$ is  definably generated.     
    \eprf

The next lemma will play an essential role in the proof, allowing the key   \lemref{fg2} to be propagated.   
  It continues to hold if $G/N$ is assumed to be nilpotent, rather than Abelian; indeed  it suffices to find a sequence of definable subgroups
 $G=H_1 \supset \cdots \supset H_k = H$ with $H_{i+1}$ normal in $H_i$, and $H_i/H_{i+1}$ Abelian, and apply the lemma inductively.

\<{lem}\lbl{fg1}  Assume $\diamond$ holds, and let $N$ be a 0-definable normal subgroup of $G$.

Then $\diamond$ holds if $G$ is replaced by $G/N$, and $X,\tG,\G$ by their images in $G/N$.

If $G/N$ is finite or Abelian,  then $\diamond$ holds if $G,X, \tG,\G$ are replaced
by $N,X\meet N, \tG \meet N, \G \meet N$.  
 \>{lem}

\prf  The first statement is straightforward; so is the second, except for the finite
generation of $N_0=N(M)$.  We proceed to show this.  In case $G/N$ is finite, so is $G_0/N_0$ since $M$ is an elementary submodel, so
finite generation is clear.  Assume therefore that $G/N$ is Abelian.


   
 
 Let $G_1= G_0 \meet \brG$.    As $G_1$ has finite index in $G_0$ by (5),  
  it is a finitely generated group.
 
 Let $N_1 = N  \meet G_1$.  Since $G_1/(N \meet G_1)$ is finite or Abelian, it is a finitely presented group.  By the remarks above, 
    $N_1$ is finitely generated as a normal subgroup of $G_1$.

 Let $g_1,\ldots,g_r$ be generators for $G_1$, and let $T_i(x)=g_i \inv x g_i$.   Then $N_1$ is finitely generated as a group with
 these operators.    Let $Y$ be a finite subset of $N_1$ such that $N_1$ is generated by $Y$ under multiplication and the operators $T_i$.  
 
By \lemref{gr1.1}, $\brG \meet N$ is   generated by an $M^*$- definable set $U$.  We may take $Y \subset U=U \inv$.   Since $\brG$ and $N$ are closed under
the operators $T_i$ (as $g_i \in \brG$ and $N$ is normal), we have $T_i(U) \subset U \cdot \cdots \cdot U = U^m$ for some $m$.  
Since $U$ is a definable subset of $\brG$, it is contained in finitely many translates of $X$.  
Now $M$ is an elementary submodel of $M^*$.  So there exists an $M$-definable set $U' \subset N$ containing $Y$, with $T_i(U') \subset U' \cdot \cdots \cdot U'$, and $U'$ contained in finitely many translates of some $X_c, c \in \Phi(M_0)$.  From the last property it follows that $U'$ is finite; so $U' \subset M$; hence $U' \subset G_0$.  Thus $U' \subset N_1$.  
Moreover the group generated by $U'$ is closed under the operators $T_i$, and contains $Y$.  So it equals $N_1$.  This shows that $N_1$ is a finitely generated group.  Since it has finite index in $N_0 = N \meet G_0$, it follows that $N_0$ too is a  finitely generated group.
\eprf

\<{lem}\lbl{fg2}  Assume $\diamond$ holds, and $G_0$ is infinite.  Then there exists a  normal subgroup
$N_0$ of $G_0$ with $G_0/N_0$ virtually Abelian, and infinite. 
 \>{lem}

\prf       Let $G_0'=G_0 \meet \brG$.   By note (5) above, $G_0'$ has finite index in $G_0$.  
 If $N_0'$ is normal in $G_0'$
with infinite virtually Abelian quotient, let $N_0$ be the intersection of the finitely many $G_0$-conjugates of $N_0'$;
then $G_0/N_0$ is infinite and virtually Abelian.   Thus proving the lemma for $G_0'$ would imply it for $G_0$.  By note (6), the hypotheses hold of
$G_0'$; so we may assume $G_0 \leq \brG$.   

Let $L= \brG / \G'$ as in note (3),  $d=\dim(L)$, and consider  the natural homomorphism $\psi: \brG \to L$.
Note that $\psi(G_0)$ has a cofinal system of $k$-approximate subgroups.  By \corref{lie-solvable}, any finite subset $w$ of $\psi(G_0)$ is contained in at most $k''$ cosets of a $d+2$-solvable subgroup $S_w$ of $L$.   Taking an ultraproduct, $L$ embeds in an ultraproduct  of itself, in such a way that   the image of $\psi(G_0)$ is contained  in at most 
$k''$ cosets of a $d+2$-solvable group $S$.   Thus $\psi(G_0)$ has a solvable subgroup $S'$ of finite index.   If $S'$ is infinite,
then it contains  a subgroup $S''$ of finite index, such that $S''/ [S'',S'']$ is infinite.  Thus $\psi \inv(S'')$ is a  finite index subgroup   of $G_0$,
and $\psi \inv ([S'',S''])$ is a normal subgroup with infinite Abelian quotient.
So we are done unless $S'$ above is finite, i.e. $\psi(G_0)$ is finite, so that a finite index subgroup $H_0$ of $G_0$ is contained in $\G'$.
We have $H_0=H(M)$ for some 0-definable subgroup $H$ of $G$.

For $g \in G$, let $ad_g(x) = g \inv x g$, and let $\tau_g = ad_g | H_0$.    Let $J = \{g \in G:  \tau_g (H_0) \leq \brG \}$.  
If $g \in J$, we may repeat the previous paragraph with $\psi \circ \tau_g$ in place of $\psi$.  Thus again we are done unless $\psi \circ \tau_g(H_0)$
is finite for any $g \in J$.   We thus assume this is the state of affairs.   

The rest of the proof is a straightforward transcription of the corresponding part of \cite{gromov}.  
By Jordan's theorem \cite{jordan}, since  $\psi \circ \tau_g (H_0)$ is a finite subgroup of the Lie group $L$, it has an Abelian subgroup $S_g$ of index $\leq \mu$,
with $\mu$ independent of $g$.   If $\psi \circ \tau_g$ can have arbitrarily large finite size for $g \in J$, taking an ultraproduct, we obtain
a homomorphism to a group with an infinite Abelian subgroup of index $\leq \mu$.  Thus in this case too the lemma is proved, and
we may assume $\psi \circ \tau_g (H_0)$ has   size $\leq \mu'$ for some fixed $\mu'$.

Let $F_1$ be a finite set of generators for $H_0$. 
Let $U$ be a neighborhood of the identity in the Lie group $L$, such that if $u \in U$ is an element of order $\leq \mu'$, then $u=1$.  (For instance
we can take a neighborhood $V$ of the Lie algebra on which the exponential map is injective, and then let $U = exp((1/\mu') V)$.)  
Since $\{1\}$ is closed and $U$ is open, there exists a definable set $D_2 \subset G$ with $\G' \subset D_2 \subset \psi \inv(U)$.  
Since $F_1 \inv \G' F_1 = \G' \subset D_2$, we can find a definable set $D$ with $\G' \subset D$ and $F_1 \inv D F_1  \subset D_2$.     Any subgroup of $D_2$ of size $\leq \mu'$  is trivial.  Now if $\tau_g(F_1) \subset D_2$, then $\tau_g(H_0) \leq \brG$, so $g \in J$; hence $\psi \circ \tau_g (H_0)$ has size $\leq \mu'$;
but $\psi \circ \tau_g (F_1)$ is a set of elements of $U$, and any such nonidentity element has order $> \mu'$; so $\psi \circ \tau_g(F_1)$
must reduce to the identity element of $L$.  Hence  if $\tau_g(F_1) \subset D_2$, then $\tau_g(H_0) \subset \G'$, and in particular $\tau_g(F_1) \subset D$.

Let $W=\{g:  \tau_g(F_1) \subset D\} = \{g:  \tau_g(F_1) \subset D_2\}$.   If $g \in W$ and $f \in F_1$, then $gf \in W$, since $(gf) \inv F_1 gf \subseteq f \inv D f \subseteq D_2$.
So $W$ is a definable, right $F_1$-invariant set.  Now in the model $M$, any definable, right $F_1$-invariant set is empty or contains $H_0$.
Since $M \prec M^*$, it follows that $W=\emptyset$ or $W$ contains $H$.   We have $1 \in W$, as $H_0 \leq \Gamma' \leq D$.
So all $H$-conjugates of $F_1$ are contained in $D$.  Note
that $D$ is contained in the union of finitely many translates of $X$.  It follows that all $H_0$-conjugates of $F_1$ are contained
in finitely many translates of some $X_c$, $c \in \Phi(M)$.     In this case each element of $F_1$ has centralizer of finite index in $H_0$;
so $H_0$ has a center of finite index; we may take $G_0'$ to be this center, and $N=1$.
\eprf

We will need some elementary group-theoretic discussion before proceeding.   We define a group $H$ to be $0$-polycyclic if it is trivial, and to be   $d+1$-polycyclic if it has a $d$-polycyclic normal subgroup $N$, with $H/N$ a finitely generated Abelian group of rank $1$.     In particular, $H$ is $d$-solvable.

   For $d \geq 0$, say a finitely generated group $H$ is almost $d$-polycyclic if  it has   subgroups $H_{2d+2} \normal H_{2d+1} \normal H_{2d}  \normal \cdots \normal H_1= H$    with $H_{2i}/H_{2i+1}$ finite ($0 \leq i \leq d$), and  $H_{2i+1}/H_{2i+2} \cong \Zz$ ($0 \leq i \leq d$), and $H_{2d+2}=0$.    By reverse induction on $i$ we see that each $H_i$ is finitely generated in this situation.

These definitions differ in that the quotients   $H_{2i}/H_{2i+1}$ are not required to be Abelian, but 
if $H$ is  almost $d$-polycyclic, then it does have a  $d$-polycyclic normal subgroup of finite index.   To show this we may pass to a finite index subgroup, so we may assume $H$ has an almost $d-1$- polycyclic normal subgroup $N$  with 
 $H/N \cong \Zz$.     Like all almost polycyclic groups, $N$ is  finitely generated.   
 Using the induction hypothesis, let $N_1$ be a $d-1$-polycyclic subgroup of $N$, with $[N:N_1]=r< \infty$.  
As $N$ is   finitely generated, it  has only finitely many subgroups of index $r$.  Let $N_2$ be their intersection.  Then $N_2$ is $d-1$-polycyclic  and is characteristic in $N$, hence normal in $H$.  $H/N_2$ contains the finite group
$N/N_2$ as a normal subgroup; within $H/N_2$, the centralizer of $N/N_2$ has the form $H_1/N_2$, with $H_1$ a finite index subgroup of $H$.
Now $H_1/(N \meet H_1) \cong \Zz$, while $(N \meet H_1) / N_2$ is a finite central subgroup of $H_1/N_2$; so $H_1/N_2$ is a finitely generated Abelian group of rank $1$.   Thus $H_1$  is $d$-polycyclic.

\<{lem} \lbl{pbf} Assume $\diamond$.  Then $G_0=G(M)$ is polycyclic-by-finite (and in particular solvable-by-finite).  \>{lem}

\prf  We use induction on $d=\dim(L)$, $L=\brG / \G'$.    If $d=0$, then $\tG/\G$ is totally disconnected, hence $G_0$ is finite by \lemref{fg4}. 
For higher $d$, we use \lemref{fg2}.  By note (2) to $\diamond$ we may pass to a finite index subgroup; so we may assume there
exists a 0-definable normal subgroup $N$ of $G$, with $G/N$ infinite Abelian.  By \lemref{fg1}, the hypotheses $\diamond$
hold for $N, X \meet N, \tG \meet N, \G \meet N$, and also for the images in $G/N$.  By \lemref{fg4} applied to $G/N$, we see that the image of $\tG$ in $G/N$ has Lie rank $\geq 1$.   By \remref{towers} (3) it follows that $\tG \meet N$  has Lie rank $<d$.  So the inductive hypothesis applies, and  $N_0=N(M)$ is polycyclic-by-finite.    Thus $G_0$ is almost polycyclic and hence also polycyclic by finite .  

\eprf

We are now essentially in the solvable case, and can quote either \cite{tao-solvable}, or \cite{bg}.    Polycyclicity is a strong additional tool,
and with it one may be able to reduce to the sum-product phenomenon for fields somewhat more rapidly;
a model-theorist is reminded here of Zilber's arguments in the 70's, connecting solvable groups of finite Morley rank with definable fields.  
 We will simply invoke \cite{bg}; thanks to Emmanuel Breuillard for pointing out a nicer path to this paper than
we had initially.  

\<{lem}\lbl{nbf} Assume $\diamond$.  Then $G_0=G(M)$ is nilpotent-by-finite.   \>{lem}

\prf    By \lemref{fg2}, $G_0$ has a normal subgroup $N_0$ with $G_0/N_0$ virtually Abelian; as we just saw, 
  $N_0$ satisfies $\diamond$, has lower Lie rank, and so inductively is nilpotent by finite.    Passing to a finite index subgroup, we may assume $G_0$ is solvable, as well as polycyclic.    Now any polycyclic group is linear over $\Cc$;
  see chapter 4 of \cite{raghunathan} for a stronger result, due to Auslander-Swan.  Hence $G_0$ can be viewed as a solvable subgroup of 
  $GL_n(\Cc)$; according to \cite{bg},  every $k$-approximate
subgroup of  $G_0$ is covered by a bounded number of cosets of a $(n-1)$-step nilpotent subgroup  of
$GL_n(C)$.   Taking the ultraproduct of the cofinal family of approximate groups associated with $G_0$,  we see that  $X$ and hence $G_0$ are covered by finitely
many cosets of an $(n-1)$-step nilpotent group, hence $G_0$ is itself virtually nilpotent.    \eprf

\prf of \thmref{gromov1} 

Let $(X_c: c \in \Phi_0)$ be the given family of $k$-approximate subgroups of $G_0$.  Consider the two-sorted structure $(G,\Phi_0,\cdot,E)$
where $(x,c) \in E$ if $x \in X_c$.  Enrich it by adding a predicate for each subgroup of $G_0$.  Further enrich the language by closing under probability quantifiers as in \secref{measurequantifiers}.
Let $M$ be the resulting structure, and let $M^*$ be a saturated elementary extension.  By saturation and by the cofinality of the $X_c$, there exists $c^* \in \Phi(M^*)$
with $G_0 \subset X_{c^*}$.  All clauses of $\diamond$ are now clear, so by \lemref{nbf}, $G_0$ is nilpotent-by-finite.
 \eprf

\>{section}

\<{thebibliography}{10}

\bibitem{BU}  Itai Ben Yaacov and Alexander Usvyatsov, Continuous first order logic and local stability.  To appear in the Transactions of the American Mathematical Society. 

\bibitem{bergman-lenstra}  Bergman, George M.; Lenstra, Hendrik W., Jr.,
Subgroups close to normal subgroups.
J. Algebra 127 (1989), no. 1, 80--97.

\bibitem{bgt} {Emmanuel Breuillard, Ben Green, Terence Tao, Approximate subgroups of linear groups,
arXiv:1005.1881} 

\bibitem{bg} Emmanuel Breuillard, Ben Green, Approximate groups, II: the solvable linear case, arXiv:0907.0927  

 \bibitem{bkt} Bourgain, J., Katz, N.H., Tao, T.C.: A sum-product estimate in
 finite Þelds and applications, Geom. Funct. Anal. 14 (2004), no. 1, 27-57. 

\bibitem{CK} Chang, C. C.; Keisler, H. J. Model theory. Third edition. Studies in Logic and the Foundations of Mathematics, 73. North-Holland Publishing Co., Amsterdam, 1990.

\bibitem{chang} 
Chang, Mei-Chu, Product theorems in ${\rm SL}\sb 2$ and ${\rm SL}\sb 3$.  J. Inst. Math. Jussieu  7  (2008),  no. 1, 1--25.

\bibitem{CH}   Cherlin, Gregory; Hrushovski, Ehud Finite structures with few types. Annals of Mathematics Studies, 152. Princeton University Press, Princeton, NJ, 2003. vi+193

\bibitem{curtis-reiner}  Curtis, Charles W.; Reiner, Irving Representation theory of finite groups and associative algebras. 
Reprint of the 1962 original. AMS Chelsea Publishing, Providence, RI, 2006. xiv+689 pp.

\bibitem{dinai}  Dinai, forthcoming.

\bibitem{eg}    Elekes, Gy\"orgy ; Kir‡ly, Zolt‡n 
On the combinatorics of projective mappings. (English summary)
J. Algebraic Combin. 14 (2001), no. 3, 183--197.


\bibitem{vddries} van den Dries, L, notes in {\tt http://www.math.uiuc.edu/$\sim$vddries/}.

\bibitem{vddries-wilkie} 
 van den Dries, L.; Wilkie, A. J.
Gromov's theorem on groups of polynomial growth and elementary logic. 
J. Algebra 89 (1984), no. 2, 349--374.

\bibitem{erdos-sz}  Erd\"os, P.; Szemer\'edi, E. On sums and products of integers.  Studies in pure mathematics,  213 - 218, BirkhŠuser, Basel, 1983.

\bibitem{gleason} Gleason, A. M.
The structure of locally compact groups. 
Duke Math. J. 18, (1951). 85--104

\bibitem{goldbring} Goldbring, I. Hilbert's Fifth Problem for Local Groups, to appear in the Annals of Mathematics, available from: http://www.math.uiuc.edu/$\sim$igoldbr2/.

\bibitem{gs}    Ben Green,  Tom Sanders,
A quantitative version of the 
idempotent theorem in harmonic analysis, 
Annals of Mathematics, 168 (2008), 1025Ð1054 
  

\bibitem{gromov}  Gromov, Mikhael, Groups of polynomial growth and expanding maps. Inst. Hautes \'Etudes Sci. Publ. Math. No. 53 (1981), 53--73.

\bibitem{halmos}   P. Halmos, ``Measure Theory,'' Van Nostrand, Princeton, NJ, 1950.

\bibitem{helfgott}  Helfgott, H. A. Growth and generation in ${\rm SL}\sb 2(\mathbb Z/p\mathbb Z)$.  Ann. of Math. (2)  167  (2008),  no. 2, 601--623. 
 
 \bibitem{helfgott3} H. A. Helfgott,  Growth in $SL_3(Z/pZ)$,  	arXiv:0807.2027
 
\bibitem{hewitt-savage} Hewitt, E., Savage, L.J.: Symmetric measures on Cartesian products. Trans. Amer. Math. Soc.,
80, 470-501 (1955)

\bibitem{PAC} Hrushovski, Ehud 
Pseudo-finite fields and related structures. Model theory and applications, 151--212,
Quad. Mat., 11, Aracne, Rome, 2002.

\bibitem{hw} Hrushovski, Ehud; Wagner, Frank,
Counting and dimensions.  Model theory with applications to algebra and analysis. Vol. 2, 161--176,
London Math. Soc. Lecture Note Ser., 350, Cambridge Univ. Press, Cambridge, 2008. 

\bibitem{hpp}  Hrushovski, Ehud; Peterzil, Ya'acov; Pillay, Anand Groups, measures, and the NIP. J. Amer. Math. Soc. 21 (2008), no. 2, 563--596 (section 7.)

\bibitem{nip2}  Hrushovski, Ehud;  Pillay, Anand, On NIP and invariant measures,   	arXiv:0710.2330 

\bibitem{nip3}  Ehud Hrushovski, Anand Pillay, Pierre Simon, Generically stable and smooth measures in NIP theories, 	arXiv:1002.4763

 \bibitem{kallenberg}
 Kallenberg, Olav On the representation theorem for exchangeable arrays. J. Multivariate Anal. 30 (1989), no. 1, 137--154.

\bibitem{jordan} Jordan, Camille, M\'emoire sur les \'equations diff\'erentielles lin\'eaires ˆ int\'egrale algbriques, Crelle 84 (1878) pp. 89-215

\bibitem{kaplansky} Kaplansky, Irving Lie algebras and locally compact groups. Reprint of the 1974 edition. Chicago Lectures in Mathematics. University of Chicago Press, Chicago, IL, 1995. xii+148 pp.

\bibitem{kim-pillay}  Kim, Byunghan; Pillay, Anand Simple theories. Joint AILA-KGS Model Theory Meeting (Florence, 1995).  Ann. Pure Appl. Logic  88  (1997),  no. 2-3, 149--164.

\bibitem{komlos-s} Koml\'os, J.; Simonovits, M. Szemer\'edi's regularity lemma and its applications in graph theory.  Combinatorics, Paul Erd\"os is eighty, Vol. 2 (Keszthely, 1993),  295--352, Bolyai Soc. Math. Stud., 2, J‡nos Bolyai Math. Soc., Budapest, 1996.

\bibitem{krauss}  Peter H. Krauss, Represetation of Symmetric Probability Models, JSL v. 34 No. 2 (June 1969), pp. 183-193

\bibitem{krauss-scott} D. Scott and P. Krauss, Assigning probabilities to logical formulas, in:  Aspects of inductive logic,
ed. J. Hintikka and P. Suppes, North Holland, Amsterday, 1966, pp. 219-259.

\bibitem{larsen-pink}  Michael J. Larsen and Richard Pink. Finite subgroups of algebraic groups, 1998.  Available from \text{http://www.math.ethz.ch/$\sim$pink/publications.html}.

\bibitem{lee} Lee, John M.,  Riemannian manifolds: An introduction to curvature. 
Graduate Texts in Mathematics, 176. 
Springer-Verlag, New York, 1997.

\bibitem{marker}  Marker, David Model theory. An introduction. Graduate Texts in Mathematics, 217. Springer-Verlag, New York, 2002

\bibitem{marker90}      David Marker,     Semialgebraic Expansions of C,  Transactions of the American Mathematical Society, Vol. 320, No. 2 (Aug., 1990), pp. 581-592                            

\bibitem{morley} Morley, Michael The L\"owenheim-Skolem theorem for models with standard part.  1971  Symposia Mathematica, Vol. V (INDAM, Rome, 1969/70)  pp. 43--52 Academic Press, London 


\bibitem{n-p}  Nikolay Nikolov, L\'aszl\'o Pyber 
 Product decompositions of
quasirandom groups,   arXiv:math/0703343  
 
\bibitem{pillay} Anand Pillay, {\tt http://www.amsta.leeds.ac.uk/$\sim$pillay/}
\bibitem{pillay-book}
Pillay, Anand,
An introduction to stability theory.
Oxford Logic Guides, 8. The Clarendon Press, Oxford University Press, New York, 1983. xii+146 pp

\bibitem{poizat} Bruno Poizat , Cours de th\'eorie des mod\`eles.  Bruno Poizat, Lyon, 1985. vi+584 pp. 
	
\bibitem{poizat2} Poizat, Bruno An introduction to algebraically closed fields and varieties.  in:  The model theory of groups (Notre Dame, IN, 1985Ð1987), 41Ð67, Notre Dame Math. Lectures, 11, Univ. Notre Dame Press, Notre Dame, IN, 1989

\bibitem{raghunathan} Raghunathan, M. S. Discrete subgroups of Lie groups.  Math. Student  2007,  Special Centenary Volume, 59--70 (2008). 
 
\bibitem{shalom-tao} Yehuda Shalom, Terence Tao, A finitary version of Gromov's polynomial growth theorem,  arXiv:0910.4148v2 [math.GR] 
 
\bibitem{lazy} Shelah, Saharon
The lazy model-theoretician's guide to stability.
Comptes Rendus de la Semaine d'ƒtude en Th\'eorie des Modles (Inst. Math., Univ. Catholique Louvain, Louvain-la-Neuve, 1975).
Logique et Analyse (N.S.) 18 (1975), no. 71-72, 241--308.

\bibitem{shelah} Shelah, S. Classification theory and the number of nonisomorphic models. Second edition. Studies in Logic and the Foundations of Mathematics, 92. North-Holland Publishing Co., Amsterdam, 1990. xxxiv+705 pp.

\bibitem{shelah-simple} Shelah, Saharon
Simple unstable theories.
Ann. Math. Logic 19 (1980), no. 3, 177--203.  


\bibitem{tao-nc}  Terence Tao, Product-set estimates for non-commutative groups, arXiv   	math/0601431 

\bibitem{tao-rings} Terence Tao,   
     The sum-product phenomenon in arbitrary rings, arXiv:0806.2497
     
\bibitem{tao-blog}  Terence Tao,  \hbox{ http://terrytao.wordpress.com/2007/03/02/open-question-noncommutative-freiman-theorem/   }

\bibitem{tao-solvable}  Terence Tao, Freiman's theorem for solvable groups,  arXiv:0906.3535  

\bibitem{t-v}   Tao, Terence; Vu, Van Additive combinatorics. Cambridge Studies in Advanced Mathematics, 105. Cambridge University Press, Cambridge,  2006. xviii+512 pp.
     
\bibitem{yamabe} Yamabe, Hidehiko
A generalization of a theorem of Gleason. 
Ann. of Math. (2) 58, (1953). 351--365.

 \bibitem{weil} Weil, Andr\'e, Foundations of algebraic geometry. American Mathematical Society, Providence, R.I. 1962 xx+363 pp.  
%
   
\bibitem{zimmer} Zimmer, Robert J., Essential results of functional analysis.
Chicago Lectures in Mathematics. University of Chicago Press, Chicago, IL, 1990. x+157
   
\>{thebibliography}

\end{document}